\documentclass[UTF-8,reqno]{amsart}
\usepackage{enumerate, bbm}
\usepackage{amssymb,url,color, booktabs}
\usepackage{tikz}
\usepackage{mathrsfs}
\usepackage{dutchcal}
\usepackage{threeparttable}
\usepackage{color}
\usepackage[colorlinks=true]{hyperref}
\hypersetup{
    linkcolor=blue,          % color of internal links
    citecolor=red,        % color of links to bibliography
    filecolor=blue,      % color of file links
    urlcolor=cyan
}

\usepackage{color}
\usepackage{ulem}

 \usepackage{scalerel} 

\definecolor{MyDarkBlue}{cmyk}{0.8,0.3,0.8,0.4}
\definecolor{yellow}{rgb}{0.99,0.99,0.70}
\definecolor{white}{rgb}{1.0,1.0,1.0}
\definecolor{black}{rgb}{0.00,0.00,0.00}

\numberwithin{equation}{section}

\newcommand{\be}{\begin{eqnarray}}
\newcommand{\ee}{\end{eqnarray}}
\newcommand{\ce}{\begin{eqnarray*}}
\newcommand{\de}{\end{eqnarray*}}
\newtheorem{theorem}{Theorem}[section]
\newtheorem{lemma}[theorem]{Lemma}
\newtheorem{remark}[theorem]{Remark}
\newtheorem{definition}[theorem]{Definition}
\newtheorem{proposition}[theorem]{Proposition}
\newtheorem{Examples}[theorem]{Example}
\newtheorem{corollary}[theorem]{Corollary}

\usepackage[nobysame]{amsrefs}
\BibSpec{article}{%
+{}{\PrintAuthors} {author}
+{,}{ \textrm} {title}
+{.}{ \textit} {journal}
+{,}{ \textbf} {volume}
+{}{ \parenthesize} {date}
+{,}{ } {pages}
+{.}{ arXiv:} {eprint}
+{.}{} {transition}
}
\BibSpec{book}{%
+{}{\PrintAuthors} {author}
+{,}{ \textit} {title}
+{.}{ \textrm} {series} %+{,}{ \textrm} {series}
+{,}{ Vol.} {volume} 
+{.}{ } {publisher}
+{,}{ } {date}
+{.}{} {transition}
}

\def\nor{|\mspace{-3mu}|\mspace{-3mu}|}

\def\var{{\mathrm{var}}}

\def\eps{\varepsilon}

\def\e{\mathrm{e}}

\def\p{\partial}

\def\[{{\Big[}}
\def\]{{\Big]}}
\def\<{{\langle}}
\def\>{{\rangle}}
\def\({{\Big(}}
\def\){{\Big)}}

\def\bx{{\mathbf{x}}}

\def\H{{\scaleto{H}{2pt}}}

\def\dif{{\mathord{{\rm d}}}}
\def\dis{{\mathord{{\rm \bf d}}}}

\def\min{{\mathord{{\rm min}}}}

\def\no{\nonumber}
\def\={&\!\!=\!\!&}

\def\bB{{\mathbf B}}
\def\bC{{\mathbf C}}

\def\cA{{\mathcal A}}
\def\cB{{\mathcal B}}
\def\cC{{\mathcal C}}

\def\cE{{\mathcal E}}

\def\cH{{\mathcal H}}

\def\cK{{\mathcal K}}
\def\cL{{\mathcal L}}
\def\cM{{\mathcal M}}

\def\cP{{\mathcal P}}

\def\cS{{\mathcal S}}

\def\cW{{\mathcal W}}

\def\mC{{\mathbb C}}
\def\mD{{\mathbb D}}
\def\mE{{\mathbb E}}

\def\mH{{\mathbb H}}
\def\mI{{\mathbb I}}

\def\mL{{\mathbb L}}

\def\mN{{\mathbb N}}

\def\mP{{\mathbb P}}
\def\mQ{{\mathbb Q}}
\def\mR{{\mathbb R}}

\def\mX{{\mathbb X}}

\def\bB{{\mathbf B}}
\def\bP{{\mathbf P}}
\def\bQ{{\mathbf Q}}

\def\1{{\mathbf{1}}}

\def\sB{{\mathscr B}}

\def\sF{{\mathscr F}}

\def\sI{{\mathscr I}}
\def\sJ{{\mathscr J}}

\def\sL{{\mathscr L}}

\def\sV{{\mathscr V}}
\def\sW{{\mathscr W}}

\def\geq{\geqslant}
\def\leq{\leqslant}
\def\ge{\geqslant}
\def\le{\leqslant}

\def\div{\mathord{{\rm div}}}

\def\var{{\mathrm{var}}}

\def\eps{\varepsilon}

\def\e{\mathrm{e}}

\def\p{\partial}

\def\[{{\Big[}}
\def\]{{\Big]}}
\def\<{{\langle}}
\def\>{{\rangle}}
\def\({{\Big(}}
\def\){{\Big)}}

\def\bx{{\mathbf{x}}}

\def\dif{{\mathord{{\rm d}}}}
\def\dis{{\mathord{{\rm \bf d}}}}

\def\min{{\mathord{{\rm min}}}}

\def\no{\nonumber}
\def\={&\!\!=\!\!&}
\def\bt{\begin{theorem}}
\def\et{\end{theorem}}
\def\bl{\begin{lemma}}
\def\el{\end{lemma}}
\def\br{\begin{remark}}
\def\er{\end{remark}}
\def\bx{\begin{Examples}}
\def\ex{\end{Examples}}
\def\bd{\begin{definition}}
\def\ed{\end{definition}}
\def\bp{\begin{proposition}}
\def\ep{\end{proposition}}
\def\bc{\begin{corollary}}
\def\ec{\end{corollary}}

\def\geq{\geqslant}
\def\leq{\leqslant}
\def\ge{\geqslant}
\def\le{\leqslant}

\def\div{\mathord{{\rm div}}}

\def\bP{{\mathbf P}}

\def\<{\langle} \def\>{\rangle}

\def\bpf{\begin{proof}}
\def\epf{\end{proof}}

\allowdisplaybreaks

\begin{document}
	
\title[DFSDEs driven by (fractional) Bronian motion and NSE]{Distribution-flow dependent SDEs driven by (fractional) Brownian motion and Navier-Stokes equations}
\date{\today}
\author{Zimo Hao, Michael R\"ockner and Xicheng Zhang}

\thanks{{\it Keywords: \rm Distribution-flow dependent SDE, Navier-Stokes equations, fractional Brownian motion}}

\thanks{
This work is supported by National Key R\&D program of China (No. 2023YFA1010103) and NNSFC grant of China (No. 12131019)  and the DFG through the CRC 1283 
``Taming uncertainty and profiting from randomness and low regularity in analysis, stochastics and their applications''. }

\address{Zimo Hao:
Fakult\"at f\"ur Mathematik, Universit\"at Bielefeld,
33615, Bielefeld, Germany\\
Email: zhao@math.uni-bielefeld.de} 

\address{Michael R\"ockner:
Fakult\"at f\"ur Mathematik, Universit\"at Bielefeld,
33615, Bielefeld, Germany\\
Email: roeckner@math.uni-bielefeld.de} 

\address{Xicheng Zhang:
School of Mathematics and Statistics, Beijing Institute of Technology, Beijing 100081, China\\
Faculty of Computational Mathematics and Cybernetics, Shenzhen MSU-BIT University, 518172 Shenzhen, China\\
Email: XichengZhang@gmail.com
 }

\begin{abstract}
Motivated by the probabilistic representation for solutions of the Navier-Stokes equations, we introduce a novel class of stochastic differential equations that depend on the entire flow of its time marginals. We establish the existence and uniqueness of both strong and weak solutions under one-sided Lipschitz conditions and for singular drifts. These newly proposed distribution-flow dependent stochastic differential equations are closely connected to quasilinear backward Kolmogorov equations and Fokker-Planck equations.
Furthermore, we investigate a stochastic version of the 2D-Navier-Stokes equation associated with fractional Brownian noise. 
We demonstrate the global well-posedness and smoothness of solutions when the Hurst parameter $H$ lies in the range $(0, \frac12)$ and the initial vorticity is a finite signed measure.
\end{abstract}

\maketitle
\setcounter{tocdepth}{2}
\tableofcontents

\section{Introduction}
Throughout this paper we fix $T>0$ and $d\in\mN$ and write 
$$
\mD_T:=\{(s,t): 0\leq s<t\leq T\}.
$$ 
Let $\cP:=\cP(\mR^d)$ be the space of all probability measures over $\mR^d$, which is endowed with the weak topology. 
Let $\cC^d_\cP:=C(\mR^d;\cP(\mR^d))$ be
the space of all continuous probability measure-valued functions from $\mR^d$ to $\cP(\mR^d)$.
Let $\{W_{t}\}_{t\in[0,T]}$ be a $d$-dimensional standard Brownian motion on some probability space $(\Omega,\sF,\mP)$.
We consider the following nonlinear stochastic differential equation (SDE), also called  distribution-flow dependent SDE (abbreviated as DFSDE): for $(s,t,x)\in\mD_T\times\mR^d$,
\begin{align}\label{DFSDE}
X_{s,t}^x=x+\int_s^t B(r,X_{s,r}^x,\mu_{r,T}^{\centerdot},\mu_{s,r}^{\centerdot})\dif r+\int_s^t \Sigma(r,X_{s,r}^x,\mu_{r,T}^{\centerdot},\mu_{s,r}^{\centerdot})\dif W_r,
\end{align}
 where $\mu^x_{s,t}=\mP\circ(X_{s,t}^x)^{-1}$ is the probability distribution measure of  $X_{s,t}^x$ {satisfying 
 \begin{align}\label{flow-1125}
   \int_{\mR^d}\mu^x_{s,r}(\dif y)\mu^y_{r,t}=\mu^x_{s,t},\quad \text{for all } 0\le s\le r\le t \le T \text{and $x\in\mR^d$},
 \end{align}
 } and 
 $$
 (B,\Sigma): [0,T]\times \mR^d\times \cC^d_\cP\times \cC^d_\cP\to(\mR^d,\mR^d\otimes\mR^d)
 $$ 
 are two Borel measurable functions. 
 
  The main feature of SDE \eqref{DFSDE} is that the coefficients depend on the distribution-flow $x\mapsto\mu^x_{s,t}$ of the solution itself, even the future distribution.
  Of course, one can regard $\mu^\centerdot_{s,t}$ as a probability kernel.
  Such type of SDEs naturally arises in the stochastic representation of Navier-Stokes equations as we shall see in the next subsection.
Before we continue the discussion, we first introduce the following notion of a solution to the above SDE:
\bd\label{Def11}
Let ${\frak F}:=(\Omega,\sF,(\sF_s)_{s\geq 0}, \mP)$ be a stochastic basis.
We call a pair of stochastic processes $((X^x_{s,t})_{(s,t,x)\in\mD_T\times\mR^d}, (W_t)_{t\in[0,T]})$ defined on ${\frak F}$ a weak solution of DFSDE \eqref{DFSDE}, if 
\begin{enumerate}[(i)]
\item $W_t$ is a standard $d$-dimensional $\sF_t$-Brownian motion;
\item For each $(s,t)\in\mD_T$,
$\mR^d\ni x \to \mu^x_{s,t}:=\mP\circ(X_{s,t}^x)^{-1}\in\cP(\mR^d)$ is weakly continuous {and the family $\mu^x_{s,t}$, $(s,t,x)\in \mD_T\times\mR^d$, satisfies \eqref{flow-1125}};
\item For each $(s,x)\in [0,T]\times\mR^d$,
\begin{align*}
\int_s^T |B(r,X_{s,r}^x,\mu_{r,T}^{\centerdot},\mu_{s,r}^{\centerdot})|\dif r+\int_s^T| \Sigma(r,X_{s,r}^x,\mu_{r,T}^{\centerdot},\mu_{s,r}^{\centerdot})|^2\dif r<\infty,\quad \mP-\text{a.s.}
\end{align*}
and the pair of processes a.e. satisfies equation \eqref{DFSDE} for all $t\in [s,T]$, $s\ge0$, $x\in\mR^d$.
\end{enumerate}
If, in addition, $X^x_{s,t}$ is adapted to the filtration generated by the Brownian motion $\sF^W_t:=\sigma\{W_r,r\in[0,t]\}$, then  it is called a strong solution of DFSDE \eqref{DFSDE}.
\ed

McKean-Vlasov SDEs, also referred to as distribution-dependent SDEs (DDSDEs), represent a significant class of stochastic differential equations where the coefficients are depending on the distribution of the solution process itself. These equations extend the scope of standard SDEs by incorporating the influence of interactions among particles or agents within a system. Initially introduced by Henry McKean  \cite{Mc67} in 1966 
in the context of nonlinear parabolic partial differential equations and later by Anatoli Vlasov \cite{Vl68} in 1968  in plasma physics, McKean-Vlasov SDEs have now attracted considerable attention across diverse fields such as mathematical finance, statistical physics, population dynamics, and mean field games (see, for example, \cite{BRTV98,BT97,CX10}).

Furthermore, DDSDEs exhibit substantial connections to vortex models like the Navier-Stokes and Euler equation (see, e.g., \cite{Osa86,FHM14}). DDSDEs with singular vortex kernels have been further developed by researchers such as Jabin-Wang \cite{JW18}, Serfaty \cite{Ser20b}, and other who at the same time contributed significantly to the advancement of propagation of chaos results.

For McKean-Vlasov SDEs, the dynamics of each individual particle is influenced by the collective behavior of the entire population, resulting in complex collective phenomena. This modeling framework allows the analysis of systems comprising a large number of interacting components, for which traditional approaches are inadequate. The general form of such McKean-Vlasov SDEs reads:
\begin{align}\label{DDSDE0}
Y_t=\xi+\int_0^t b(r,Y_r,\mu_r)\dif r+\int_0^t\sigma(r,Y_r,\mu_r)\dif W_r,
\end{align}
where $\mu_r$ denotes the distribution of $Y_r$. 
{We emphasize that these McKean-Vlasov SDEs differ from \eqref{DFSDE} in several key aspects. First, \eqref{DDSDE0} represents a single equation with a given initial condition $\xi$, whereas \eqref{DFSDE} describes a system of SDEs. Even if we consider \eqref{DDSDE0} with $\xi \sim \delta_x$, where $x \in \mathbb{R}^d$, and arbitrary starting times $s \geq 0$, the corresponding $\mu^x_{s,t}$ does not satisfy \eqref{flow-1125}. Instead, it only fulfills the following flow property:
$$
\mu^x_{s,t} = \mu^{\mu^x_{s,r}}_{r,t},
$$
which is the law of $Y_t$, where $Y$ is the solution of \eqref{DDSDE0} started at time $r \geq s$ with $\xi \sim \mu^x_{s,r}$.
}

Considerable attention has been paid to exploring the well-posedness of DDSDEs \eqref{DDSDE0} with singular drifts. Mishura and Veretenikov \cite{MV16} established the strong well-posedness of DDSDEs \eqref{DDSDE0} if the coefficient $b$ is only measurable and of at most linear growth, and additionally is Lipschitz continuous with respect to the distribution $\mu$, while $\sigma$ is assumed to be uniformly non-degenerate and Lipschitz continuous both in the spatial and measure vaiable. Later, R\"ockner and Zhang \cite{RZ21} extended this to cases involving local $L^q_tL^p_x$-drift. Additionally, Lacker \cite{La21} used the relative entropy method and Girsanov's theorem to obtain well-posedness results for DDSDEs with linear growth and $\sigma=\mI$, further extended by Han \cite{Ha22} to situations involving $L^q_tL^p_x$-drifts. Zhao \cite{Zhao2020} used heat kernel estimates and the Schauder-Tychonoff fixed-point theorem to establish well-posedness results for a more general class of DDSDEs with singular coefficients.

Through Zvonkin's transformation and the entropy method, the authors in \cite{HRZ22} proved the strong well-posedness for DDSDEs in cases where $\sigma$ is independent of $\mu$ and $b$ belongs to certain mixed $L^q_tL^p_x$-spaces. For specific cases, such as where $\sigma=\mI$ and $b(t,y,\mu)=b*\mu(t,y)$, both weak and strong well-posedness have been proved by various researchers \cite{CJM22, CJM23, HRZ23}. The Nemytskii-type DDSDEs, where 
$$
(b,\sigma)(t,y,\mu)=(b,\sigma)\Big(t,y,\tfrac{\mu(\dif y)}{\dif y}(y)\Big),
$$ 
 has been studied conducted by Barbu and R\"ockner \cite{BR18, BR20, BR22, BR22b}, and subsequently also in \cite{IR23}. Moreover, for kinetic cases of DDSDEs, further studies can be found in \cite{JW16, HRZ23, HZZZ21, IPRT24} and the references therein.

{Now let us return to the SDEs of type \eqref{DFSDE}. One of the main motivations for studying such equations arises from the Navier-Stokes equation, which provides an example of \eqref{DFSDE} through its stochastic representation.
}

\subsection{Motivation}
Consider the following Navier-Stokes equation on $\mR^d$ with $d=2,3$:
\begin{align}\label{NS0}
\left\{\begin{aligned}
&\p_t u=\Delta u-u\cdot \nabla u-\nabla p,\\
&\div u=0,\quad u_0=\varphi,
\end{aligned}\right.
\end{align}
where $u:[0,T]\times\mR^d\to\mR^d$ is the velocity field and $p$ stands for the pressure, and $\varphi$ is the initial velocity.
In \cite{CI08}, Constantin and Iyer presented the following probabilistic representation:
\begin{align}\label{NSF}
\left\{\begin{aligned}
&X_t^x=x+\int_0^t u(s,X_s^x)\dif s+\sqrt2W_t,\quad t\ge 0,\\
&u(t,x)={\bf P}\mE [\nabla^{\rm t}_xY^x_t\cdot \varphi(Y^x_t)],
\end{aligned}\right.
\end{align}
where $Y^x_t$ is the inverse of the flow mapping $x\to X_t^x$, $\nabla^{\rm t}$ denotes the transpose of the Jacobi matrix $(\nabla X)_{ij}:=\p_{x_j}X^i$,
and $\bP:=\mI-\nabla\Delta^{-1}\div$ is the Leray projection onto the space of divergence free vector fields. In particular, there is a one-to-one correspondence between \eqref{NS0} and \eqref{NSF}
when $u$ is smooth. Here an interesting question is how irregular $\varphi$ may be such that \eqref{NSF} admits a unique solution.

Now, for a velocity field $u$, let us consider its vorticity 
\begin{align*}
w=\text{curl} u=\begin{cases}
\p_2u_1-\p_1u_2,\quad &d=2;\\
\nabla \times u,\quad &d=3.
\end{cases}
\end{align*}
It is well-known that $u$ can be recovered from $w$ by the Biot-Savart law, i.e.,
$$
u=K_d*w,\ \ d=2,3,
$$ 
with
\begin{align}\label{BSK1}
K_2(x):=(x_2,-x_1)/(2\pi|x|^2),\ K_3(x)h=(x\times h)/(4\pi|x|^3).
\end{align}
Let $u$ be a smooth solution of \eqref{NS0}.
By direct calculations, we have
\begin{align}\label{1125:02}
\begin{split}
w(t,x)=\begin{cases}
\mE \left((\text{curl}\varphi)(Y^x_t)\det(\nabla_x Y^x_t)\right),\ & d=2,\\
\mE \left(\nabla^{\rm t}_xY^x_t\cdot (\text{curl}\varphi)(Y^x_t)\right),\ & d=3.
\end{cases}
\end{split}
\end{align}
By a change of variables and since $\det(\nabla Y^x_t)=1$, we get (cf. \cite{Zh15})
\begin{align}\label{1125:03}
\begin{split}
    u(t,x)=(K_d*w(t))(x)=
\left\{\begin{aligned}
&\mE \left(\int_{\mR^2}K_2(x-X_t^y)\cdot (\text{curl}\varphi)(y)\dif y\right),&d=2,\\
&\mE \left(\int_{\mR^3}K_3(x-X_t^y)\cdot\nabla_y X^y_t\cdot (\text{curl}\varphi)(y)\dif y\right),&d=3.
\end{aligned}\right.
\end{split}
\end{align}
In particular, for $d=2$, if we let
\begin{align*}
B(x,\mu^\centerdot):=\int_{\mR^2}(K_2*\mu^y)(x)\text{curl}\varphi(y)\dif y,
\end{align*}
and $\mu^y_t:=\mP\circ(X^y_t)^{-1}$, $t\in[0,T]$, then $X^x_t$ solves the following SDE:
\begin{align}\label{FNS0}
X_t^x=x+\int_0^t B(X^x_s,\mu^\centerdot_s)\dif s+\sqrt2W_t,
\end{align}
which leads to the system \eqref{DFSDE}.
{
\br
It should be noted that in both \eqref{1125:02} and \eqref{1125:03}, $w$ and $u$ do not depend linearly on the initial velocity $\varphi$, as the solution $X^y_t$ to the SDE \eqref{NSF} also depends on the initial velocity. 
\er

}

 DFSDE \eqref{FNS0} was introduced by Chorin \cite{Ch73} as the random vortex method to simulate viscous incompressible fluid flows for smooth kernels. Then it was further developed by Beale-Majda \cite{BM81}, Marchioror-Pulvirenti \cite{MP82} and Goodman \cite{Go87}. In particular, Long \cite{Lo88} showed the optimal convergence rate of the related particle system for mollifying kernels $K_2$. Later,
the interaction particle system and propagation of chaos related to \eqref{FNS0} have been attracted the attention of more and more investigators (see  \cite{FHM14, JW18}). 
However, the solvability of \eqref{FNS0} has not been tackled in the above references until the recent papers \cite{Zh23, CJM23, HRZ23} (See below for a further discussion).

If $d=3$,
formally, %$(X_t^x,\nabla X^x_t)$ can be regarded as a functional of the path $X_{[0,t]}^x$ and
\begin{align}\label{FNS00}
X_t^x=x+\int_0^t \bar\mE \left(\int_{\mR^3}K_3(X_s^x-\bar{X}_s^y)\cdot\nabla \bar{X}^y_s\cdot (\text{curl}\varphi)(y)\dif y\right)\dif s+\sqrt2W_t,
\end{align}
where $\bar{X}^y_t$ is an independent copy of ${X}^y_t$ and $\bar\mE$ is the expectation w.r.t. $\bar X^y_\cdot$ (see \cite{Zh15} and \cite{QSZ22} 
for its numerical simulations under the assumption of smoothness on the interaction kernel). 
To write down the above SDE in the form of \eqref{DFSDE}, we introduce a matrix-valued process $U^x_t:=\nabla X^x_t$.
Formally, $U$ solves the following linear ODE:
$$
U^x_t=\mI_{3\times 3}+\int_0^t \bar\mE \left({U^x_s\cdot\nabla}\int_{\mR^3} K_3(\cdot-\bar{X}_s^y)\cdot {\bar U}^y_s\cdot (\text{curl}\varphi)(y)\dif y\right)(X_s^x)\dif s.
$$
Let $(\mu^x)_{x\in\mR^3}$ be a family of probability measures over $\mR^3\times\cM^3$, where $\cM^3$ stands for the space of all $3\times 3$-matrices. Now let us introduce
$$
B(x,\mu):=\int_{\mR^3}\int_{\mR^3\times\cM^3}K_3(x-z)\cdot M \mu^y(\dif z\times\dif M)\cdot (\text{curl}\varphi)(y)\dif y.
$$
Then we obtain the following closed SDE
\begin{align}\label{3NS}
\left\{
\begin{aligned}
&X^x_{t}=x+\int_0^t B(X^x_{r},\mu^\centerdot_{r})\dif r+\sqrt2 W_t,\\
&U^x_{t}=\mI_{3\times 3}+\int_0^t{(U^x_{r}\cdot \nabla} B)(\cdot,\mu^\centerdot_{r})(X^x_{r})\dif r,
\end{aligned}
\right.
\end{align}
where $\mu^x_{t}:=\mP\circ(X^x_{t},U^x_{t})^{-1}\in\cP(\mR^3\times\cM^3)$ for $x\in\mR^3$.
{
\br\label{rmk:12}
We note that
\begin{align*}
    &\quad U^x_s\cdot\nabla\left(\int_{\mR^3} K_3(\cdot-\bar{X}_s^y)\cdot {\bar U}^y_s\cdot (\text{curl}\varphi)(y)\dif y\right)(X_s^x)\\
    &\ne \int_{\mR^3}(U^x_s\cdot \nabla K_3)(X_s^x-\bar{X}_s^y)\cdot {\bar U}^y_s\cdot (\text{curl}\varphi)(y)\dif y.
\end{align*}
%The derivation of the DFSDE \eqref{3NS} presented above lacks rigor. 
Specifically, for the gradients of the Biot-Savart kernel $K_d$, we note that $|\nabla K_d(x)|\lesssim |x|^{-d}\notin L^1_{loc}(\mR^d)$, which leads that
\begin{align*}
    (\nabla K_d)* f(x):=\lim_{\eps\to 0}\int_{|x-y|>\eps}\nabla K_d(x-y)f(y)\dif y
\end{align*}
is a Calder\'on-Zygmund operator. 

However, $\nabla (K_d*f)\ne (\nabla K_d)*f$.
For $d=2$ and $f:\mR^2\to\mR$, we have the following expression:
\begin{align*}
    \nabla_j (K_2*f)^i=(\nabla_j K_2^i)*f+\frac12 {\rm sign} (i-j)f, \quad i,j=1,2.
\end{align*}
For $d=3$, $f:\mR^3\to\mR^3$ and any $h\in\mR^3$, the expression is:
\begin{align*}
    h\cdot\nabla (K_3*f)(x)=& (h\cdot\nabla K_3)*f(x)+\frac13f(x)\times h,
\end{align*}
where explicitly,
\begin{align*}
    (h\cdot\nabla K_3)*f(x)
    ={\rm p.v.}\frac{3}{4\pi}\int_{\mR^3}\frac{[(x-y)\times f(y)]\otimes(x-y)}{|x-y|^5}h\dif y+{\rm p.v.}\frac{1}{4\pi}\int_{\mR^3}\frac{f(y)\times h}{|x-y|^3}\dif y.
\end{align*}
Further details and derivations of these results can be found in \cite[Section 2.4.3, p. 76]{MB02}.

%{\blue Specially, for 
%\begin{align*}
    %B(x,\mu)=\int_{\mR^d} K_3(x-z) f(z)\dif z,\quad %f(x)=\int_{\mR^3\times 
    %\cM^3} M\rho^y(x,M)\cdot ({\rm curl} \varphi)%(y)\dif y\dif M,
%\end{align*}
%where we assume that $\rho^y(z,M):=\mu^y(\dif z\times \dif M)/(\dif z\times \dif M)$ exists. Then the equation \eqref{}}

In the present paper, we only consider the case $d=2$. A %rigorous derivation and 
detailed investigation of \eqref{3NS} for $d=3$ will be addressed in future work.
\er
}

On the other hand, if we set $\tilde u(t,x):=-u(T-t,x)$ and $\tilde p(t,x):=p(T-t,x)$, then $\tilde u$ solves the following backward Navier-Stokes equation:
\begin{align*}
\begin{cases}
\p_t \tilde u+\Delta \tilde u+\tilde u\cdot \nabla \tilde u+\nabla \tilde p=0,\\
\div \tilde u=0,\quad \tilde u_T=\varphi.
\end{cases}
\end{align*}
In \cite{Zh12}, the author provided a probabilistic representation for $\tilde u$ as well:
\begin{align}\label{NSB}
\left\{\begin{aligned}
&\tilde X_{s,t}^x=x+\int_s^t \tilde u(r,\tilde X_{s,r}^x)\dif r+\sqrt2(W_t-W_s),\quad (s,t)\in\mD_T,\\
&\tilde u(t,x)={\bf P}\mE[\nabla^{\rm t} \tilde X_{t,T}^x\cdot \varphi(\tilde X_{t,T}^x)].
\end{aligned}\right. 
\end{align}
As above, in the two dimensional case, we have
\begin{align*}
\tilde w(t,x):= \text{curl} \tilde u(t,x)=\mE[ (\text{curl}\varphi)(\tilde X_{t,T}^x)]=\<\text{curl}\varphi,\tilde \mu_{t,T}^x\>,
\end{align*}
where $\tilde \mu_{s,t}^x:=\mP\circ(X_{s,t}^x)^{-1}$. By the Biot-Savart law, we have
\begin{align*}
\tilde u(t,x)=(K_2* \tilde w(t))(x)= \int_{\mR^2} K_2(x-y) \<\text{curl}\varphi,\tilde \mu_{t,T}^y\>\dif y.
\end{align*} 
Thus \eqref{NSB} is transformed into the following DFSDE:
\begin{align}\label{DFSDE2}
\tilde X_{s,t}^x=x+\int_s^t B(\tilde X_{s,r}^x,\tilde \mu_{r,T}^\centerdot)\dif r+\sqrt2(W_t-W_s),
\end{align}
where
\begin{align*}
B(x,\mu^\centerdot)=K_2*\left(\int_{\mR^2}\text{curl}\varphi(y)\mu^\centerdot(\dif y)\right)(x).
\end{align*}
In particular, SDE \eqref{DFSDE2} is exactly an example of DFSDE \eqref{DFSDE} with $B(r,x,\mu^\centerdot,\nu^\centerdot)=B(x,\mu^\centerdot)$.
For the three dimensional case, as in \eqref{3NS}, we have the following representation:
\begin{align}
\left\{
\begin{aligned}
&\tilde X^x_{s,t}=x+\int_s^t  B(\tilde X^x_{s,r},\mu^\centerdot_{r,T})\dif r+\sqrt2(W_t-W_s),\\
&\tilde U^x_{s,t}=\mI_{3\times 3}+\int_s^t (\tilde U^x_{s,r}\cdot\nabla) B(\cdot,\mu^\centerdot_{r,T})(\tilde X^x_{s,r})\dif r,
\end{aligned}
\right.
\end{align}
where $\mu^x_{s,t}:=\mP\circ(\tilde X^x_{s,t},\tilde U^x_{s,t})^{-1}\in\cP(\mR^3\times\cM^3)$, and
$$
B(x,\mu):=\int_{\mR^3}K_3(x-y)\left(\int_{\mR^3\times\cM^3}M^{\rm t}\cdot (\text{curl}\varphi)(z)\mu^y(\dif z\times\dif M)\right)\dif y.
$$
We must point out that \eqref{FNS0} and \eqref{DFSDE2} are essential different as we discuss in the next subsection.

\subsection{Main results}

Our first result is about the strong well-posedness of DFSDE \eqref{DFSDE} with regular coefficients. More precisely, 
let $\cC\cP_1$ be the space of all continuous probability measure-valued functions from $\mR^d$ to $\cP_1(\mR^d)$ with finite first order moment (see Section 2 for more details about
the space $\cC\cP_1$). We assume that $B$ and $\Sigma$ satisfy the following assumptions:
\begin{itemize}
\item[{\bf(H$_0$)}] 
For each $t\in[0,T]$, the function 
$$
\mR^d\times\cC\cP_1\times\cC\cP_1\ni (x,\mu^\centerdot,\nu^\centerdot )\mapsto (B,\Sigma)(t,x,\mu^\centerdot,\nu^\centerdot)\in(\mR^d,\mR^d\otimes\mR^d)\mbox{is  continuous,}
$$
and there are constants $\kappa_0, \kappa_2,\kappa_3,\kappa_4>0$ and $\kappa_1\in\mR$ such that
for any $(t,x,\mu,\nu)\in [0,T]\times\mR^d\times \cC\cP_1\times\cC\cP_1$,
\begin{align}\label{AA05}
\<x, B(t,x,\mu^\centerdot,\nu^\centerdot)\>+2\|\Sigma(t,x,\mu^\centerdot,\nu^\centerdot)\|^2_{\rm HS}\le \kappa_0+\kappa_1|x|^2
+\kappa_2(\|\mu^\centerdot\|^2_{\cC\cP_1}+\|\nu^\centerdot\|^2_{\cC\cP_1}),
\end{align}
and for any $(t,x_i,\mu_i,\nu_i)\in [0,T]\times\mR^d\times \cC\cP_1\times\cC\cP_1$, $i=1,2$,
\begin{align}\label{AA06}
\begin{split}
\qquad&\<x_1-x_2, B(t,x_1,\mu_1^\centerdot,\nu_1^\centerdot)-B(t,x_2,\mu_2^\centerdot,\nu_2^\centerdot)\>
+2\|\Sigma(t,x_1,\mu_1^\centerdot,\nu_1^\centerdot)-\Sigma(t,x_2,\mu_2^\centerdot,\nu_2^\centerdot)\|^2_{\rm HS}\\
&\qquad\qquad\le \kappa_3|x_1-x_2|^2+\kappa_4(1+|x_1|^2+|x_2|^2)\left({\rm d}^2_{\cC\cP_1}(\mu^\centerdot_1,\mu^\centerdot_2)
+{\rm d}^2_{\cC\cP_1}(\nu^\centerdot_1,\nu^\centerdot_2)\right),
\end{split}
\end{align}
where $\|\cdot\|_{\rm HS}$ stands for the Hilbert-Schmit norm and the distance ${\rm d}_{\cC\cP_1}$ is defined in \eqref{DD1} below. 
\end{itemize}

Our first main result is the following strong well-posedness, which is proven by freezing the distribution-flow and Picard's iteration.
\bt\label{Thm21}
Under {\bf(H$_0$)}, there is a unique strong solution to DFSDE \eqref{DFSDE} in the sense of Definition \ref{Def11}. 
Moreover, there is a constant $C_T=C_T(\kappa_i)>0$ such that for all $(s,t,x)\in\mD_T\times\mR^d$, 
$$
\mE|X^{x}_{s,t}|^2\leq C_T(1+|x|^2),
$$
and if $\kappa_1<0$ and $\kappa_1+2\kappa_2<0$, then
\begin{align}\label{AA106}
\mE|X^{x}_{s,t}|^2\leq \e^{\kappa_1 (t-s)}|x|^2+ (\kappa_0+\kappa_5)(\e^{\kappa_1(t-s)}-1)/\kappa_1,
\end{align}
where $\kappa_5:=2\kappa_2(|\kappa_1|+\kappa_0)/(|\kappa_1|-2\kappa_2)$.
\et

Our second main result is about the well-posedness of {DFSDE related to the} 2D-Navier-Stokes equation {driven by} the fractional Brownian motion ({\rm f}Bm).
Recall that a  Gaussian process $(W^H_t)_{t\geq 0}$  is called an {\rm f}Bm with Hurst parameter $H\in(0,1)$ if for any $0\le s\le t$,
$$
\mE(W^{H}_tW^{H}_s)=\tfrac12(t^{2H}+s^{2H}-|t-s|^{2H}).
$$
Clearly, $W^H$ has the following self-similarity:
for $\lambda>0$,
$$
(W^H_t)_{t\geq 0}\stackrel{(d)}{=}(\lambda^{-H}W^H_{\lambda t})_{t\geq 0}.
$$
Consider the following {DFSDE related to the} 2D-Navier-Stokes equation {driven by} {\rm f}Bm:
\begin{align}\label{FNS10}
X_t^x=x+\int_0^t \int_{\mR^2}(K_2*\mu_s^y)(X^x_s)\nu_0(\dif y)\dif s+ W^H_t,
\end{align}
where $K_2$ is the Biot-Savart law given in \eqref{BSK1}, $\mu^y_t=\mP\circ (X^y_t)^{-1}$, $\nu_0$ is a finite signed measure on $\mR^2$ and $W^H=(W^{H,1}, W^{H,2})$ with that $W^{H,i}, i=1,2$ are two independent {\rm f}Bms with the same Hurst parameter $H$. We have the following weak well-posedness  {(see Theorem \ref{thm:51} below for a detailed statement and its proof)}.
\bt\label{thm1}
Let $H\in(0,\frac12)$.
For any %initial 
vorticity $\nu_0$ being a finite singed measure, there is a unique weak solution $X^\cdot_t$ to SDE \eqref{FNS10}.
{Moreover, for any $p\in(1,2)$ and $\eps>0$, there is a constant $C>0$ such that for all $0<t\le T$,
\begin{align*}
\nor u(t)-K_2*\nu_0\nor_{p}\le C t^{[H(\frac{2}{p}-1)]\wedge[\frac{1-2H}{1-H}]-\eps}.
\end{align*}
Here the localized $L^p$ norm $\nor\cdot\nor_p$ is defined in \eqref{local-1125}.
Furthermore, if we let
\begin{align*}
u(t,x):=\int_{\mR^2} \mE K_2(x-X_{t}^y)\nu_0(\dif y),
\end{align*}
then $u\in C((0,T]; C^\infty_b(\mR^2))$. }
\et
\br
Since {\rm f}Bm is neither a Markov process nor a martingale, one can not say that $u$ solves any PDE.
%But by the probabilistic representation, we shall call $u$ the solution of 2D $H$-fractional  Navier-Stokes equation with  initial vorticity $\nu_0$.
By the change of variable,
the above $u$ has the following scaling property: for $\lambda>0$, if we let
$$
u_\lambda(t,x):=\lambda^{1/H-1} u(\lambda^{1/H} t,\lambda x),
$$
then $u_\lambda(t,x)=\int_{\mR^2} \mE K_2(x-X_{t}^{y;\lambda})\nu^\lambda_0(\dif y)$, where 
\begin{align}\label{SU0}
\nu^\lambda_0(\dif y)=\lambda^{1/H-2}\nu_0(\dif (\lambda y)),
\end{align}
and
$X_{t}^{x;\lambda}$ solves the following DFSDE:
$$
X_{t}^{x;\lambda}=x+\int_0^t\!\! \int_{\mR^2}(K_2*\mu_s^{y;\lambda})(X^{x;\lambda}_s)\nu^\lambda_0(\dif y)\dif s+W^H_t.
$$
If $\nu_0(\dif y)=\varrho(y)\dif y$, then \eqref{SU0} reduces to $\nu^\lambda_0(\dif y)=\lambda^{1/H}\varrho(\lambda y)\dif y$.
\er

Note that Theorem \ref{thm1} does not include  $H=\frac12$. Next we consider the following backward version of {DFSDE related to the} Navier-Stokes equation {driven by} Brownian motion:
\begin{align}\label{NSSDE03}
X_{s,t}^x=x+\int_s^t\!\!\int_{\mR^2}K_2(X^x_{s,r}-y)\mu^y_{r,T}(g)\dif y\dif r+\sqrt2(W_t-W_s),
\end{align}
In this case, we also have
\bt\label{thm2}
Let $p_0\in(1,2)$ and $g\in \mL^{p_0}$.
For each $s\in[0,T]$ and $x\in\mR^2$, there is a unique strong solution $X^x_{s,t}$ to DFSDE \eqref{NSSDE03}.
Moreover, if we let
\begin{align*}
u(s,x)=\int_{\mR^2} K_2(x-y)\mE g(X_{s,T}^y)\dif y,
\end{align*}
then $u\in C([0,T); C^\infty_b(\mR^2))$ solves the following backward Navier-Stokes equation:
$$
\p_su+\Delta u+u\cdot \nabla u+\nabla p=0,\ \ u(T)=K_2*g.
$$  
\et
{The proof of this theorem is provided in the proof of Theorem \ref{thm:54}.}
\br
When $g\in \mL^{p_0}$ with $p_0>2$, the well-posedness of DDSDE \eqref{NSSDE03} was obtained in \cite{Zh15} by Zvonkin's transformation. Here we regard 
\eqref{NSSDE03} as an abstract distribution-flow SDE. 
\er
\subsection{Related works} 

In this section, we review some literature relevant for our main Theorems \ref{thm1} and \ref{thm2}. We begin by considering a classical DDSDE with a singular Biot-Savart interaction kernel, often referred to as the 2D random vortex model:
\begin{align}\label{0224:02}
X_{t}=X_0+\int_0^t(K_2*\mu_s)(X_{s})\dif s+\sqrt 2W_t,
\end{align}
where $K_2$ is the Biot-Savart kernel, $\mu_t$ denotes the law of $X_t$, and $\mP\circ X_0^{-1}(\dif y)=\mu_0(\dif y)$. Suppose $\mu_t(\dif y)=\rho_t(y)\dif y$. Then, by It\^o's formula, $\rho_t$ satisfies the following vorticity form of the Navier-Stokes equation \eqref{NS0}:
\begin{align}\label{1125:00}
    \p_t\rho=\Delta\rho-\div((K_2*\rho)\rho).
\end{align}
In other words, $u(t,x):=K_2*\rho_t(x)$ solves the Navier-Stokes equation \eqref{NS0}.
In \cite{Zh23}, the second author utilized the De-Giorgi method to establish the existence of a weak solution to \eqref{0224:02} when $X_0=x$, while the uniqueness remains an open question. {Weak existence for \eqref{0224:02} for $X_0=x$ was also proved in \cite{BRZ23} using a nonlinear variant of the superposition principle (see \cite{Tr16}). Furthermore, letting $X^{\zeta}_{r,t}$ denote the weak solution to \eqref{0224:02} started at $r\ge0$, with $X_0\sim \zeta$, it was proved in \cite{BRZ23}, that the path law $P_{(s,\zeta)}$, $s\ge0$, $\zeta\in \cP_1$, of $(X^{\zeta}_{r,t})_{t\ge r}$ form a nonlinear Markov process in the sense of McKean \cite{Mc66}. Moreover , it is proved in \cite{BRZ23} that if $X_0$ has a density $\rho_0\in L^4$, then \eqref{0224:02} has a strong solution and that pathwise uniqueness holds for \eqref{0224:02} in the class of all solutions having time marginal law densities in $L^{4/3}([0,T];L^{4/3})$. Furthermore,} if the initial data $X_0$ admits a density $\rho_0\in \mL^{1+}$ with respect to the Lebesgue measure, weak and strong well-posedness for \eqref{0224:02} were established in \cite{CJM23} and \cite[Theorem 6.4]{HRZ23}.

It's important to note that if the initial vorticity of the Navier-Stokes equation is not a probability measure, then there is no one-to-one correspondence between \eqref{0224:02} and \eqref{NS0}. To address this issue, we consider the forward and backward random vortex models, as introduced in Section 1.1:
\begin{align}\label{0224:01}
X^x_{s,t}=x+\int_s^t\!\! \int_{\mR^{4}}K_2(X^x_{s,r}-y)g(z)\mu_{s,r}^z(\dif y)\dif z+W^H_t-W_s^H
\end{align}
and
\begin{align}\label{0224:00}
X^x_{s,t}=x+\int_s^t\!\! \int_{\mR^{4}}K_2(X^x_{s,r}-y)g(z)\mu_{r,T}^y(\dif z)\dif y+W_t^H-W_s^H,
\end{align}
where $g$ represents the initial vortex, $\mu^x_{s,t}$ denotes the time marginal law of the solution $X^x_{s,t}$, and $W^H_s$ is the fractional Brownian motion with $H\in(0,\frac12]$ (see Section \ref{Sec:fBM} for details).

For the forward DFSDE \eqref{0224:01}, we focus on $H\in(0,\frac12)$. Recent advancements in regularization by averaging paths (see \cite{CG16,GG22}) and the stochastic sewing lemma (see \cite{Le20,Le23})  have led to an increased interest in the well-posedness of SDEs driven by fractional Brownian motion of the form:
\begin{align*}
\dif X_t=b(t,X_t)\dif t+\dif W_t^H,
\end{align*}
where $b\in \mL^q_t\mL^q_x$. Several classical results have been established, such as those by Nualart-Ouknine in \cite{NO03} and L\^e in \cite{Le20}, where
Nualart-Ouknine, using the Girsanov transformation, established weak well-posedness for $p,q\ge2$ and $1/q+Hd/p<1/2$, and L\^e in \cite{Le20} 
extended this result by introducing the stochastic sewing lemma and gave a new proof for the weak well-posedness.
The strong well-posedness was obtained as well in \cite{Le20} when $1/q+Hd/p<1/2-H$. 
Furthermore, Galeati and Gubinelli in \cite{GG22} employed the averaging paths technique to achieve path-by-path well-posedness for $q=\infty$ and $Hd/p<1/4-H$. 
In fact, upon assuming $X_t^\eps := \eps^{-H} X_{\eps t}$ and $b_\eps(t,x) := \eps^{1-H}b(\eps t, \eps^H x)$ with some $\eps>0$, we have
\begin{align*}
\dif X_t^\eps = b_\eps(t,X_t^\eps)\dif t + \eps^{-H}\dif W_{\eps t}^H,
\end{align*}
where $\eps^{-H}W_{\eps t}^H$ remains an {\rm f}Bm with the Hurst index $H$. 
Note that the scaling hypothesis $\lim_{\eps\to 0}\|b_\eps\|_{\mL^q_t\mL^p_x}=0$ leads to
\begin{align}\label{0224:04}
\tfrac{1}{q}+\tfrac{Hd}{p}<1-H.
\end{align}
Under condition \eqref{0224:04}, Butkovsky, L\^e and Mytnik \cite{BLM23} recently established the existence of a solution for $q=\infty$, and when $q\in(1,2]$, 
Galeati and Gerencs\'er \cite{GG23} demonstrated the strong well-posedness. It is noteworthy that when $b$ is independent of time $t$, as in the case of the Biot-Savart kernel, the condition in \cite{GG23} simplifies to $\tfrac{1}{2}+\tfrac{Hd}{p}<1-H$, aligning with the strong well-posedness result in \cite{Le20}. It is important to note that the review here is primarily focused on the $\mL^q_t\mL^p_x$-drift. Indeed, 
both \cite{BLM23} and \cite{GG23} cover various measure and distributional cases respectively. More recently, Butkovsky and Gallay in \cite{BG23}, employing a combination of the stochastic sewing lemma and John-Nirenberg's inequality, established the existence of solutions under the condition $(1-H)/q+Hd/p<1-H$, which is considerably weaker.
Beyond these results, a lot of related works exists, and interested readers can refer to the comprehensive overview in \cite{GG23}.

In the context of the following DDSDE driven by {\rm f}Bm
\begin{align*}
\dif X_t=(b*\mu_t)(t,X_t)\dif t+\dif W_t^H,
\end{align*}
where $\mu_t$ represents the time marginal law of $X_t$, the authors in \cite{GHM23} and \cite{GG23} established strong well-posedness for $b\in \mL^q_t\bC^{\alpha}$ with  $\alpha>1+1/(Hq)-1/H$ and $q\in(1,2]$.
Here, $\bC^{\alpha}$ denotes the Besov space. Moreover, Han \cite{Ha22} used the entropy method to present a concise proof of the main results in \cite{GHM23}.

Following this review, we examine the condition on $H$ for the Biot-Savart kernel by applying the aforementioned results. Notably, $b=K_2\in L^{2-}_{loc}\cap \bC^{-1}$. 
Therefore, the restriction $p\ge2$ in \cite{NO03} precludes its application to the Biot-Savart law. Moreover, the conditions in \cite{Le20} and \cite{GG23} also imply that $H$ must be strictly less than $1/4$. 
Consequently, it is natural to inquire whether the well-posedness holds for \eqref{0224:01} in the range $H\in[1/4,1/2)$. We address this question in Theorem \ref{thm1} 
by establishing weak well-posedness for \eqref{0224:01} across all $H\in(0,1/2)$. Additionally, we define a solution to the 2D fractional Navier-Stokes equation with an arbitrary initial vortex measure $\nu_0$ and show its smoothness for $t>0$ by the Malliavin calculus.

For the backward DFSDE \eqref{0224:00}, limited results are available in the literature. In Theorem \ref{thm2}, we establish the unique strong flow solution $X^x_{s,t}$ for $H=1/2$ and any $\mL^{1+}$ initial data $\varphi_0$. However, for $H\ne 1/2$, investigating the well-posedness of \eqref{0224:00} becomes challenging, as our methodology heavily relies on the Markov property of Brownian motion. Notably, the Girsanov transformation cannot be used for backward DFSDE \eqref{0224:00} with singular kernels, leaving this as an open question.

\subsection{Organization of the paper}

In Section 2, we provide preliminary results concerning the space of probability kernels and {\rm f}Bms. Some of these results are novel and are crucial in proving Theorems \ref{Thm21}, \ref{thm1}, and \ref{thm2}.

Section 3 is dedicated to proving Theorem \ref{Thm21} using standard Picard's iteration. Additionally, we offer several examples to illustrate our main results. While the proof itself is not particularly challenging, it serves as a foundation for our future investigations into various issues such as ergodicity and propagation of chaos.

Section 4 focuses on demonstrating weak and strong well-posedness for a broad class of DFSDEs driven by {\rm f}Bms. 
For  $H\in(0,\frac12)$, we employ Girsanov's transformation and the entropy method, while for $H=\frac12$, we rely on PDE estimates.

In Section 5, we utilize the results obtained in Section 4 to prove Theorems \ref{thm1} and \ref{thm2}. To establish the smoothness of the velocity field, we employ Malliavin calculus when 
$H\in(0,\frac12)$ and PDE techniques for $H=\frac12$.

In the Appendix, we provide detailed proofs of certain technical results for the convenience of the readers.

Throughout this paper, we shall use the following convention and notations: The letter $C$ with or without subscripts will denote an unimportant constant, whose value
may change from line to line. We also use $:=$ to indicate a definition and set
$$
a\wedge b:=\max(a,b),\ \ a\vee b:=\min(a,b),\ \ a^+:=0\vee a.
$$ 
By $A\lesssim_C B$ and $A\asymp_C B$
or simply $A\lesssim B$ and $A\asymp B$, we respectively mean that for some constant $C\geq 1$,
$$
A\leq C B,\ \ C^{-1} B\leq A\leq CB.
$$
Below we collect some frequently used notations for the readers' convenience.
\begin{itemize}
\item For $p\in[1,\infty]$, $p'$ denotes the conjugate index of $p$, i.e., $\frac1p+\frac1{p'}=1$.
\item $\cP_1$: The space of all probability measures with finite first order moment.
\item $\cC\cP_1$: The space of $\cP_1$-valued continuous functions on $\mR^d$ w.r.t. the Wasserstein-1 distance.
\item $\cC\cP_0$: The space of $\cP$-valued continuous functions on $\mR^d$ w.r.t. the total variation distance.
\item $\cL^p\cP_s$: The space of sub-probability kernels defined in \eqref{SPK1}.
\item $\tilde\cL^p\cP$: The space of probability kernels defined in \eqref{SPK2}.
\item $\mH^q_T$: The space of all absolutely continuous function $f:[0,T]\to\mR^d$ with $f(0)=0$ and $\dot f\in L^q([0,T];\mR^d)=:\mL^q_T$.
\end{itemize}

\section{Prelimiaries}

In this section, we first introduce several spaces of flow probability measures associated with the Wasserstein-1 metric, the total variation distance, and localized 
$L^p$-probability kernels. Then, we also recall the definition and basic properties of fractional Brownian motions ({\rm f}Bms). In particular, we demonstrate an important exponential estimate for the functional of {\rm f}Bm, which is crucial for solving  singular SDEs driven by {\rm f}Bms using Girsanov's theorem.
 
\subsection{Flow probability measure space}
Let $\cP_1:=\cP_1(\mR^d)$ be the space of all probability measures with finite first order moment and $\cC\cP_1$ the space of $\cP_1$-valued continuous functions on $\mR^d$ w.r.t. the Wasserstein-1 
distance $\cW_1$ defined by
$$
\cW_1(\mu,\nu):=\inf_{\mP\circ X^{-1}=\mu,\mP\circ Y^{-1}=\nu}{\mE|X-Y|}. 
$$
Note that by the duality of Monge-Kantorovich (cf. \cite[(6.3)]{Vi06}), 
$$
\cW_1(\mu,\nu)=\sup_{\|g\|_{\rm Lip}\leq 1}|\mu(g)-\nu(g)|.
$$
For two $\mu^\centerdot,\nu^\centerdot\in \cC\cP_1$, we introduce a distance between $\mu^\centerdot$ and $\nu^\centerdot$ by
\begin{align}\label{DD1}
{\rm d}_{\cC\cP_1}(\mu^\centerdot,\nu^\centerdot):=\sup_{x\in\mR^d}\frac{\cW_1(\mu^x,\nu^x)}{1+|x|}.
\end{align}
For simplicity, we write
$$
\|\mu^\centerdot\|_{\cC\cP_1}:={\rm d}_{\cC\cP_1}(\mu^\centerdot,\delta_0)=\sup_{x\in\mR^d}\frac{\int_{\mR^d}|y|\mu^x(\dif y)}{1+|x|}.
$$
Moreover, the total variation distance $\|\cdot\|_{\var}$ is defined by
\begin{align*}
\|\mu-\nu\|_{\var}:=\sup_{A\in\sB(\mR^d)}|\mu(A)-\nu(A)|.
\end{align*}
Let $\cC\cP_0$ be the space of all continuous probability kernels $x\mapsto \mu^x$ w.r.t $\|\cdot\|_{\rm var}$ with the distance
$$
\|\mu^\centerdot-\nu^\centerdot\|_{\cC_\var}:=\sup_{x\in\mR^d}\|\mu^{x}-\nu^{x}\|_{\var}.
$$
Under these distances, it is easy to see that $\cC\cP_1$ and $\cC\cP_0$ are complete metric spaces. We would like to point out that the distance ${\rm d}_{\cC\cP_1}$ is only used in the study of DFSDEs with
regular coefficients.

Next we introduce some localized $L^p$-spaces for later use.  Let $(D_i)_{i\in\mN}$ be the set of all unit cubes with center at the integer lattice so that
\begin{align}\label{BC7}
0\in D_1,\ \ \cup_{i\in\mN} D_i=\mR^d,\ \ D_i\cap D_j=\emptyset,\ \ i\not=j.
\end{align}
For $z=(z_1,\cdots,z_d)\in\mR^d$, we shall write
$$
D_z:=D_1+z=\big\{x: -1/2\leq x_i-z_i<1/2\big\}.
$$
For $p\in[1,\infty]$, let $\mL^p=L^p(\mR^d)$ be the usual $L^p$-space {with norm $\|\cdot\|_p$}. We also introduce the Banach spaces
\begin{align}\label{local-1125}
    \tilde\mL^p:=\left\{f\in L^p_{loc}(\mR^d): \nor f\nor_p:=\sup_i\|\1_{D_i}f\|_p<\infty\right\}
\end{align}
and
$$
\bar\mL^p:=\left\{f\in L^p_{loc}(\mR^d): \nor f\nor^*_p:=\sum_i\|\1_{D_i}f\|_p<\infty\right\}.
$$
By a finite covering technique,  there is a constant $C_1=C_1(p,d)>1$ such that
\begin{align}\label{BC4}
C_1^{-1}\nor f\nor_p\leq\sup_z\|\1_{D_z} f\|_p\leq C_1\nor f\nor_p.
\end{align}
The advantage of using localized space $\tilde\mL^p$ is  the following inclusion: for $p_1\geq p_2$,
$$
\tilde\mL^{p_1}\subset\tilde\mL^{p_2}.
$$
This is quite convenient for treating singular potentials like $|x|^{-\alpha}$, where $\alpha\in(0,d)$, since it does not belong to any $L^p$-space, but belongs to $\tilde\mL^p$ for $p<\frac d\alpha$.
About the spaces $\tilde\mL^p$ and $\bar\mL^p$, we have the following properties, that are similar to the classical $L^p$-spaces. For the readers' convenience, 
we provide detailed proofs in  Appendix A.
\bp\label{Le20}
\begin{enumerate}[(i)]
\item For each $p\in[1,\infty]$, it holds that $\bar\mL^p\subset\mL^p\subset\tilde\mL^p$, and 
\begin{align}\label{AW1}
\nor f\nor_p\asymp\sup_{\nor g\nor^*_{p'}\leq 1}\int_{\mR^d} f(x)g(x)\dif x,\ \ \nor g\nor^*_p\asymp\sup_{\nor f\nor_{p'}\leq 1}\int_{\mR^d} f(x)g(x)\dif x,
\end{align}
where $p'$ is the conjugate index of $p$.
\item For any $p,q,r\in[1,\infty]$ with $1+\frac1r=\frac1p+\frac1q$,  the following Young's inequalities hold: for some $C=C(d,p,q,r)>0$,
\begin{align}\label{AW2}
\nor f*g\nor_r\leq C\nor f\nor_p\nor g\nor^*_q,\quad \ \nor f*g\nor^*_r\leq C\nor f\nor^*_p\nor g\nor^*_q.
\end{align}
\end{enumerate}
\ep

Finally, we introduce a space of  probability kernels that will be used in the study of backward DFSDEs.
Let  $\cK\cP$ be the set of all probability kernels from $\mR^d$ to $\cP$ and $\cK\cP_s$ the set of all sub-probability kernels from $\mR^d$ to $\cP_s$,
where $\cP_s$ is the space of all sub-probability measures over $\mR^d$.
For given $p\in[1,\infty]$, we introduce two subclasses
\begin{align}\label{SPK1}
\cL^p\cP_s:=\Big\{\mu^\centerdot\in \cK\cP_s: \|\mu^\centerdot\|_p:=\sup_{\|\phi\|_p\leq 1}\|\mu^\centerdot(\phi)\|_p<\infty\Big\}
\end{align}
and
\begin{align}\label{SPK2}
\tilde\cL^p\cP:=\Big\{\mu^\centerdot\in \cK\cP: \nor \mu^\centerdot\nor_p:=\sup_{\nor\phi\nor_p\leq 1}\nor \mu^\centerdot(\phi)\nor_p<\infty\Big\}.
\end{align}
Similarly, we also introduce the subclasses $\cL^p\cP$ and $\tilde\cL^p\cP_s$.
It is easy to see that 
$$
\cL^p\cP\subset \cL^p\cP_s,\ \ \tilde\cL^p\cP\subset \tilde\cL^p\cP_s.
$$ 

\br
{Here $\sup_{\|\phi\|_p\leq 1}\|\mu^\centerdot(\phi)\|_p<\infty$ means that for any $\phi\in L^p$, $\mu^x(\phi)=\int_{\mR^d}\phi(y)\mu^x(\dif y)$ is well-defined for Lebesgue almost {all} $x\in\mR^d$ and that $\mu^\centerdot(\phi)$ belongs to $L^p(\mR^d)$. 
This condition implies that $\mu^x(\phi)=\int_{\mR^d}\phi(y)\mu^x(\dif y)$ is independent of the representative of $\phi$ in $L^p(\mR^d)$. Indeed, if $\phi(x)=\tilde{\phi}(x)$ for Lebesgue almost {all} $x\in\mR^d$, then $\mu^\centerdot(\phi)=\mu^\centerdot(\tilde{\phi})$ in $L^{p}(\mR^d)$. This follows from the estimate $ \|\mu^\centerdot(\phi-\tilde{\phi})\|_p\lesssim\|\phi-\tilde{\phi}\|_p=0$.}
\er
%In a natural way, $\tilde\cL^p\cP$ {\blue and $\tilde\cL^p\cP_s$} can be regarded as a subset of the space $\cL(\tilde\mL^p,\tilde\mL^p)$ that consists of all bounded linear operators from $\tilde\mL^p$ to $\tilde\mL^p$.
%Moreover, in fact 
One would like to point out that such classes of kernels naturally appear in the study of stochastic Lagrangian flows (see \cite{Zh13}). More precisely, consider the following SDE:
$$
X^x_t=x+\int^t_0b(X^x_s)\dif s+W_t,
$$
where $b$ is a divergence free Lipschitz vector field. Let $\mu^x_t$ be the law of the unique solution $X^x_t$. It is well-known that
for any $p\in[1,\infty]$ and $t\geq 0$  (see \cite{Zh13}),
$$
\|\mu^\centerdot_t(f)\|_p\leq\|f\|_p,\ \  f\in\mL^p.
$$

We have the following important properties that are also proven in Appendix A.

\bp\label{Le21}
(i) Let $\mu^\centerdot\in\cK\cP_s$. For any $p\in[1,\infty)$,
we have
\begin{align}\label{AW202}
\|\mu^\centerdot\|_p=\sup_{\phi\in C_c(\mR^d),\|\phi\|_p\leq 1}\| \mu^\centerdot(\phi)\|_p=\sup_{\phi\in C^\infty_c(\mR^d),\|\phi\|_p\leq 1}\| \mu^\centerdot(\phi)\|_p,
\end{align}
and
\begin{align}\label{AW02}
\nor \mu^\centerdot\nor_p=\sup_{\phi\in C_c(\mR^d),\nor\phi\nor_p\leq 1}\nor \mu^\centerdot(\phi)\nor_p=\sup_{\phi\in C^\infty_c(\mR^d),\nor\phi\nor_p\leq 1}\nor \mu^\centerdot(\phi)\nor_p.
\end{align}
%When substituting the norm $\nor\cdot\nor_p$ with $\|\cdot\|_p$, the equality \eqref{AW02} still holds.\\
(ii) ${\cL}^p\cP_s$ and $\tilde{\cL}^p\cP_s$ are complete metric spaces with respect to the distance
\begin{align*}
\|\mu^\centerdot-\nu^\centerdot\|_p:=\sup_{\| \phi\|_p\le 1}\| \mu^\centerdot(\phi)-\nu^\centerdot(\phi)\|_p,\ \ 
\nor \mu^\centerdot-\nu^\centerdot\nor_p:=\sup_{\nor \phi\nor_p\le 1}\nor \mu^\centerdot(\phi)-\nu^\centerdot(\phi)\nor_p.
\end{align*}
Moreover, with the above distances,  $\tilde{\cL}^p\cP$ is still complete, but ${\cL}^p\cP$ is not complete.
 \ep

\br
Note that $\|\mu^\centerdot-\nu^\centerdot\|_\infty=\nor \mu^\centerdot-\nu^\centerdot\nor_\infty=\|\mu^\centerdot-\nu^\centerdot\|_{\cC_\var}$ for $\mu^\centerdot,\nu^\centerdot\in\cC\cP_0$.
\er
\subsection{Fractional Brownian motion and Girsanov's theorem}\label{Sec:fBM}
In this section, we recall the definition and basic properties of {\rm f}Bm and the related Girsanov theorem (see \cite{DU99, NO02}).

A $d$-dimensional Gaussian process $(W^H_t)_{t\geq 0}$ defined on some probability space $(\Omega,\sF, \mP)$ is called an {\rm f}Bm with Hurst parameter $H\in(0,1)$ if for any $0\le s\le t$,
$$
\mE(W^{H,i}_tW^{H,j}_s)=\tfrac12(t^{2H}+s^{2H}-|t-s|^{2H})\1_{i=j},\ \ i,j=1,\cdots,d.
$$
For fixed $r\ge0$, it is easy to see that for any $t,s\ge0$,
\begin{align}\label{0115:00}
\mE((W^{H,i}_{t+r}-W^{H,i}_{r})(W^{H,j}_{s+r}-W^{H,j}_{r}))=\mE(W^{H,i}_{t}W^{H,j}_{s}).
\end{align}
This means that $(W^{H}_{t+r}-W^{H}_{r})_{t\ge0}$ is another  standard {\rm f}Bm.
The value of $H$ tells the behavior of {\rm f}Bm:
when $H=1/2$, the process is exactly a standard $d$-dimensional Brownian motion;
when $H>1/2$,  the increments of the process are positively correlated;
when $H<1/2$,   the increments of the process are negatively correlated. 

In what follows, we fix $H\in(0,\frac12]$ and introduce two constants used below
$$
q_H:=\tfrac1{1-H},\ \ c_H:=\sqrt{2H/((1-2H)\cB(1-2H,H+\tfrac12))}\1_{H\in(0,\frac12)}+\1_{H=\frac12},
$$
where $\cB(\alpha,\beta)$ is the usual Beta function defined by
\begin{align}\label{Bet}
\cB(\alpha,\beta):=\int^1_0(1-s)^{\alpha-1}s^{\beta-1}\dif s,\ \ \alpha,\beta>0.
\end{align}
It is well-known that {\rm f}Bm has the following representation (cf. \cite[Corollary 3.1]{DU99}):
\begin{align}\label{1230:05}
W^H_t=\int_0^t K_H(t,s)\dif W_s,
\end{align}
where $W$ is a standard $d$-dimensional Brownian motion and $K_H$ is given by
$$
K_H(t,s)=c_H\left((t/s)^{H-\frac12}(t-s)^{H-\frac12}+\left(\tfrac12-H\right)s^{\frac12-H}\int_s^tr^{H-\frac32}(r-s)^{H-\frac12}\dif r\right).
$$ 
Let $\mC_T:=C([0,T];\mR^d)$ be the space of all continuous functions from $[0,T]$ to $\mR^d$.  It is also well-known that 
there is a continuous functional $\Phi:\mC_T\to\mC_T$ so that (cf. \cite[Proposition A.1]{BG23})
\begin{align}\label{WW1}
W_t=\Phi(W^H_\cdot)(t),\ \ t\in[0,T].
\end{align}
{\bf Convention:} If there is no special declaration, we always write
$$
\sF_t:=\sigma\{W_s: s\leq t\}=\sigma\{W^H_s: s\leq t\}.
$$

From the very definition, it is easy to see that
\begin{align}\label{BC5}
K_H(t,s)\geq c_H(t-s)^{H-\frac12}\left((t/s)^{H-\frac12}+\left(\tfrac12-H\right)s^{\frac12-H}\int_s^tr^{H-\frac32}\dif r\right)=c_H(t-s)^{H-\frac12},
\end{align}
and by \cite[Theorem 3.2]{DU99},
$$
K_H(t,s)\leq c'_{H,T}(t-s)^{H-\frac12} s^{H-\frac12}.
$$
To state the Girsanov theorem for {\rm f}Bm, we introduce a function
$$
\widetilde K_H(t,s):=t^{H-\frac12}(t-s)^{-\frac12-H} s^{\frac12-H},\ \ 0\leq s<t.
$$
By the integration by parts and elementary calculus, one sees that
\begin{align}\label{DD1}
\int^t_s K_H(t,r)\widetilde K_H(r,s)\dif r\equiv1,\ \ 0\leq s<t.
\end{align}
For given $q\in[1,\infty)$, let $\mH^q_T$ be the space of all absolutely continuous function $f:[0,T]\to\mR^d$ with $f(0)=0$ and $\dot f\in L^q([0,T];\mR^d)=:\mL^q_T$,
which is a Banach space under the norm 
$$
\|f\|_{\mH^q_T}:=\|\dot f\|_{\mL^q_T}.
$$
Now, for any function $f\in C^1([0,T];\mR^d)$, we define
\begin{align*}
\widetilde{\bf K}_Hf(t):=\int_0^t\widetilde K_H(t,s)\dot f(s)\dif s.
\end{align*}
\bl
The operator $\widetilde{\bf K}_H$ can be extended to a bounded linear operator
from $\mH^{q_{\H}}_T$ to $\mL^2_T$, and there is a constant $C=C(H)>0$ such that for all $f\in\mH^{q_{\H}}_T$,
\begin{align}\label{AS2}
\|\widetilde{\bf K}_Hf\|_{\mL^2_T}\leq C\|f\|_{\mH^{q_{\H}}_T}
\end{align}
and
\begin{align}\label{AS3}
f(t)=\int^t_0K_H(t,s)\widetilde{\bf K}_H f(s)\dif s.
\end{align}
\el 
\begin{proof}
Estimate \eqref{AS2} follows by $\widetilde K_H(t,s)\leq(t-s)^{-\frac12-H}$ and Hard-Littlewood's inequality \cite[Theorem 1.7]{BCD}. Equality  \eqref{AS3} follows by Fubini's theorem and \eqref{DD1}.
\end{proof}

We have the following Girsanov theorem (see \cite[Theorem 4.9]{DU99}).
\bt\label{B:thm00}
Recall $q_H=\frac1{1-H}$.
Let $h(\cdot,\omega)\in \mH^{q_{\H}}_T$  be an $\sF_s$-adapted process satisfying
\begin{align}\label{WW2}
\mE\exp\left(\|h\|^2_{\mH^{q_{\H}}_T}\right)<\infty.
\end{align}
Then $\widetilde{W}^H_t:=W^H_t+h(t)$ is a new {\rm f}Bm with Hurst parameter $H$ under the new probability measure $\mQ:=Z_T\mP$ with
\begin{align*}
Z_T:=\exp\left(-\int_0^T(\widetilde{\bf K}_H h)(s)\dif W_s-\frac12 \|\widetilde{\bf K}_Hh\|^2_{\mL^2_T}\right),
\end{align*}
where $W$ defined by \eqref{WW1} is a $d$-dimensional standard Brownian motion.
\et
\begin{proof}
By \eqref{AS3}, we have
$$
\widetilde{W}^H_t=W^H_t+h(t)=\int^t_0 K_H(t,s)(\dif W_s+\widetilde{\bf K}_H h(s)\dif s).
$$
By \eqref{WW2} and Novikov's criterion, $\mE Z_T=1$.
Thus by the classical Girsanov theorem, $\widetilde W_t:=W_t+\int^t_0\widetilde{\bf K}_H h(s)\dif s$ is still a standard Brownian motion under $\mQ$, and therefore, $\widetilde{W}^H_t=\int^t_0K_H(t,s)\dif\widetilde W_s$ is an {\rm f}Bm under $\mQ$.
\end{proof}

Now we prove the following basic estimate. The new point is that  we are using the localized $L^p$-space.
%A general version can be found in \cite[Lemma 3.10]{BG23}. 
\bl\label{Le22}
Let $H\in(0,\frac12]$. For any $p\in[1,\infty]$ and $j\in\mN_0$, there is a constant $C=C(j,p,d,H)>0$ such that for all $0\le s\le t$ and $f\in \tilde{\mL}^p$,
 \begin{align*}
 | \mE^{\sF_s}(\nabla^j f(W_t^H))|\lesssim_C (t-s)^{%\frac{d}{2q}
-\frac{Hd}{p}-jH}%s^{-\frac{d}{2q}}
\nor f\nor_p.
\end{align*}
 \el
 \begin{proof}
 Note that by the representation \eqref{1230:05},
 $$
 \mE^{\sF_s}(W^H_t)=\int^s_0 K_H(t,r)\dif W_r
 $$
and
$$
W^H_{s,t}:=W^H_t-\mE^{\sF_s}(W^H_t)=\int^t_sK_H(t,r)\dif W_r.
$$
Clearly, $W^H_{s,t}$ is independent of $\sF_s$ and
$$
 \mE^{\sF_s}(W^H_t)\sim N(0, \sigma^H_{s,t}),\ \ W^H_{s,t}\sim N(0, \lambda^H_{s,t}),\ \ 
$$
where
$$
\sigma^H_{s,t}:=\int^s_0|K_H(t,r)|^2\dif r\stackrel{\eqref{BC5}}{\geq} c_H^2\int^s_0|t-r|^{2H-1}\dif r,
$$
and
\begin{align}\label{Lam1}
\lambda^H_{s,t}:=\int^t_s|K_H(t,r)|^2\dif r\stackrel{\eqref{BC5}}{\geq} c_H^2\int^t_s|t-r|^{2H-1}\dif r=\frac{c^2_H}{2H}(t-s)^{2H}.
\end{align}
By the independence of $W^H_{s,t}$ and $\sF_s$, we have
$$
\mE^{\sF_s}[\nabla^j f(W^H_t)]=\mE^{\sF_s}[\nabla^j f(W^H_{s,t}+\mE^{\sF_s}(W^H_t))]= F^{(j)}_{s,t}(\mE^{\sF_s}(W^H_t)),
$$
where
$$
F^{(j)}_{s,t}(y):=\mE [\nabla^j f(W^H_{s,t}+y)].
$$
By Lemma \ref{LeA1}, we have
$$
|\mE^{\sF_s}[\nabla^j f(W^H_t)]|\leq\|F^{(j)}_{s,t}\|_\infty\lesssim (\lambda^H_{s,t})^{(d/p-j)/2}\nor f\nor_p.
$$
Combining the above calculations, we obtain the desired estimate.
 \end{proof}

Below for simplicity of notations, we always write
$$
\mL^q_T\tilde\mL^p:=L^q([0,T]; \tilde\mL^p).
$$
As a result of the previous estimate, we have the following Krylov-type estimate.
\bl
For any $p,q\in[1,\infty]$ with $\alpha:=1-(1/q+Hd/p)>0$, there is a constant $C_0=C_0(d,p,q,H)>0$ such that for all $f\in \mL^q_T\tilde{\mL}^p$, $k\in\mN_0$ and $0\le s< t$,
\begin{align}\label{0202:00}
\mE^{\sF_s}\left(\int_s^tf(r,W_r^H)(t-r)^{k\alpha}\dif r\right)\le C_0k^{-\alpha}(t-s)^{(k+1)\alpha}\nor f\nor_{\mL^q_T\tilde{\mL}^p}.
\end{align}

\el
\begin{proof}
Since $\alpha=1-(1/q+Hd/p)>0$, we must have $q>1$ and $q'Hd/p<1$, where $q'=q/(q-1)$.
By Lemma \ref{Le22} and H\"older's inequality, one sees that
\begin{align*}
\sI&:=\mE^{\sF_s}\left(\int_s^t f(r,W_r^H)(t-r)^{k\alpha}\dif r\right)\\
&\lesssim_C \int_s^t\nor f(r)\nor_p(r-s)^{-\frac{Hd}{p}}(t-r)^{k\alpha}\dif r\\
&\leq \nor f\nor_{\mL^q_T\tilde{\mL}^p}\left(\int_s^t(r-s)^{-q'\frac{Hd}{p}}(t-r)^{q'k\alpha}\dif r\right)^{1/q'}.
\end{align*}
Let $\cB$ be the Beta function defined in \eqref{Bet}.
By a change of variable and Lemma \ref{Beta} , we get
\begin{align*}
\sI&\le C\nor f\nor_{\mL^q_T\tilde{\mL}^p}(t-s)^{(k+1)\alpha}\cB\left(1-q'\tfrac{Hd}{p},q'k\alpha+1\right)^{1/q'}\\
&\le C(p,q,H,d)\nor f\nor_{\mL^q_T\tilde{\mL}^p}(t-s)^{(k+1)\alpha}k^{-\frac1{q'}+\frac{Hd}{p}}.
\end{align*}
This completes the proof.
\end{proof}
Then we have the following moment estimate.
\bl
For any $p,q\in[1,\infty]$ with $\alpha:=1-(1/q+Hd/p)>0$, there is a constant $C_1=C_1(d,H,p,q)>0$ such that for all $f\in \mL^q_T\tilde{\mL}^p$ and $m\in\mN$,
\begin{align}\label{0202:02}
\left\| \int_0^tf(s,W_s^H)\dif s\right\|_{L^m(\Omega)}\le C_1t^{\alpha}m^{1-\alpha}\nor f\nor_{\mL^q_T\tilde{\mL}^p},\ \forall t\in(0,T]. 
\end{align}
\el
\begin{proof}
Without loss of generality, we assume that $f\geq 0$. For simplicity of notation, we write
$$
h(t):=f(t, W^H_t).
$$
By the symmetric of integral and \eqref{0202:00}, we have
\begin{align*}
\mE \left| \int_0^t h(s)\dif s\right|^{m}&=m!\mE \int_0^t h(s_1)\int_{s_1}^t h(s_2)\cdots \int_{s_{m-1}}^th(s_m)\dif s_m\cdots \dif s_2 \dif s_1\\
&=m!\mE \int_0^t h(s_1)\int_{s_1}^t h(s_2)\cdots \left[\mE^{\sF_{s_{m-1}}}\int_{s_{m-1}}^th(s_m)\dif s_m\right]\cdots\dif s_2 \dif s_1\\
&\leq C_0 m!\nor f\nor_{\mL^q_T\tilde{\mL}^p}\mE \int_0^t h(s_1)\cdot\cdot\cdot \int_{s_{m-2}}^t(t-s_{m-1})^{\alpha} h(s_{m-1})\dif s_{m-1}\cdot\cdot\cdot \dif s_1.
\end{align*}
Then, by \eqref{0202:00} and induction, we have for any $k=1,...,m-1$,
\begin{align*}
\mE \left| \int_0^t h(s)\dif s\right|^m&\le C_0^k m! ((k-1)!)^{-\alpha}\nor f\nor_{\mL^q_T\tilde{\mL}^p}^k\mE \int_0^t h(s_1)\cdot\cdot\cdot\\
&\qquad\qquad\times \mE^{\sF_{s_{m-k-1}}}\int_{s_{m-k-1}}^t(t-s_{m-k})^{k\alpha}h(s_{m-k})\dif s_{m-k}\cdot\cdot\cdot \dif s_1\\
&\le C_0^m (m!)^{1-\alpha} t^{m\alpha}\nor f\nor_{\mL^q_T\tilde{\mL}^p}^m.
\end{align*}
Finally, by Stirling's formula, we get
\begin{align*}
\mE \left| \int_0^t h(s)\dif s\right|^m\le C_1^m m^{m(1-\alpha)}t^{m\alpha} \nor f\nor_{\mL^q_T\tilde{\mL}^p}^m.
\end{align*}
The proof is complete.
\end{proof}

Now we can show the following important Khasminskii's type estimate.
\bt\label{Le24}
Let $q_1,p_1\in[\frac1{1-H},\infty]$ with $\alpha:=\frac12-(\frac1{q_1}+\frac {Hd}{p_1})>0$. 
Then for any $\lambda>0$, 
there is a constant $C=C(\lambda,p_1,q_1,d,H)>0$  such that for all  $b\in\mL^{q_1}_T\tilde\mL^{p_1}$,
\begin{align}\label{0624:00}
\mE\exp\left\{\lambda\|\widetilde{\bf K}_H\sI_b\|^2_{\mL^2_T}\right\}\leq \exp\left\{C\left(1+\kappa_b^{2/\alpha}\right)\right\},
\end{align}
where $\sI_b(t):=\int^t_0 b(s, W^H_s)\dif s$ and $\kappa_b:=\nor b\nor_{\mL^{q_1}_T\tilde\mL^{p_1}}$.
\et
\begin{proof}
By \eqref{AS2}, we have
\begin{align}\label{AS87}
\mE\exp\left\{\lambda\|\widetilde{\bf K}_H\sI_b\|^2_{\mL^2_T}\right\}
\leq\mE\exp\left\{C_0\lambda\|\sI_b\|^2_{\mH^{q_{\H}}_T}\right\}=\sum_{m=0}^\infty\frac{(C_0\lambda)^m}{m!} \mE\|\sI_b\|^{2m}_{\mH^{q_{\H}}_T}.
\end{align}
Observing that
\begin{align*}
\mE\|\sI_b\|^{2m}_{\mH^{q_{\H}}_T}=\mE\left(\int_0^T |b(s,W_s^H)|^{q_{\H}}\dif s\right)^{2m/q_{\H}}\le \left[\mE\left(\int_0^T |b(s,W_s^H)|^{q_{\H}}\dif s\right)^{2m}\right]^{1/q_{\H}}, 
\end{align*}
and $\alpha:=\frac12-(\frac1{q_1}+\frac {Hd}{p_1})>0$ and $q_1,p_1\in [q_{H},\infty]$, by \eqref{0202:02} with $(p,q)=(\frac {p_1}{q_{\H}},\frac {q_1}{q_{\H}})$, we have
\begin{align*}
\mE\|\sI_b\|^{2m}_{\mH^{q_{\H}}_T}&\le \left[(C_1T^{\alpha})^{2m}(2m)^{m(1-2\alpha)q_{\H}}\nor |b|^{q_{\H}}\nor_{\mL^{q_1/q_{\H}}_T\tilde\mL^{p_1/q_{\H}}}^{2m} \right]^{1/q_{\H}}\\
&=C_2^{m} m^{m(1-2\alpha)}   \nor b\nor_{\mL^{q_1}_T\tilde\mL^{p_1}}^{2m}
=C_2^{m} m^{m(1-2\alpha)}\kappa_b^{2m}.
\end{align*}
Substituting this into \eqref{AS87} and by Stirling's formula, we have
\begin{align*}
\mE\exp\left\{\lambda\|\widetilde{\bf K}_H\sI_b\|^2_{\mL^2_T}\right\}
&\le \sum_{m=0}^\infty\frac{(C_0\lambda C_2)^{m} m^{m(1-2\alpha)}\kappa_b^{2m}}{m!}\le \sum_{m=0}^\infty\frac{C_3^m \kappa_b^{2m}}{(m!)^{2\alpha}}.
\end{align*}
The proof is complete by $\alpha>0$ and Lemma \ref{Seri}.
\end{proof}

{
\br
Similar estimates to those in \eqref{0202:02} and \eqref{0624:00} have been established for VMO processes, as demonstrated in \cite[Corollary 3.5]{Le22}.
\er
}

\section{Well-posedness of DFSDE: Regular coefficients case}
In this section, we show the strong well-posedness (existence of strong solution and pathwise uniqueness) of the DFSDE \eqref{DFSDE} and prove Theorem \ref{Thm21}.

We first prepare the following standard result for later use.
\bl\label{Le31}
Let $\mD_T\ni(s,t)\mapsto\mu_{s,t}\in\cC\cP_1$ be a measurable function with
$$
\gamma^\mu_T:=\sup_{0\leq s\leq t\leq T}\|\mu_{s,t}\|_{\cC\cP_1}<\infty. 
$$
Consider the following classical SDE:
\begin{align}\label{AA01}
X^{x,\mu}_{s,t}=x+\int_s^t B(r,X^{x,\mu}_{s,r},\mu^{\centerdot}_{r,T},\mu^{\centerdot}_{s,r})\dif r+\int_s^t \Sigma(r,X^{x,\mu}_{s,r},\mu^{\centerdot}_{r,T},\mu^{\centerdot}_{s,r})\dif W_r.
\end{align}
Under {\bf(H$_0$)}, there is a unique strong solution to the above SDE with the estimate:
\begin{align}\label{AA1}
\mE|X^{x,\mu}_{s,t}|^2\leq \e^{\kappa_1 (t-s)}|x|^2+\frac{\kappa_0(\e^{\kappa_1(t-s)}-1)}{\kappa_1}
+\kappa_2\int^t_s\e^{\kappa_1(t-r)}(\|\mu_{r,T}\|^2_{\cC\cP_1}+\|\mu_{s,r}\|^2_{\cC\cP_1})\dif r.
\end{align}
Moreover, let $\nu_{s,t}$ be another $\cC\cP_1$-valued function with $\gamma^\nu_T:=\sup_{0\leq s\leq t\leq T}\|\nu_{s,t}\|_{\cC\cP_1}<\infty.$ Then
there is a constant $C_T=C_T(\kappa_i,\gamma^{\mu}_T,\gamma^\nu_T)>0$ such that for all $(s,t,x)\in\mD_T\times\mR^d$,
\begin{align}\label{AA2}
\frac{\mE|X^{x,\mu}_{s,t}-X^{x,\nu}_{s,t}|^2}{1+|x|^2}\leq C_T\int^t_s\left({\rm d}^2_{\cC\cP_1}(\mu_{r,T},\nu_{r,T})
+{\rm d}^2_{\cC\cP_1}(\mu_{s,r},\nu_{s,r})\right)\dif r.
\end{align}
\el
\begin{proof}
Under {\bf(H$_0$)}, the strong well-posedness  to SDE \eqref{AA01} are well-known (see \cite[]{Ro}).
We only show the estimates \eqref{AA1} and \eqref{AA2}. Fix $s\in[0,T)$.
By It\^o's formula and \eqref{AA05}, we have
\begin{align*}
\dif_t(\e^{-\kappa_1 (t-s)}|X^{x,\mu}_{s,t}|^2)
&=\e^{-\kappa_1(t-s)}\Big[\big(\<X^{x,\mu}_{s,t}, B(t,X^{x}_{s,t},\mu^\centerdot_{t,T},\mu^\centerdot_{s,t})\>+2\|\Sigma(t,X^{x,\mu}_{s,t},\mu^\centerdot_{t,T},\mu^\centerdot_{s,t})\|^2_{\rm HS}\\
&\qquad\qquad\qquad-\kappa_1|X^{x,\mu}_{s,t}|^2\big)\dif t+\<X^{x,\mu}_{s,t}, \Sigma(t,X^{x,\mu}_{s,t},\mu^{\centerdot}_{t,T},\mu^{\centerdot}_{s,t})\dif W_t\>\Big]\\
&\leq
\e^{-\kappa_1(t-s)}\big(\kappa_0+\kappa_2(\|\mu_{t,T}\|^2_{\cC\cP_1}+\|\mu_{s,t}\|^2_{\cC\cP_1})\big)\dif t+\dif M_t,
\end{align*}
where 
$$
t\mapsto M_t:=\int^t_0\e^{-\kappa_1(r-s)}\<X^{x,\mu}_{s,r}, \Sigma(r,X^{x,\mu}_{s,r},\mu^{\centerdot}_{r,T},\mu^{\centerdot}_{s,r})\dif W_r\>
$$ 
is a continuous local martingale.
By a standard stopping time technique and integrating both sides with respect to the time variable $t$ from $s$ to $t$, we derive that
$$
\e^{-\kappa_1 (t-s)}\mE|X^{x,\mu}_{s,t}|^2\leq |x|^2+\int^t_s\e^{-\kappa_1(r-s)}(\kappa_0+\kappa_2(\|\mu_{r,T}\|^2_{\cC\cP_1}+\|\mu_{s,r}\|^2_{\cC\cP_1}))\dif r.
$$
From this we get \eqref{AA1}.
For \eqref{AA2}, by It\^o's formula again and \eqref{AA06}, we have
\begin{align*}
\mE|X^{x,\mu}_{s,t}-X^{x,\nu}_{s,t}|^2\leq &\int^t_s\Big(\kappa_3\mE|X^{x,\mu}_{s,r}-X^{x,\nu}_{s,r}|^2
+\kappa_4(1+\mE|X^{x,\mu}_{s,r}|^2+\mE|X^{x,\nu}_{s,r}|^2)\\
&\qquad\times
\left({\rm d}^2_{\cC\cP_1}(\mu_{r,T},\nu_{r,T})
+{\rm d}^2_{\cC\cP_1}(\mu_{s,r},\nu_{s,r})\right)\Big)\dif r.
\end{align*}
By Gronwall's inequality, we get for all $t\in[s,T]$,
\begin{align}\label{AA4}
\mE|X^{x,\mu}_{s,t}-X^{x,\nu}_{s,t}|^2\lesssim\int^t_s
(1+\mE|X^{x,\mu}_{s,r}|^2+\mE|X^{x,\nu}_{s,r}|^2)
\left({\rm d}^2_{\cC\cP_1}(\mu_{r,T},\nu_{r,T})
+{\rm d}^2_{\cC\cP_1}(\mu_{s,r},\nu_{s,r})\right)\Big)\dif r.
\end{align}
Note that by \eqref{AA1},
$$
\sup_{t\in[s,T]}\mE|X^{x,\mu}_{s,t}|^2\leq C(T,\kappa_0,\kappa_1,\kappa_2)(1+|x|^2)(1+\gamma^\mu_T).
$$
Substituting this into \eqref{AA4}, we obtain the desired estimate.
\end{proof}

Now we can give 
\begin{proof}[Proof of Theorem \ref{Thm21}]
We use the method of freezing the distribution. 
Let $\mu^{x,0}_{s,t}:=\delta_x$ for all $(s,t,x)\in \mD_T\times\mR^d$.
For $n\in\mN$, by Lemma \ref{Le31}, we can recursively define the approximation solution $X^{x,n}_{s,t}$ by
\begin{align}\label{SD3}
X^{x,n+1}_{s,t}=x+\int_s^t B(r,X^{x,n+1}_{s,r},\mu^{\centerdot,n}_{r,T},\mu^{\centerdot,n}_{s,r})\dif r+\int_s^t \Sigma(r,X^{x,n+1}_{s,r},\mu^{\centerdot,n}_{r,T},\mu^{\centerdot,n}_{s,r})\dif W_r.
\end{align}
By \eqref{AA1} we have
\begin{align}\label{AA8}
\mE|X^{x;n+1}_{s,t}|^2\leq \e^{\kappa_1 (t-s)}|x|^2+\frac{\kappa_0(\e^{\kappa_1(t-s)}-1)}{\kappa_1}
+\kappa_2\int^t_s\e^{\kappa_1(t-r)}(\|\mu^{\centerdot, n}_{r,T}\|^2_{\cC\cP_1}+\|\mu^{\centerdot, n}_{s,r}\|^2_{\cC\cP_1})\dif r.
\end{align}
Noting that
\begin{align}\label{AA9}
\|\mu^{\centerdot,n}_{s,t}\|^2_{\cC\cP_1}=\sup_{x\in\mR^d }\left(\frac{\mE|X^{x,n}_{s,t}|}{1+|x|}\right)^2
\leq \sup_{x\in\mR^d }\frac{\mE|X^{x,n}_{s,t}|^2}{(1+|x|)^2}=:f_n(s,t),
\end{align}
we have
\begin{align*}
f_{n+1}(s,t)\leq \e^{\kappa_1 (t-s)}+\frac{\kappa_0(\e^{\kappa_1(t-s)}-1)}{\kappa_1}
+\kappa_2\int^t_s\e^{\kappa_1(t-r)}(f_{n}(r,T)+f_{n}(s,r))\dif r.
\end{align*}
For $m\in\mN$, if we let 
$$
F_m(s,t):=\sup_{n\leq m+1}f_n(s,t),
$$
then for each $0\leq s\leq t\leq T$,
\begin{align*}
F_m(s,t)&\lesssim 1+ \int_s^t [F_m(s,r)+F_m(r,T)]\dif r
\leq 1+\int_s^T F_m(r,T)\dif r+\int_s^t F_m(s,r)\dif r. 
\end{align*}
By \eqref{AA9} and Gronwall's inequality (see Lemma \ref{lemA03}), we have
\begin{align}\label{AA08}
\sup_m\sup_{(s,t)\in\mD_T}\|\mu^{\centerdot,m}_{s,t}\|^2_{\cC\cP_1}\leq
\sup_{m\in\mN}\sup_{(s,t)\in\mD_T}F_m(s,t)<\infty.
\end{align}
Next, we show that the sequence $\{X^{x,n}_{s,\cdot}\}_{n=1}^\infty$ is a Cauchy sequence in {a suitable norm}. 
By \eqref{AA2}, we have for any $n,m\in\mN$, $x\in\mR^d$ and $0\le s\le t\le T$,
\begin{align}\label{SD1}
&\frac{\mE|X^{x,n+1}_{s,t}-X^{x,m+1}_{s,t}|^2}{1+|x|^2}
\lesssim \int_s^t  \left({\rm d}^2_{\cC\cP_1}(\mu^{\centerdot,n}_{r,T},\mu^{\centerdot,m}_{r,T})
+{\rm d}^2_{\cC\cP_1}(\mu^{\centerdot,n}_{s,r},\mu^{\centerdot,m}_{s,r})\right)\dif r.
\end{align}
Noting that
$$
H^{n,m}_{s,t}:={\rm d}^2_{\cC\cP_1}(\mu^{\centerdot,n}_{s,t},\mu^{\centerdot,m}_{s,t})\leq \sup_{x\in\mR^d}\frac{\mE|X^{x,n}_{s,t}-X^{x,m}_{s,t}|^2}{(1+|x|)^2},
$$
by \eqref{SD1} we have
$$
H^{n+1,m+1}_{s,t}\lesssim \int_s^t  \left[H^{n,m}_{s,r}+H^{n,m}_{r,T}\right]\dif r.
$$
Thus, by \eqref{AA08} and Fatou's lemma, we derive that for all $t\in[s,T]$,
$$
\varlimsup_{n,m\to\infty}H^{n+1,m+1}_{s,t}\lesssim \int_s^t  \left[\varlimsup_{n,m\to\infty}H^{n,m}_{s,r}+\varlimsup_{n,m\to\infty}H^{n,m}_{r,T}\right]\dif r,
$$
which implies by Gronwall's inequality that
\begin{align}\label{SD2}
\varlimsup_{n,m\to\infty}H^{n,m}_{s,t}=\varlimsup_{n,m\to\infty}\sup_{x\in\mR^d}\left(\frac{\cW_1(\mu^{x,n}_{s,t},\mu^{x,m}_{s,t})}{1+|x|}\right)^2=0.
\end{align}
Substituting this into \eqref{SD1}, we obtain
$$
\varlimsup_{n,m\to\infty}\sup_{(t,x)\in[s,T]\times\mR^d}\frac{\mE|X^{x,n+1}_{s,t}-X^{x,m+1}_{s,t}|^2}{1+|x|^2}=0.
$$
In particular, for each fixed $(s,x)\in[0,T)\times\mR^d$, there is an adapted process $\{X^{x}_{s,t}\}_{t\in[s,T]}$ so that
$$
\lim_{n\to\infty}\sup_{t\in[s,T]}\mE|X^{x,n}_{s,t}-X^{x}_{s,t}|^2=0,
$$
and by \eqref{SD2}, for each $s\leq t$, there is a family of probability measures $(\mu^x_{s,t})_{x\in\mR^d}\in\cC\cP_1$ so that
$$
\lim_{n\to\infty}\sup_{x\in\mR^d}\frac{\cW_1(\mu^{x,n}_{s,t},\mu^x_{s,t})}{1+|x|}
=\lim_{n\to\infty}{\rm d}_{\cC\cP_1}(\mu^{\centerdot,n}_{s,t},\mu^{\centerdot}_{s,t})=0.
$$
Since $\mu^{x,n}_{s,t}=\mP\circ (X^{x,n}_{s,t})^{-1}$, we have
$$
\mu^x_{s,t}=\mP\circ (X^{x}_{s,t})^{-1},\ \ \forall x\in\mR^d.
$$
Finally, by the continuity of $(x,\mu^\centerdot,\nu^\centerdot )\mapsto B(t,x,\mu^\centerdot,\nu^\centerdot)$ and taking limits for equation \eqref{SD3}, one sees that
for each $(s,t,x)\in\mD_T\times\mR^d$,
$$
X^{x}_{s,t}=x+\int_s^t B(r,X^{x}_{s,r},\mu^{\centerdot}_{r,T},\mu^{\centerdot}_{s,r})\dif r+\int_s^t \Sigma(r,X^{x}_{s,r},\mu^{\centerdot}_{r,T},\mu^{\centerdot}_{s,r})\dif W_r.
$$
Thus we obtain the existence. The pathwise uniqueness is derived by the same argument.
Moreover, by \eqref{AA8}, we have
\begin{align}\label{AA81}
\mE|X^{x}_{s,t}|^2\leq \e^{\kappa_1 (t-s)}|x|^2+\frac{\kappa_0(\e^{\kappa_1(t-s)}-1)}{\kappa_1}
+\kappa_2\int^t_s\e^{\kappa_1(t-r)}(\|\mu^{\centerdot}_{r,T}\|^2_{\cC\cP_1}+\|\mu^{\centerdot}_{s,r}\|^2_{\cC\cP_1})\dif r,
\end{align}
which implies by \eqref{AA9} that
\begin{align*}
\|\mu^{\centerdot}_{s,t}\|^2_{\cC\cP_1}\leq \e^{\kappa_1 (t-s)}+\frac{\kappa_0(\e^{\kappa_1(t-s)}-1)}{\kappa_1}
+\kappa_2\int^t_s\e^{\kappa_1(t-r)}(\|\mu^{\centerdot}_{r,T}\|^2_{\cC\cP_1}+\|\mu^{\centerdot}_{s,r}\|^2_{\cC\cP_1})\dif r.
\end{align*}
By Gronwall's inequality  (see Lemma \ref{lemA03}), we have
\begin{align}\label{AA82}
\|\mu^{\centerdot}_{s,t}\|^2_{\cC\cP_1}\leq C_T.
\end{align}
If $\kappa_1<0$ and $2\kappa_2<|\kappa_1|$, then
$$
\|\mu^{\centerdot}_{s,t}\|^2_{\cC\cP_1}\leq 1+\frac{\kappa_0}{|\kappa_1|}
+\frac{2\kappa_2}{|\kappa_1|}\sup_{0\leq s\leq t\leq T}\|\mu^{\centerdot}_{s,t}\|^2_{\cC\cP_1},
$$
which implies that
\begin{align}\label{AA83}
\sup_{0\leq s\leq t\leq T}\|\mu^{\centerdot}_{s,t}\|^2_{\cC\cP_1}\leq (1+\frac{\kappa_0}{|\kappa_1|})/(1-\frac{2\kappa_2}{|\kappa_1|}).
\end{align}
Substituting \eqref{AA82} and \eqref{AA83} into \eqref{AA81}, we get the desired estimates.
\end{proof}

Now we provide simple examples to illustrate the assumptions.
\bx
Let $b_1, b_2:\mR^d\to\mR^d$ satisfy
\begin{align*}
(1+|x|)\, b_1,\nabla b_1\in L^1\ \mbox{and } \nabla b_2\in L^\infty.
\end{align*}
Let $\varphi_1,\varphi_2:\mR^d\to\mR$ be two Borel measurable functions with
\begin{align*}
\varphi_1,\nabla \varphi_1\in L^\infty\ \text{and}\ (1+|x|)\,\varphi_2\in L^1.
\end{align*}
For $\mu,\nu\in\cC\cP_1$, we introduce
\begin{align*}
B(x,\mu^\centerdot,\nu^\centerdot):=\int_{\mR^{d}}b_1(x-y)\mu^y(\varphi_1)\dif y+\int_{\mR^{d}}(b_2*\nu^z)(x)\varphi_2(z)\dif z=:B_1(x,\mu^\centerdot)
+B_2(x,\nu^\centerdot).
\end{align*}
Then one sees that  \eqref{AA06} and \eqref{AA05} hold. Indeed,  for $B_1(x,\mu^\centerdot)$, we have
\begin{align*}
|B_1(x_1,\mu^\centerdot_1)-B_1(x_2,\mu^\centerdot_2)|
&\leq \int_{\mR^{d}}|b_1(x_1-y)-b_1(x_2-y)|\,|\mu_1^y(\varphi_1)|\dif y\\
&\quad+\int_{\mR^{d}}|b_1(x_2-y)|\,|(\mu^y_1-\mu^y_2)(\varphi_1)|\dif y\\
&\leq \|\varphi_1\|_\infty|x_1-x_2|\int_{\mR^{d}}\int^1_0|\nabla b_1(x_1-y+\theta(x_2-x_1))|\dif\theta\dif y\\
&\quad+\left(\int_{\mR^{d}}|b_1(x_2-y)|(1+|y|)\dif y\right)\|\nabla\varphi_1\|_\infty{\rm d}_{\cC\cP_1}(\mu^\centerdot_1,\mu^\centerdot_2)\\
&\leq \|\varphi_1\|_\infty\|\nabla b_1\|_{1}|x_1-x_2|\\
&\quad+(|x_2|+1)\|(1+|\cdot|)b_1\|_{1}\|\nabla\varphi_1\|_\infty{\rm d}_{\cC\cP_1}(\mu^\centerdot_1,\mu^\centerdot_2).
\end{align*}
For $B_2(x,\nu^\centerdot)$, we have
\begin{align*}
|B_2(x_1,\nu^\centerdot_1)-B_2(x_2,\nu^\centerdot_2)|
&\leq \int_{\mR^{d}}\int_{\mR^d}|b_2(x_1-y)-b_2(x_2-y)|\nu^z_1(\dif y)|\varphi_2(z)|\dif z\\
&\quad+\int_{\mR^{d}}|b_2*\nu^z_1-b_2*\nu^z_2|(x_2)|\varphi_2(z)|\dif z\\
&\leq |x_1-x_2|\|\nabla b_2\|_\infty\|\varphi_2\|_{1}
+\|\nabla b_2\|_\infty\int_{\mR^{d}}\cW_1(\nu^z_1,\nu^z_2)|\varphi_2(z)|\dif z\\
&\leq \|\nabla b_2\|_\infty\Big(\|\varphi_2\|_{1}|x_1-x_2|+\|(1+|\cdot|)\varphi_2 \|_{1}{\rm d}_{\cC\cP_1}(\nu^\centerdot_1,\nu^\centerdot_2)\Big).
\end{align*}
For \eqref{AA05}, for any $\eps\in(0,1)$, by Young's inequality we have
\begin{align*}
\left\<x, B(x,\mu^{\centerdot},\nu^\centerdot)\right\>&\le |x|\|B(\cdot,\mu^{\centerdot},\nu^\centerdot)\|_\infty\\
&\le |x|\left(\|b_1\|_{1}\|\varphi_1\|_\infty+\int_{\mR^d}\|b_2*\nu^z\|_\infty |\varphi_2(z)|\dif z\right)\\
&\le |x|\left(\|b_1\|_{1}\|\varphi_1\|_\infty+\int_{\mR^d}\|\nabla b_2\|_\infty\|\nu^\centerdot\|_{\cC\cP_1}(1+|z|) |\varphi_2(z)|\dif z\right)\\
&\le \eps |x|^2+\left(\|b_1\|_1\|\varphi_1\|_\infty+\|\nabla b_2\|_\infty\|\nu^\centerdot\|_{\cC\cP_1}\|(1+|\cdot|) \varphi_2(\cdot)\|_{1}\right)^2/(4\eps).
\end{align*}
 Hence, \eqref{AA05} holds with 
$$
\kappa_1=\eps\quad \text{and} \quad\kappa_2=\|\nabla b_2\|^2_\infty\|(1+|\cdot|) \varphi_2(\cdot)\|^2_{1}/(4\eps).
$$
In particular, for any $\lambda>\eps+\|\nabla b_2\|^2_\infty\|(1+|\cdot|) \varphi_2(\cdot)\|^2_{1}/(2\eps)$, by \eqref{AA106}, we have uniform moment estimate  in time  
for the solution of \eqref{AA01} with diffusive coefficient $\mI$ and drift $B(x,\mu^{\centerdot},\nu^\centerdot)-\lambda x$.
\ex
\bx\label{Ex33}
Let $\sigma:[0,T]\times\mR^d\times\mR\to\mR$ satisfy
\begin{align*}
|\sigma(t,x,r)-\sigma(t,x',r')|\leq C(|x-x'|+|r-r'|).
\end{align*}
Let $\phi_\eps$ be a family of mollifiers.
For $\mu\in\cC\cP_1$, we introduce
\begin{align*}
\Sigma_\eps(t,x,\mu^\centerdot):=\sigma\left(t,x, \int_{\mR^d}\phi_\eps(x-y)\<\mu^y,\varphi\>\dif y\right).
\end{align*}
Then it is easy to see that  \eqref{AA06} and \eqref{AA05} hold. In particular, the following SDE admits a unique solution
$$
X^{x,\eps}_{s,t}=x+\int_s^t \Sigma_\eps(r,X^{x,\eps}_{s,r},\mu^{\centerdot,\eps}_{r,T})\dif W_r.
$$
An open question is whether we can take limits $\eps\to 0$ so that we can give a probability representation $u(s,x)=\mE\varphi(X^{x,0}_{s,T})$ for local quasi-linear PDE:
$$
\p_s u+\frac12\sum_{i,j,k}(\sigma_{ik}\sigma_{jk})(s,x,u)\p_i\p_j u=0,\ \ u(T)=\varphi,
$$
where $X^{x,0}_{s,T}$ solves the following nonlinear-SDE:
$$
X^{x,0}_{s,t}=x+\int_s^t \sigma(r,X^{x, 0}_{s,r},\mu^{x, 0}_{r,T}(\varphi))\dif W_r.
$$
We will study this problem in a future work.
\ex

\section{Well-posedness of DFSDE: Singular drift case}

In this section, we consider the DFSDE driven by fractional Brownian motion with a fixed value of $H\in(0,\frac12]$.
In Subsection \ref{41}, we focus on the well-posedness of SDEs driven by {\rm f}Bm using Girsanov's theorem. Specifically, we extend the results of Nualart-Ouknine \cite{NO02} to the case where the drift term 
 belongs to $\mL^q_T\tilde\mL^p$ with $p,q$ in the range of $[\frac1{1-H},\infty]$ and satisfying the condition $\frac1q+\frac{Hd}p<\frac12$. Notably, allowing $q$
 to be smaller than $2$ is crucial for applications to the 2D-Navier-Stokes equation associated with {\rm f}Bm.
 In Subsection \ref{42}, we establish the weak well-posedness for DFSDEs driven by {\rm f}Bm using the entropy method.
In Subsection  \ref{43}, we examine a backward DFSDE driven by Brownian motion by utilizing It\^o's formula and PDE's estimates. 
This analysis will be instrumental in demonstrating the well-posedness of the backward Navier-Stokes equation with $\mL^{1+}$-initial vorticity.

\subsection{SDEs driven by {\rm f}Bm}\label{41}

Let $\mC_T:=C([0,T];\mR^d)$ be the space of all continuous functions from $[0,T]$ to $\mR^d$, which is endowed with the uniform convergence topology.  
The canonical process on $\mC_T$ is defined by
$$
w_t(\omega):=\omega_t,\ \ \omega\in\mC_T.
$$
Let $\sB_t:=\sigma\{w_s: s\leq t\}$ be the natural filtration. Let $b:[0,T]\times \mC_T\to\mR^d$ be a $\sB_t$-progressively measurable vector field.
In this section we consider the following SDE:
\begin{align}\label{SDE3}
X_t=X_0+\int^t_0b(s,X_\cdot)\dif s+W^H_t,
\end{align}
where $W^H$ is an {\rm f}Bm with $H\in(0,\frac12]$.
To emphasize the dependence on $b$, we shall call SDE \eqref{SDE3} as SDE$_b$. We introduce the following definition of a weak solution to SDE$_b$.
\bd\label{Def3}
Let  $\nu\in\cP(\mR^d)$. We call a probability measure  $\bP\in\cP(\mC_T)$ a weak solution of SDE$_b$ starting from the initial distribution $\nu$ if $\bP\circ w_0^{-1}=\nu$ and
\begin{align}\label{SK2}
t\mapsto w_t-w_0-\int^t_0 b(s,w_\cdot)\dif s=:\sW^b_t
\end{align}
is an {\rm f}Bm with Hurst parameter $H$ under $\bP$. The set of all weak solutions of SDE$_b$ with initial distribution $\nu$ is denoted by $\cS(b,\nu)$.
We call the uniqueness in law holds for SDE$_b$ if any two $\bP_1,\bP_2\in\cS(b,\nu)$ are the same.
\ed

Recall that for two $\bP_1,\bP_2\in\cP(\mC_T)$, the relative entropy is defined by
$$
\cH(\bP_1|\bP_2):=
\left\{
\begin{aligned}
&\mE^{\bP_1}\ln(\dif \bP_1/\dif\bP_2),&\ \bP_1\ll\bP_2,\\
&\infty,&\mbox{otherwise.}
\end{aligned}
\right.
$$
By Csisz\'ar-Kullback-Pinsker's inequality (abbreviated as CKP's inequality) (see \cite[(22.25)]{Vi06}), we have
\begin{align}\label{CKP}
\|\bP_1-\bP_2\|_{\rm var}\leq \sqrt{2\cH(\bP_1|\bP_2)}.
\end{align}
We now prepare the following result about the relative entropy
(see Lacker \cite[Lemma 4.1]{La21} for a version of Brownian case). 
\bl\label{Le45}
Let $\nu\in\cP(\mR^d)$ and $b_i:[0,T]\times\mC_T\to\mR^d$, $i=1,2$ be two progressively measurable vector fields
and $\bP_i\in\cS(b_i,\nu)$.
Suppose that the  uniqueness in law holds for SDE$_{b_2}$ and
\begin{align}\label{In1}
\mE^{\bP_1}\exp\left\{\lambda\|\sI_{b_1-b_2}\|^2_{\mH^{q_{\H}}_T}\right\}<\infty,\ \ \forall\lambda>0,
\end{align}
where $\sI_{b_1-b_2}(t):=\int^t_0(b_1-b_2)(s,w_\cdot)\dif s$.
Then for some $C=C(H)>0$, it holds that
\begin{align*}
\cH(\bP_1|\bP_2)\lesssim_C\mE^{\bP_1}\left(\|\sI_{b_1-b_2}\|^2_{\mH^{q_{\H}}_T}\right).
\end{align*}
\el
\begin{proof}
For $i=1,2$, by definition, one has
$$
\sW^{b_i}_t:=w_t-w_0-\int^t_0 b_i(s,w_\cdot)\dif s\mbox{is an {\rm f}Bm with respect to $\bP_i$.}
$$
Let $W_t:=\Phi(\sW^{b_1}_\cdot)(t)$ (see \eqref{WW1}). Then $W$ is a standard Brownian motion under $\bP_1$.
Write
$$
Z_T:=\exp\left(-\int_0^T(\widetilde{\bf K}_H\sI_{b_1-b_2})(s)\dif W_s-\frac12 \|\widetilde{\bf K}_H\sI_{b_1-b_2}\|^2_{\mL^2_T}\right).
$$
By \eqref{AS2}, \eqref{In1}  and Novikov's criterion,  $\mE^{\bP_1}Z_T=1$.
Thus, by Girsanov's theorem (Theorem \ref{B:thm00}),
$$
t\mapsto \sW^{b^1}_t+\sI_{b_1-b_2}(t)=\sW^{b^2}_t
$$
 is still an {\rm f}Bm under $\bQ:=Z_T\bP_1$. 
Thus $\bQ\in\cS(b_2,\nu)$.
By the uniqueness in law of SDE$_{b_2}$, we have
$$
Z_{T}\bP_1=\bQ=\bP_2.
$$
Hence,
\begin{align}\label{1230:06}
\cH(\bP_1|\bP_2)=-\int_{\mC_T}\ln Z_{T}\dif \bP_1
=\frac12 \mE^{\bP_1}\left(\|\widetilde{\bf K}_H\sI_{b_1-b_2}\|^2_{\mL^2_T}\right)
\stackrel{\eqref{AS2}}{\lesssim}_C\mE^{\bP_1}\|\sI_{b_1-b_2}\|^2_{\mH^{q_{\H}}_T}.
\end{align}
Thus we complete the proof.
\end{proof}

By Theorem \ref{Le24} and Girsanov's theorem, it is by now standard to show the following result  (see \cite[Theorem 2]{NO02}).
\bt\label{Th25}
Suppose that for some $(p_1,q_1)\in[\frac1{1-H},\infty]^2$ with $\alpha:=\frac12-(\tfrac{1}{q_1}+\tfrac{Hd}{p_1})>0$,
$$
\kappa_b:=\nor b\nor_{\mL^{q_1}_T\tilde \mL^{p_1}}<\infty.
$$
Then for any $x\in\mR^d$, there is a unique weak solution $\bP_x\in\cS(b,\delta_x)$  in the class that
$$
\bP_x\left(\int^T_0|b(s,w_s)|^{q_{\H}}\dif s<\infty\right)=1. 
$$
Moreover, we have the following conclusions:
\begin{enumerate}[(i)]
\item For any $1<p\leq q\leq\infty$,
there is a constant $C=C(T,H,d,p_1,q_1,p,q)>0$ such that for all $f\in\tilde\mL^p$ and $t\in(0,T]$,
\begin{align}\label{0123:01}
\nor\mE^{\bP_\cdot} f(w_t)\nor_q\leq \exp\left\{C\left(1+\kappa_b^{2/\alpha}\right)\right\} t^{\frac{dH}{q}-\frac{dH}{p}}\nor f\nor_p.
\end{align}
\item For any $p,q\in(1,\infty]$  with $\beta:=1-(\frac1q+\frac{Hd}p)>0$, there is a  $C=C(T,H,d,p_1,q_1,\kappa_b,p,q)>0$ such that for all $t\in[0,T]$, $x\in\mR^d$ and $m\geq 1$,
\begin{align}\label{0123:41}
\left\|\int^t_0f(s,w_s)\dif s\right\|_{L^m(\mC_T;\bP_x)}\leq C t^{\beta}m^{1-\beta}\nor f\nor_{\mL^q_T\tilde{\mL}^p}.
\end{align}
\item For any $\gamma\in(0,1)$, $p\in[1,\infty]$ and $p_b,q_b\in(q_\H,\infty]$ satisfying 
$$
\beta_H:=1-q_H(\tfrac1q_b+\tfrac{Hd}{p_b})>0,
$$
there is a constant $C=C(\gamma,d,H,p,p_b,q_b,\kappa_b,T)>0$ such that for all $t\in[0,T]$, $f\in \bB^{\gamma}_p(\mR^d)$ and $x\in\mR^d$,
\begin{align}\label{0318:00}
\|\mE^{\bP_x}f(\cdot-w_t)-f(\cdot-x)\|_p\le C t^{H\gamma}\|f\|_{\bB^{\gamma}_p}+t^{\beta_H}\|f\|_{p}\nor b\nor^{q_\H}_{\mL^{q_b}_T\widetilde{\mL}^{p_b}},
\end{align}
where $\bB^{\gamma}_p(\mR^d)$ denotes the Sobolev space consisting of all functions $f$ with
\begin{align}\label{0320:00}
\|f\|_{\bB^{\gamma}_p}:=\sup_{h\ne0}\frac{\|f(\cdot+h)-f(\cdot)\|_{p}}{|h|^\gamma}+\|f\|_p<\infty.
\end{align}
\end{enumerate}
\et
\begin{proof}
{\bf (Existence)}
Let $x\in\mR^d$ and $W^H$ be an {\rm f}Bm over a probability space $(\Omega,\sF,\mP)$. Define
$$
\tilde W^H_t:=W^H_t-\int^t_0 b(s, W^H_s+x)\dif s=:W^H_t-\sI^x_b(t).
$$
Since $b\in \mL^{q_1}_T\tilde \mL^{p_1}$ and $\alpha=\frac12-(\tfrac{1}{q_1}+\tfrac{Hd}{p_1})>0$, by Theorem \ref{Le24}, for any $\lambda>0$, there is a
constant $C=C(\lambda,p_1,q_1,d,H)>0$ such that
\begin{align}\label{BC8}
\sup_{x\in\mR^d}\mE\exp\left\{\lambda\|\widetilde{\bf K}_H\sI^x_b\|^2_{\mL^2_T}\right\}\leq\exp\left\{C\left(1+\kappa_b^{2/\alpha}\right)\right\}.
\end{align}
By Theorem \ref{B:thm00}, $\tilde W^H$ is an {\rm f}Bm under $\mQ_x=Z^x_T\mP$, where
$$
Z^x_t=\exp\left(-\int_0^t(\widetilde{\bf K}_H \sI^x_b)(s)\dif W_s-\frac12 \|\widetilde{\bf K}_H\sI^x_b\|^2_{\mL^2_T}\right),
$$
is an exponential martingale. Here $W_t:=\Phi(W^H_\cdot)(t)$ (see \eqref{WW1}).
Now if we let $X^x_t:=W^H_t+x$, then 
$$
X^x_t=x+\int^t_0b(s, X^x_s)\dif s+\tilde W^H_t.
$$
In particular, $\bP_x:=\mQ_x\circ (X^x_\cdot)^{-1}\in\cS(b,\delta_x)$ is a solution of SDE$_b$. 

{\bf (Uniqueness)} 
For $i=1,2$, let $\bP_i\in\cS(b,\delta_x)$ so that $\sW^b$ is an {\rm f}Bm with Hurst parameter $H$ under $\bP_i$, and
$$
\bP_i(w_0=x)=1,\ \ \bP_i\left(\int^T_0|b(s,w_s)|^{q_{\H}}\dif s<\infty\right)=1. 
$$
Define a $\sB_t$-stopping time by
$$
\tau_n:=\inf\left\{t\in[0,T]: \int^t_0 |b(s,w_s)|^{q_{\H}}\dif s>n\right\}.
$$
Then $\lim_{n\to\infty}\bP_i(\tau_n=T)=1$. Let $\sI_b(t):=\int^t_0b(s,w_s)\dif s$. 
Note that
$$
\|\widetilde{\bf K}_H\sI_b(\cdot\wedge\tau_n)\|_{\mL^2_T}\leq C\| \sI_b(\cdot\wedge\tau_n)\|_{\mH^{q_{\H}}_T}\leq Cn^{1/q_{\H}}.
$$
By Girsanov's theorem,
\begin{align}\label{AQ1}
\sW^b_t+\sI_b(t\wedge\tau_n)=w_t-w_0-\int^t_{t\wedge\tau_n}b(s,w_s)\dif s
\end{align}
is still an {\rm f}Bm under the new probability $\bQ^n_i:=Z^n_T\bP_i$, where
$$
Z^n_T:=\exp\left(-\int_0^T(\widetilde{\bf K}_H \sI_b(\cdot\wedge\tau_n))(s)\dif W_s-\frac12 \|\widetilde{\bf K}_H\sI_b(\cdot\wedge\tau_n)\|^2_{\mL^2_T}\right).
$$
In particular, for any $t_1<t_2<\cdots<t_m\leq T$ and $\Gamma_i\in\sB(\mR^d)$, by \eqref{AQ1} we have
\begin{align*}
\bP_1(w_{t_1}\in\Gamma_1,\cdots, w_{t_m}\in\Gamma_m; \tau_n=T)
&=\int_{\mC_T}\1_{\{w_{t_1}\in\Gamma_1,\cdots, w_{t_m}\in\Gamma_m; \tau_n=T\}}/Z^n_T\bQ^n_1(\dif \omega)\\
&=\int_{\mC_T}\1_{\{w_{t_1}\in\Gamma_1,\cdots, w_{t_m}\in\Gamma_m; \tau_n=T\}}/Z^n_T\bQ^n_2(\dif \omega)\\
&=\bP_2(w_{t_1}\in\Gamma_1,\cdots, w_{t_m}\in\Gamma_m; \tau_n=T).
\end{align*}
Letting $n\to\infty$, we conclude the proof of uniqueness.

{\bf (Proofs of (i) and (ii))}
Let $1<p\leq q\leq\infty$ and $p_0\in(1,p)$. Set $\gamma:=\frac p{p_0}$ and $\frac1\gamma+\frac1{\gamma'}=1$. 
By H\"older's inequality, we have
\begin{align*}
\mE^{\bP_x}f(w_t)&=\mE^{\mQ_x} (f(X^x_t))=\mE^\mP (Z^x_T f(W^H_t+x))\leq(\mE^\mP (Z^x_T)^{\gamma'})^{1/\gamma'}(\mE^\mP|f(W^H_t+x)|^\gamma)^{1/\gamma}.
\end{align*}
Noting that by \eqref{BC8} and H\"older's inequality,
$$
\sup_{t\in[0,T]}\sup_{x\in\mR^d}\mE^\mP|Z^x_t|^{\gamma'}\leq\exp\left\{C\left(1+\kappa_b^{2/\alpha}\right)\right\}=:C_0,
$$
by Lemma \ref{LeA1}, we further have
\begin{align*}
\nor\mE^{\bP_\cdot}f(w_t)\nor_q&\lesssim_{C_0} \nor\mE^\mP|f(W^H_t+\cdot)|^\gamma\nor^{1/\gamma}_{q/\gamma}\lesssim
t^{\frac{dH}{q}-\frac{dH}{\gamma p_0}}\nor |f|^\gamma \nor^{1/\gamma}_{p_0}=t^{\frac{dH}{q}-\frac{dH}{p}} \nor f \nor_{p}.
\end{align*}
Thus we get \eqref{0123:01}. For \eqref{0123:41}, it is similar by \eqref{0202:02}  and H\"older's inequality.

{\bf (Proof of (iii)):}  Let $\bQ_x\in \cS(0,\delta_x)$ and note that
\begin{align*}
\|\mE^{\bP_x}f(\cdot-w_t)-f(\cdot-x)\|_p\le \|\mE^{\bQ_x}f(\cdot-w_t)-f(\cdot-x)\|_p+\|\mE^{\bP_x}f(\cdot-w_t)-\mE^{\bQ_x}f(\cdot-w_t)\|_p.
\end{align*}
Since $\bQ_x$ is the law of {\rm f}Bm starting from $x\in\mR^d$, we have 
\begin{align*}
&\|\mE^{\bQ_x}f(\cdot-w_t)-f(\cdot-x)\|_p\leq\left \|\int_{\mR^d}f(\cdot-x-z)p_t^H(z)\dif z-f(\cdot-x)\right\|_p\\
&\qquad\leq\int_{\mR^d}\|f(\cdot-z)-f(\cdot)\|_p p_t^H(z)\dif z\lesssim \|f\|_{\bB^{\beta}_p}\int_{\mR^d}|z|^\beta p_t^H(z)\dif z
\lesssim t^{\beta H}\|f\|_{\bB^{\beta}_p},
\end{align*}
where $p_t^H(z)=(2\pi \lambda^H_{0,t})^{-d/2}e^{|z|^2/\lambda^H_{0,t}}$ and $\lambda^H_{0,t}$ is defined in \eqref{Lam1}. 

Moreover, by \eqref{CKP}, \eqref{1230:06} and \eqref{0123:41}, we have
\begin{align*}
\|\mE^{\bP_x}f(\cdot-w_t)-\mE^{\bQ_x}f(\cdot-w_t)\|_p
&\le\|f\|_{p}\|\bP_x\circ w_t^{-1}-\bQ_x\circ w_t^{-1}\|_{\rm var}\lesssim \|f\|_{p}\left(\mE^{\bQ_x}\|\sI_{b}\|_{\mH^{q_\H}_t}^2\right)^{1/2}\\
&\lesssim t^{\beta_H}\|f\|_{p}\nor |b|^{q_\H}\nor_{\mL^{q_b/q_\H}_T\widetilde{\mL}^{p_b/q_\H}}=t^{\beta_H}\|f\|_{p}\nor b\nor^{q_\H}_{\mL^{q_b}_T\widetilde{\mL}^{p_b}}.
\end{align*}
The proof is complete.
\end{proof}
\br
In comparison with \cite[Theorem 2]{NO03}, we relax the condition $p_1,q_1\ge2$ in \cite[Theorem 2]{NO03} to $p_1,q_1\ge 1/(1-H)$ in Theorem \ref{Th25}.  
This allows us to treat the Biot-Savart kernel in subsequent discussions.
\er

It should be noted that recently, Butkovsky and Gallay \cite{BG23} showed the existence of the weak solution under the  weaker assumption for $b\in \mL^q_T\mL^p$ with
\begin{align*}
\tfrac{1-H}{q}+\tfrac{Hd}{p}<1-H\Leftrightarrow\tfrac{Hd}{p}<(1-H)(1-\tfrac1q),
\end{align*}
which coincides with the result in \cite{Kr20} for $H=1/2$. But the uniqueness in this case is still open. 
However, 
based on the entropy estimate in Lemma \ref{Le45}, we have the following partial result.

\bt\label{Th44}
Let $H\in(0,\frac12)$ and $p_1,q_1\in[\frac1{1-H},\infty)$ satisfy $\frac{Hd}{p_1}+\frac{1-H}{q_1}<(1-H)^2$. Assume $b\in \mL^{q_1}_T\mL^{p_1}$. Let $b_n$ be
a sequence of bounded smooth function converging to $b$ in $\mL^{q_1}_T\mL^{p_1}$ as $n\to\infty$. For $x\in\mR^d$, let $X^n$ be the unique strong solution of
$$
X^n_t=x+\int^t_0b_n(s, X^n_s)\dif s + W^H_t,\ \ t\in[0,T].
$$
Then the law $\bP_n$ of $X^n$ in $\mC_T$ weakly converges to a solution $\bP\in\cS(b,\delta_x)$.  We call such $\bP$ a regular solution of SDE$_b$. Moreover, for any two
$b_1, b_2\in \mL^{q_1}_T\mL^{p_1}$, letting $\bP_i$ be the unique regular solution of SDE$_{b_i}$ starting from $x$, where $i=1,2$, we have
\begin{align}\label{SK6}
\cH(\bP_1|\bP_2)\lesssim_C\mE^{\bP_1}\|\sI_{b_1-b_2}\|^2_{\mH^{q_{\H}}_T}.
\end{align}
\et
\begin{proof}
Let $\frac{Hd}p+\frac{1-H}q<1-H$ and $m\in\mN$.
By \cite[Lemma 3.11]{BG23}, there is a constant $C>0$ such that for all $n\in\mN$ and $f\in\mL^q_T\mL^p$,
\begin{align}\label{SK1}
\mE^{\bP_n}\left|\int^T_0f(s,w_s)\dif s\right|^m=\mE\left|\int^T_0f(s,X^n_s)\dif s\right|^m\leq C\|f\|^m_{\mL^q_T\mL^p}.
\end{align}
By Lemma \ref{Le45} and the above Krylov's estimate with $q=q_1/q_H$ and $p=p_1/q_H$, where $q_H=1/(1-H)$, we have
\begin{align*}
\cH(\bP_n|\bP_m)\lesssim\mE^{\bP_n}\|\sI_{b_n-b_m}\|^2_{\mH^{q_{\H}}_T}
\lesssim \||b^n-b^m|^{q_{\H}}\|^{2/q_{\H}}_{\mL^q_T\mL^p}=\|b_n-b_m\|^{2}_{\mL^{q_1}_T\mL^{p_1}},
\end{align*}
which implies by CKP's inequality \eqref{CKP},
$$
\lim_{n,m\to\infty}\|\bP_n-\bP_m\|^2_{\rm var}\leq 2\lim_{n,m\to\infty}\cH(\bP_n|\bP_m)\lesssim\lim_{n,m\to\infty}\|b_n-b_m\|^{2}_{\mL^{q_1}_T\mL^{p_1}}=0.
$$
Let $\bP\in\cP(\mC_T)$ be the limit point so that
\begin{align}\label{SK3}
\lim_{n\to\infty}\|\bP_n-\bP\|^2_{\rm var}=0.  
\end{align}
It is easy to see $\bP\in\cS(b,\delta_x)$ by taking limits. In fact, it suffices to show that for any $k\in\mN$, $t_1\leq t_2\leq \cdots\leq t_k$ and $f\in C^1_b(\mR^{kd})$,
\begin{align}\label{SK5}
\mE^\mP f(W^H_{t_1},\cdots,W^H_{t_k})=\mE^{\bP} f(\sW^{b}_{t_1},\cdots,\sW^{b}_{t_k}),
\end{align}
where $W^H$ is an {\rm f}Bm on some probability space $(\Omega,\sF,\mP)$ and $\sW^b_\cdot$ is defined by \eqref{SK2}.
Since $\sW^{b_n}_\cdot$ is an {\rm f}Bm w.r.t. $\bP_n$, we have
\begin{align}\label{SK4}
\mE^\mP f(W^H_{t_1},\cdots,W^H_{t_k})= \mE^{\bP_n} f(\sW^{b_n}_{t_1},\cdots,\sW^{b_n}_{t_k}).
\end{align}
By \eqref{SK1}, we have
\begin{align*}
&|\mE^{\bP_n} f(\sW^{b_n}_{t_1},\cdots,\sW^{b_n}_{t_k})-\mE^{\bP_n} f(\sW^{b}_{t_1},\cdots,\sW^{b}_{t_k})|\\
&\quad\leq\|\nabla f\|_\infty\sum_{j=1}^k\mE^{\bP_n}\left(\int^{t_k}_0|b_n-b|(s,w_s)\dif s\right)\\
&\quad\lesssim\|b_n-b\|_{\mL^{q_1}_T\mL^{p_1}}\to 0,\ \ n\to\infty.
\end{align*}
Moreover, by \eqref{SK3} we also have
$$
\lim_{n\to\infty}|\mE^{\bP_n} f(\sW^{b}_{t_1},\cdots,\sW^{b}_{t_k})-\mE^{\bP} f(\sW^{b}_{t_1},\cdots,\sW^{b}_{t_k})|=0.
$$
By taking limits for \eqref{SK4}, we obtain \eqref{SK5}.

For $i=1,2$, let $\bP_i$ be the unique regular solution of SDE$_{b_i}$ with the same starting point $x\in\mR^d$. 
Let $b^{(n)}_i$ be the smooth approximation sequence of $b_i$, and
$\bP^{(n)}_i$ the law of the solution of the associated approximation SDE. By Lemma \ref{Le45} we have
$$
\cH(\bP^{(n)}_1|\bP^{(n)}_2)\lesssim_C\mE^{\bP^{(n)}_1}\big\|\sI_{b^{(n)}_1-b^{(n)}_2}\big\|^2_{\mH^{q_{\H}}_T}.
$$
Since $\cH(\mu|\nu)$ is lower semi-continuous w.r.t. $\mu,\nu$, by taking limits and as above, we get \eqref{SK6}.
\end{proof}

\br
Assume $q_1=\infty$ and $H\in(1-\frac1{\sqrt2}, \frac12)$. The condition  $dH/p_1<(1-H)^2$ is worse than $dH/p_1<1/2$.
Thus in this case, the result in Theorem \ref{Th25} is better than Theorem \ref{Th44}.
\er

\subsection{DFSDEs driven by {\rm f}Bm}\label{42}
In this subsection we consider the following DFSDE driven by {\rm f}Bm:
\begin{align}\label{DFSDE0}
X^{x}_{s,t}=x+\int_s^t B(r,X^{x}_{s,r},\mu^{\centerdot}_{r,T},\mu^{\centerdot}_{s,r})\dif r+W^H_t-W^H_s,\ \ (s,t)\in\mD_T,
\end{align}
where $B:[0,T]\times\mR^d\times\cC\cP_0\times \cC\cP_0\to\mR^d$ satisfies the following assumption:
\begin{itemize}
\item[{\bf (H$^s_1$)}] Let $(p_1,q_1)\in[\frac1{1-H},\infty]^2$ satisfy $\tfrac{1}{q_1}+\tfrac{Hd}{p_1}<\tfrac12$. 
There is a  $\kappa_1>0$,
\begin{align}\label{A1}
\nor B(\cdot,\mu^{\centerdot},\nu^{\centerdot})\nor_{\mL^{q_1}_T\tilde\mL^{p_1}}\le\kappa_1,\ \ \forall \mu^{\centerdot},\nu^{\centerdot}\in\cC\cP_0,
\end{align}
and there is a function $\ell\in \mL^{q_1}_T$ such that for all $\mu^{\centerdot}_i,\nu^{\centerdot}_i\in\cC\cP_0$, $i=1,2$,
\begin{align}\label{A2}
\nor B(t,\cdot,\mu^{\centerdot}_1,\nu^{\centerdot}_1)-B(t,\cdot,\mu^{\centerdot}_2,\nu^{\centerdot}_2)\nor_{p_1} \le \ell(t)(\|\mu^{\centerdot}_1-\mu^{\centerdot}_2\|_{\cC\cP_0}+\|\nu^{\centerdot}_1-\nu^{\centerdot}_2\|_{\cC\cP_0}).
\end{align}
\end{itemize}

\bt\label{Th46}
Under {\bf (H$^s_1$)}, there is a unique weak solution to DFSDE \eqref{DFSDE0}. Moreover, for any $p>1$, there is a constant $C=C(p,d,\kappa_1)>0$ such that for all $(s,t)\in\mD_T$,
\begin{align}\label{SQ1}
\sup_{x\in\mR^d}\mE f(X^{x}_{s,t})\lesssim_C (t-s)^{-Hd/p}\nor f\nor_p,
\end{align}
and for any $p,q\in(1,\infty]$ with $\alpha:=1-(\frac1q+\frac{Hd}p)>0$, there is a   $C=C(T,H,d,p_1,q_1,\kappa_1,p,q)>0$ such that for all  $m\geq 1$ and $(s,t)\in\mD_T$,
\begin{align}\label{0123:51}
\sup_{x\in\mR^d}\left\|\int^t_sf(r,X^{x}_{s,r})\dif r\right\|_{L^m(\Omega)}\lesssim_Cm^{1-\alpha}\nor f\nor_{\mL^q_T\tilde{\mL}^p}.
\end{align}
\et
\begin{proof}
We use the method of freezing the distribution-flow as Theorem \ref{Thm21}. 
Let $\mu^{\centerdot,0}_{s,t}=\delta_x$. For $n\in\mN$, define the following approximation sequence:
\begin{align*}
X^{x,n+1}_{s,t}=x+\int_s^t B(r,X^{x,n+1}_{s,r},\mu^{\centerdot,n}_{r,T},\mu^{\centerdot,n}_{s,r})\dif r+W^H_t-W^H_s,
\end{align*}
where $\mu^{x,n}_{s,t}$ is the law of $X^{x,n}_{s,t}$.
By Theorem \ref{Th25}, there is a unique weak solution to the above approximation SDE, and by \eqref{0123:41}, for any $p,q\in[1,\infty]$ with $\alpha:=1-(\frac1q+\frac{Hd}p)>0$, there is a constant $C=C(T,H,d,p_1,q_1,\kappa_1,p,q)>0$ such that for all  $m\geq 1$ and $(s,t)\in\mD_T$,
\begin{align}\label{0123:61}
\sup_{n\in\mN}\sup_{x\in\mR^d}\left\|\int^t_sf(r,X^{x,n}_{s,r})\dif r\right\|_{L^m(\Omega)}\lesssim_C m^{1-\alpha}\nor f\nor_{\mL^q_T\tilde{\mL}^p}.
\end{align}
For simplicity of notations, for any $n,k\in\mN$, we write
$$
b_{n,k}(r,x):=B(r,x,\mu^{\centerdot,n}_{r,T},\mu^{\centerdot,n}_{s,r})-B(r,x,\mu^{\centerdot,k}_{r,T},\mu^{\centerdot,k}_{s,r}).
$$
Noting that by Lemma \ref{Le45}, 
\begin{align*}
&\cH(\mu^{x,n+1}_{s,t} |\mu^{x,k+1}_{s,t})\lesssim_C\mE\left(\int_s^t\big|b_{n,k}(r,X^{x,n}_{s,r})\big|^{q_{\H}}\dif r\right)^{2/q_{\H}},
\end{align*}
by CKP's inequality \eqref{CKP} and \eqref{0123:61} with $(p,q)=(\frac{p_1}{q_{\H}},\frac{q_1}{q_{\H}})\in[1,\infty]$ and $m=\frac2{q_{\H}}>1$, we have 
\begin{align*}
\|\mu^{\centerdot,n+1}_{s,t}-\mu^{\centerdot,k+1}_{s,t}\|_{\cC\cP_0}^2
&\leq 2\sup_{x\in\mR^d}\cH(\mu^{x,n+1}_{s,t} |\mu^{x,k+1}_{s,t})
\lesssim  \left\|\int_s^t\big|b_{n,k}(r,X^{x,n}_{s,r})\big|^{q_{\H}}\dif r\right\|_{L^{2/q_{\H}}(\Omega)}^{2/q_{\H}}\\
&\lesssim \nor\1_{[s,t]} |b_{n,k}|^{q_{\H}}\nor^{2/q_{\H}}_{\mL^{q_1/q_{\H}}\tilde\mL^{p_1/q_{\H}}_x}
=\nor\1_{[s,t]} b_{n,k} \nor^{2}_{\mL^{q_1}\tilde\mL^{p_1}_x}\\
&\!\!\!\stackrel{\eqref{A2}}{\le} \left(\int_s^t\ell(r)^{q_1}\left[\|\mu^{\centerdot,n}_{r,T}-\mu^{\centerdot,k}_{r,T}\|_{\cC\cP_0}+\|\mu^{\centerdot,n}_{s,r}-\mu^{\centerdot,k}_{s,r}\|_{\cC\cP_0}\right]^{q_1}\dif r\right)^{2/q_1}.
\end{align*}
By Gronwall's inequality in Lemma \ref{lemA03}, we derive that for each $0\leq s\leq t\leq T$,
\begin{align*}
\varlimsup_{n,k\to\infty}\|\mu^{\centerdot,n}_{s,t}-\mu^{\centerdot,k}_{s,t}\|_{\cC\cP_0}=0.
\end{align*}
Hence, there is a $\mu^{\centerdot}_{s,t}\in \cC\cP_0$ such that
\begin{align*}
\lim_{n\to\infty}\|\mu^{\centerdot,n}_{s,t}-\mu^{\centerdot}_{s,t}\|_{\cC\cP_0}=0.
\end{align*}
Thus, for each $x\in\mR^d$, by \eqref{A1} and Theorem \ref{Th25}, there is a unique weak solution $X^{x}_{s,t}$ to SDE
\begin{align*}             
X^{x}_{s,t}=x+\int_s^t B(r,X^{x}_{s,r},\mu^{\centerdot}_{r,T},\mu^{\centerdot}_{s,r})\dif r+W^H_t-W^H_s.
\end{align*}
 By the same argument as above, we have
 \begin{align*}
&\|\bP\circ(X^{x}_{s,t})^{-1}-\mu^{x,n}_{s,t}\|_{\rm var}^2\le 2\cH(\bP\circ(X^{x}_{s,t})^{-1}|\mu^{x,n}_{s,t} )\\
&\quad\le  \left(\int_s^t\ell(r)^{q_1}\left[\|\mu^{\centerdot,n}_{r,T}-\mu^{\centerdot}_{r,T}\|_{\cC\cP_0}+\|\mu^{\centerdot,n}_{s,r}-\mu^{\centerdot}_{s,r}\|_{\cC\cP_0}\right]^{q_1}\dif r\right)^{2/q_1}\to 0,
\end{align*}
 as $n\to\infty$, which implies that
 \begin{align*}
\bP\circ(X^{x}_{s,t})^{-1}=\mu^{x}_{s,t}.
\end{align*}
By \eqref{0123:01} and \eqref{0123:41}, we have \eqref{SQ1} and \eqref{0123:51}. The proof is complete.
 \end{proof}

\bx\label{ex:0119:00}
Let $b_1\in \mL^1$ 
and $b_2\in \mL^p$ with some $p>d$. Let $\varphi_1\in \mL^\infty$ and $\varphi_2\in \mL^1$.
For $\mu,\nu\in\cC\cP_0$, we introduce
\begin{align*}
B(t,x,\mu^\centerdot,\nu^\centerdot):=\int_{\mR^{d}}b_1(x-y)\mu^y(\varphi_1)\dif y+\int_{\mR^{d}}(b_2*\nu^z)(x)\varphi_2(z)\dif z.
\end{align*}
Then it is easy to see that  {\bf (H$^s_1$)} holds with $(q_1,p_1)=(\infty,p)$. Indeed, 
\begin{align*}
\nor B(t,\cdot,\mu^\centerdot,\nu^\centerdot)\nor_p&\le\left\|\int_{\mR^{d}}b_1(\cdot-y)\mu^y(\varphi_1)\>\dif y\right\|_\infty+\left\|\int_{\mR^{d}}(b_2*\nu^z)(\cdot)\varphi_2(z)\dif z\right\|_p\\
&\le  \| b_1*\mu^\centerdot(\varphi_1)\|_\infty+\|\varphi_2\|_1 \sup_z \|b_2*\nu^z\|_p\\
&\le\|b_1\|_1\|\varphi_1\|_\infty+\|\varphi_2\|_1\|b_2\|_p.
\end{align*}
Moreover, we also have
\begin{align*}
&\nor B(t,\cdot,\mu^{\centerdot,1},\nu^{\centerdot,1})-B(t,\cdot,\mu^{\centerdot,2},\nu^{\centerdot,2})\nor_p
\le \|b_1\|_1\|\varphi_1\|_\infty  \|\mu^{1}-\mu^{2}\|_{\cC\cP_0}+\|b_2\|_p\|\varphi_2\|_1 \|\nu^{1}-\nu^{2}\|_{\cC\cP_0}.
\end{align*}
In Section \ref{sec:05}, we shall use Theorem \ref{Th46} to study the  2D-Navier-Stokes equation with {\rm f}Bm.
\ex

\subsection{Backward DFSDEs driven by Brownian motion}\label{43}
In this section, we consider the following backward DFSDE driven by Brownian motion:
\begin{align}\label{0117:05}
X^x_{s,t}=x+\int_s^t B(r,X^x_{s,r},\mu^\centerdot_{r,T})\dif r+\sqrt{2}(W_t-W_s),
\end{align}
where for some $p_0\in(1,\infty)$,  
$$
B: [0,T]\times\mR^d\times\tilde{\cL}^{p_0}\cP_s \mbox{(or $\cL^{p_0}\cP_s$)}\to\mR^d
$$
is a measurable vector field.
We first consider the following classical SDE
\begin{align}\label{SDE12}
X_{s,t}^x=x+\int_s^t b(r,X_r^x)\dif r+\sqrt{2}(W_t-W_s),\ \ t\in[s,T].
\end{align}
The following results were partly obtained in \cite{XXZZ21}.
\bt\label{Pr51}
Let  $(p_1,q_1)\in(2,\infty)^2$ satisfy $\alpha:=1-(\frac2{q_1}+\frac{d}{p_1})>0$. Assume that
$$
\kappa_b:=\nor b\nor_{\mL^{q_1}_T\tilde\mL^{p_1}}<\infty.
$$
\begin{enumerate}[(i)]
\item For each $(s,x)\in[0,T)\times\mR^d$, there is a unique strong solution $X^x_{s,t}=:X^x_{s,t}(b)$ to  SDE \eqref{SDE12}. Moreover, $x\mapsto X^x_{s,t}(b)$ is weakly differentiable and for any $p\geq 1$,
\begin{align}\label{AZ82}
\sup_{(t,x)\in[s,T]\times\mR^d}\mE|\nabla X^x_{s,t}(b)|^p<\infty.
\end{align}
\item For any $p\in(1,p_1]$, there are constants $C_1, C_2>0$ only depending on $T,d,p_1,q_1,p$ 
such that for all $(s,t)\in\mD_T$,
\begin{align}\label{AZ8}
\nor \mP\circ(X^\cdot_{s,t}(b))^{-1}\nor_{p} \lesssim_{C_1}1+ (t-s)^{\alpha/2}\exp\left\{C_2\left(1+\kappa_b^{2/\alpha}\right)\right\}.
\end{align}
\item Let $p\in(1,p_1]$. For any $b_1,b_2\in \mL^{q_1}_T\tilde\mL^{p_1}$, there is a constant $C_3=C_3(T,d,p_1,q_1,p)>0$ such that for all $(s,t)\in\mD_T$,
\begin{align}\label{0228:02}
\nor \mP\circ (X^\cdot_{s,t}(b_1))^{-1}-\mP\circ (X^\cdot_{s,t}(b_2))^{-1}\nor_p\lesssim_{C_3}\int_s^t (t-r)^{-\frac{1+d/p_1}{2}}\nor b_1(r)-b_2(r)\nor_{p_1}\dif r.
\end{align}
\item If  $\div b=0$, then for any $f\in L^1$,
\begin{align}\label{AZ18}
\int_{\mR^d}\mE |f(X^x_{s,t})|\dif x=\|f\|_1.
\end{align}
\end{enumerate}
\et
\begin{proof}
(i) The existence and uniqueness of strong solutions and estimate \eqref{AZ82} follow by \cite[Theorem 1.1]{XXZZ21}.

(ii) For \eqref{AZ8}, we fix $t\in(0,T]$ and $\phi\in C^\infty_c(\mR^d)$. Let $u(s,x):=\mE \phi(\sqrt2W_{t-s}+x)$.  Then $u\in C^1([0,t];  C^2_b(\mR^d))$ and
$$
\p_s u+\Delta u=0,\ \ u(t)=\phi.
$$
By  It\^o's formula, we have
$$
\mE\phi(X^{x}_{s,t})=\mE u(t,X^{x}_{s,t})=u(s,x)+\mE\int^t_s(b\cdot\nabla u)(r, X^{x}_{s,r})\dif r.
$$
Let $p\in[1,p_1]$ and $p_2,q_2\in[1,\infty)$ be defined by
\begin{align}\label{FS1}
\tfrac1{q_2}+\tfrac1{q_1}=1,\ \ \tfrac1{p_2}+\tfrac1{p_1}=\tfrac1{p}.
\end{align}
By Lemma \ref{LeA1} with $j=0$, \eqref{0123:01} with $H=\frac12$ and H\"older's inequality, we have
\begin{align}
\nor\mE\phi(X^{\cdot}_{s,t})\nor_{p}
&\leq\nor u(s,\cdot)\nor_p+\int^t_s\nor\mE (b\cdot\nabla u)(r, X^{\cdot}_{s,r})\nor_p\dif r\no\\
&\lesssim_C\nor \phi\nor_p+\exp\left\{C\kappa_b^{2/\alpha}\right\}\int^t_s\nor (b\cdot\nabla u)(r)\nor_p\dif r\no\\
&\lesssim_C\nor \phi\nor_p+\exp\left\{C\kappa_b^{2/\alpha}\right\}\nor b\nor_{\mL^{q_1}_T\tilde\mL^{p_1}}
\nor\nabla u\nor_{\mL^{q_2}_{[s,t]}\tilde\mL^{p_2}}.\label{SD0}
\end{align}
Note that by \eqref{AZ09} with $q=p_2$ and $p=p_1$,
\begin{align*}
\nor \nabla u\nor_{\mL^{q_2}_{[s,t]}\tilde\mL^{p_2}}&=\left(\int_s^t \nor \mE\nabla\phi(\sqrt2W_{t-r}+\cdot)\nor_{p_2}^{q_2}\dif r\right)^{1/q_2}\\
&\lesssim
\nor \phi\nor_p\left(\int^t_s(t-r)^{q_2(d/p_2-d/p-1)/2}\dif r\right)^{1/q_2}\lesssim \nor \phi\nor_p (t-s)^{(1-2/q_1-d/p_1)/2}.
\end{align*}
Substituting this into \eqref{SD0} and by \eqref{AW02}, we derive the estimate \eqref{AZ8}.

(iii) For simplicity we set $X^{\cdot,i}_{s,t}:=X^{\cdot}_{s,t}(b_i)$, $i=1,2$. 
We fix $t\in[0,T]$ and consider the following backward PDE:
\begin{align*}
\p_s u+\Delta u+b_1\cdot\nabla u=0,\quad u(t)=\phi\in C^\infty_c(\mR^d).
\end{align*}
Since $b_1\in \mL^{q_1}_T\tilde\mL^{p_1}$, by Theorem \ref{ThA}, there is a unique solution $u$ to the above equation.
Then by the generalized It\^o's formula (see \cite{XXZZ21}), we have
\begin{align*}
\mE u(t,X^{x,2}_{s,t})&=u(s,x)+\mE\int^t_s(\p_s u+\Delta u+b_2\cdot\nabla u)(r, X^{x,2}_{s,r})\dif r\\
&=\mE\int^t_s((b_2-b_1)\cdot\nabla u)(r, X^{x,2}_{s,r})\dif r,
\end{align*}
and
$$
\mE u(t,X^{x,1}_{s,t})=u(s,x)+\mE\int^t_s(\p_s u+\Delta u+b_1\cdot\nabla u)(r, X^{x,1}_{s,r})\dif r=u(s,x).
$$
Hence,
\begin{align*}
\mE \phi(X^{x,2}_{s,t})-\mE \phi(X^{x,1}_{s,t})
=\mE\int^t_s((b_2-b_1)\cdot\nabla u)(r, X^{x,2}_{s,r})\dif r,
\end{align*}
and by  \eqref{0123:01},
\begin{align*}
\nor\mE \phi(X^{\cdot,2}_{s,t})-\mE \phi(X^{\cdot,1}_{s,t})\nor_{p}
&\leq \int^t_s\nor\mE((b_2-b_1)\cdot\nabla u)(r, X^{\cdot,2}_{s,r})\nor_{p}\dif r\\
&\lesssim \int^t_s\nor( (b_2-b_1)\cdot\nabla u)(r, \cdot)\nor_{p}\dif r.
\end{align*}
Let $p_2,q_2$ be defined by \eqref{FS1}.
By  H\"older's inequality and Theorem \ref{ThA}, we have
\begin{align*}
\nor\mE \phi(X^{\cdot,1}_{s,t})-\mE \phi(X^{\cdot,2}_{s,t})\nor_{p}
&\lesssim_C\int^t_s\nor b_1-b_2\nor_{p_1}\nor\nabla u(r, \cdot)\nor_{p_2}\dif r\\
&\lesssim_C\nor\phi\nor_p\int^t_s\nor b_1-b_2\nor_{p_1}(t-r)^{-\frac{1+d/p-d/p_2}{2}}\dif r,
\end{align*}
which gives \eqref{0228:02} by taking the supremum of $\phi\in \bC^\infty_c$.

(iv) Let $b_n(t,x):=b(t,\cdot)*\rho_n(x)$ be the mollifying approximation of $b$. For each $x\in\mR^d$, let $X^{x,n}_{s,t}$ be the unique solution of approximation SDE
$$
X_{s,t}^{x,n}=x+\int_s^t b_n(r,X_{s,r}^{x,n})\dif r+\sqrt{2}(W_t-W_s).
$$
It is well known that (see \cite{XXZZ21})
$$
\lim_{n\to\infty}\mE|X^{x,n}_{s,t}-X^x_{s,t}|=0. 
$$
Since $\div b_n=0$, we have
$$
\mE\int_{\mR^d}f(X_{s,t}^{x,n})\dif x=\int_{\mR^d}f(x)\dif x.
$$
By taking limits, we obtain that for any $0\leq f\in C^\infty_c(\mR^d)$,
$$
\mE\int_{\mR^d}f(X_{s,t}^{x})\dif x=\int_{\mR^d}f(x)\dif x.
$$
The proof is complete by a further approximation.
\end{proof}

\br\label{rmk:0228}
If we replace all the norms of $\nor\cdot\nor_p$ by $\|\cdot\|_p$, then the results in Theorem \ref{Pr51} still hold.
\er

Now we make the following assumption on $B$:
\begin{enumerate}[{\bf (H$^s_2$)}]

\item For some $(p_1,q_1)\in(2,\infty)$ with $\frac2{q_1}+\frac{d}{p_1}<1$ and $p_0\in(1,p_1]$,
there is a function $\ell\in \mL^{q_1}_T$ such that for some $\beta\geq 0$ and all $t\in[0,T]$ and $\mu^\centerdot\in\tilde{\cL}^{p_0}\cP_s$,
\begin{align}\label{ZQ21}
\nor B(t,\cdot,\mu^{\centerdot})\nor_{p_1}\le \ell(t)(1+\nor \mu^{\centerdot}\nor^\beta_{p_0}),
\end{align}
and for all $t\in[0,T]$ and $\mu^\centerdot, \nu^\centerdot\in\tilde{\cL}^{p_0}\cP_s$,
\begin{align}\label{ZQ22}
\nor B(t,\cdot,\mu^{\centerdot})-B(t,\cdot,\nu^{\centerdot})\nor_{p_1}\le \ell(t)\nor \mu^{\centerdot}-\nu^{\centerdot}\nor_{p_0}.
\end{align}
\end{enumerate}

Now we can prove the main result of this section.

\bt\label{Th410}
Under {\bf (H$^s_2$)}, there is a time $T\in(0,1)$ such that for each $x\in\mR^d$, 
there is a unique strong solution to DFSDE \eqref{0117:05} on $[0,T]$.  When $\beta=0$, the time can be taken arbitrarily large.
Moreover, if we replace all the norms in {\bf (H$^s_2$)}  by $\|\cdot\|_p$, then the conclusion still holds.
\et

\begin{proof}
Let $\mu^{\centerdot,1}_{r,T}=\delta_x$. For $n\in\mN$, 
we consider the following Picard iteration to DFSDE \eqref{0117:05}:
\begin{align*}
X^{x,n+1}_{s,t}=x+\int_s^t B(r, X^{x,n+1}_{s,r},\mu^{\centerdot,n}_{r,T})\dif r+\sqrt2(W_t-W_s),
\end{align*}
where $\mu^{x,n}_{r,T}$ is the law of $X^{x,n}_{r,T}$.
By \eqref{AZ8} with $p=p_0$, we have
$$
\kappa_{\mu^{n+1}_{\cdot,T}}=\sup_{s\in[0,T]}\nor\mu^{\centerdot,n+1}_{s,T}\nor_{p_0}\leq C_1\left( 1+T^{\alpha/2}\exp\left\{C_2\kappa_{\mu^n_{\cdot,T}}^{2\beta/\alpha}\right\}\right),
$$
where $\alpha:=1-(\frac2{q_1}+\frac{d}{p_1})>0$.
If $\beta=0$, then
$$
\sup_n\kappa_{\mu^{n}_{\cdot,T}}<\infty.
$$
 If $\beta>0$, then one can choose a time $T$ small enough so that 
 $$
 C_1 T^{\alpha/2}\exp\left\{C_2(2C_1)^{2\beta/\alpha}\right\}\leq 1\ \mbox{and } \ \sup_n\kappa_{\mu^{n}_{\cdot,T}}\leq 2C_1.
 $$
 Now by \eqref{0228:02} we have for any $p\in(1,p_1]$,
\begin{align}\label{ZQ2}
\nor\mu^{\centerdot,n+1}_{s,t}-\mu^{\centerdot,m+1}_{s,t}\nor^{q_1'}_{p}\lesssim_{C_3}
\int_s^t (t-r)^{-q_1'\frac{1+d/p_1}{2}}\nor\mu^{\centerdot,n}_{r,T}-\mu^{\centerdot,m}_{r,T}\nor_{p_0}^{q_1'}\dif r.
\end{align}
In particular, if we choose $p=p_0$ and set
$$
h(t):=\varlimsup_{n,m\to\infty}\nor\mu^{\centerdot,n}_{t,T}-\mu^{\centerdot,m}_{t,T}\nor_{p_0}^{q_1'},
$$
then by Fatou's lemma,
$$
h(s)\leq \int_s^T (T-r)^{-q_1'\frac{1+d/p_1}{2}}h(r)\dif r.
$$
Since $q_1'\frac{1+d/p_1}{2}<1$, by Gronwall's inequality of Voltera's type, we get
$$
h(s)=\varlimsup_{n,m\to\infty}\nor\mu^{\centerdot,n}_{s,T}-\mu^{\centerdot,m}_{s,T}\nor_{p_0}^{q_1'}=0.
$$
By Proposition \ref{Le21}, for each $s\in[0,T]$, there is a probability kernel $\mu^\centerdot_{s,T}\in\tilde{\cL}^{p_0}\cP$ so that
$$
\lim_{n\to\infty}\nor\mu^{\centerdot,n}_{s,T}-\mu^{\centerdot}_{s,T}\nor_{p_0}=0.
$$
Now for each $(s,x)\in[0,T]\times\mR^d$, let $X^x_{s,t}$ be the unique strong solution of the following SDE:
\begin{align}\label{0228:00}
X^{x}_{s,t}=x+\int_s^t B(r, X^{x}_{s,r},\mu^{\centerdot}_{r,T})\dif r+\sqrt2(W_t-W_s).
\end{align}
By using \eqref{0228:02} again, it is easy to see that for each $s\in[0,T]$ and Lebesgue almost all $x\in\mR^d$,
$$
\mP\circ(X^x_{s,T})^{-1}=\mu^x_{s,T},
$$
 which gives the existence, and the uniqueness is from the stability estimate \eqref{0228:02}.

Now, let us replace all the norm $\nor\cdot\nor_p$ by $\|\cdot\|_p$. Then By the same argument, we have
\begin{align*}
\varlimsup_{n,m\to\infty}\|\mu^{\centerdot,n}_{s,T}-\mu^{\centerdot,m}_{s,T}\|_{p_0}=0.
\end{align*}
 The only difference from the localized $L^p$ case is that by Proposition \ref{Le21}, we can only find a sub-probability kernel $\mu^{\centerdot}_{s,T}\in \cL^{p_0}\cP_s$ such that
$$
\lim_{n\to\infty}\|\mu^{\centerdot,n}_{s,T}-\mu^{\centerdot}_{s,T}\|_{p_0}=0.
$$ 
 But this don't prevent us from considering the SDE \eqref{0228:00}. By using \eqref{0228:02} again, for each $s\in[0,T]$ and Lebesgue almost all $x\in\mR^d$,
$$
\mP\circ(X^x_{s,T})^{-1}=\mu^x_{s,T},
$$
which implies that $\mu^x_{s,T}$ is in fact a probability kernel. 
The proof is complete.
\end{proof}
\bx
Let $p_0\in(1,\infty)$ and $\phi\in\tilde\mL^{p_0}$. Let $p\in(1,\infty]$ satisfy
\begin{align}\label{GS1}
1+(\tfrac 1d\wedge\tfrac1{p_0})-\tfrac1{p_0}>\tfrac1p.
\end{align}
Suppose that $K\in(\bar\mL^p)^d$ and $g:\mR\to\mR$ is a Lipschitz function. Consider the following example:
\begin{align*}
B(x,\mu^\centerdot):=\int_{\mR^d}K(x-y)g(\mu^y(\phi))\dif y.
\end{align*}
By \eqref{GS1}, one can choose $p_1>d\vee p_0$ so that $1+\frac1{p_1}=\frac1p+\frac1{p_0}$. Thus by \eqref{AW2}, we have
$$
\nor B(\cdot,\mu^\centerdot)\nor_{p_1}\leq \nor K\nor^*_{p}\nor g(\mu^\centerdot(\phi))\nor_{p_0}
\leq \nor K\nor^*_{p}(|g(0)|+\|g\|_{\rm Lip}\nor \mu^\centerdot\nor_{p_0}\nor\phi\nor_{p_0}),
$$
and
\begin{align*}
\nor B(\cdot,\mu^\centerdot)-B(\cdot,\nu^\centerdot)\nor_{p_1}
&\leq \nor K\nor^*_{p}\nor g(\mu^\centerdot(\phi))-g(\nu^\centerdot(\phi))\nor_{p_0}
\leq \nor K\nor^*_{p}\|g\|_{\rm Lip}\nor\phi\nor_{p_0}\nor \mu^\centerdot-\nu^\centerdot\nor_{p_0}.
\end{align*}
Hence, Theorem \ref{Th410} can be applied to this case. In particular, if $g$ is bounded Lipschitz, then we have a global solution.
\ex

\section{{DFSDEs driven by {\rm f}Bm related to }the 2D-Navier-Stokes equations }\label{sec:05}

In this section we apply the previous results to prove Theorems \ref{thm1} and \ref{thm2}.
Let $\nu_0$ be a finite signed measure.
Consider the following {DFSDE driven by} fractional equation {related to the 2D-Navier-Stokes}:
\begin{align}\label{NSSDE2}
X_t^x=x+\int_0^t B_{\nu_0}(X^x_s,\mu^\centerdot_s)\dif s+W^H_t,
\end{align}
where
\begin{align*}
B_{\nu_0}(x,\mu^\centerdot)=\int_{\mR^2}(K_2*\mu^y)(x)\nu_0(\dif y),\ \ K_2(x)=\frac{(x_2,-x_1)}{2\pi|x|^2}.
\end{align*}
Theorems \ref{thm1} is an immediate consequence of the following result.
\bt\label{thm:51}
Let $H\in(0,\frac12)$.
For any vorticity $\nu_0$ being a finite singed measure, there is a unique strong solution $X^\cdot_t$ to DFSDE \eqref{NSSDE2}.
Moreover, if we let
\begin{align*}
u(t,x):=\int_{\mR^2} \mE K_2(x-X_{t}^y)\nu_0(\dif y)=B_{\nu_0}(x,\mu^\centerdot_t),
\end{align*}
then for any $p>1$ and $j\in\mN$, there is a constant $C>0$ such that for all $t\in(0,T]$,
\begin{align}\label{SU5}
\|\nabla^j u(t)\|_p\lesssim_Ct^{-2H(p-1)/p-(j-1)H}\|\nu_0\|^j_{\rm var}.
\end{align}
Moreover, for any $p\in(1,2)$ and $\eps>0$, there is a constant $C>0$ such that {for all $0<t\le T$,}
\begin{align*}
\nor u(t)-u(0)\nor_{p}\lesssim_C t^{[H(\frac{2}{p}-1)]\wedge[\frac{1-2H}{1-H}]-\eps},
\end{align*}
and {for all $0<s<t\le T$,}
\begin{align*}
\|u(t)-u(s)\|_\infty\lesssim_Cs^{-2H}|t-s|^{\frac{H}{3}-\eps}.
\end{align*}
\et
\begin{proof}
Note that $K_2=K_2\1_{D_0}+K_2\1_{D_0^c}$, where $D_0$ is the unit cube with center $0$. 
For any $p\in(1,2)$, by Minkowskii's inequality and $\mL^p+\mL^\infty\subset\tilde\mL^p$, we clearly have
\begin{align}
\nor B_{\nu_0}(\cdot,\mu^\centerdot)\nor_p
&\leq\left\|\int_{\mR^2}((K_2\1_{D_0})*\mu^y)(\cdot)\nu_0(\dif y)\right\|_p+\left\|\int_{\mR^2}((K_2\1_{D^c_0})*\mu^y)(\cdot)\nu_0(\dif y)\right\|_\infty\no\\
&\leq\int_{\mR^2}\|(K_2\1_{D_0})*\mu^y\|_p|\nu_0|(\dif y)+\int_{\mR^2}\|(K_2\1_{D^c_0})*\mu^y\|_\infty|\nu_0|(\dif y)\no\\
&\leq\left(\|K_2\1_{D_0}\|_p+\|K_2\1_{D^c_0}\|_\infty\right)\|\nu_0\|_{\rm var},\label{Con2}
\end{align}
and
\begin{align*}
\nor B_{\nu_0}(\cdot,\mu_1^\centerdot)-B_{\nu_0}(\cdot,\mu_2^\centerdot)\nor_p
\leq\left(\|K_2\1_{D_0}\|_p+\|K_2\1_{D^c_0}\|_\infty\right)\|\nu_0\|_{\rm var}\|\mu^{\centerdot}_1-\mu^{\centerdot}_2\|_{\cC\cP_0}.
\end{align*}
Let $H\in(0,\frac12)$. One can choose $p_1<(\frac1{1-H},2)$ and $q_1=\infty$ so that $\frac{Hd}{p_1}+\frac{1}{q_1}<\frac12$.
Thus one sees that {\bf (H$^s_1$)} holds for the above $B_{\nu_0}$.
By Theorem \ref{Th46},  there is a unique weak solution $X^x_{t}$ to DFSDE \eqref{NSSDE2}.
Moreover, for any $p>1$, by \eqref{SQ1}, there is a constant $C=C(T,p,H)>0$ such that for all $f\in\tilde\mL^p$ and $t\in(0,T]$,
$$
\sup_{x}|\mE f(X_{t}^x)|\lesssim_C \nor f\nor_p t^{-2H/p},
$$
which in turn implies that $X^x_t$ admits a density $\rho^x_t(\cdot)\in \bar\mL^{p/(p-1)}$ with
$$
\sup_x\nor\rho^x_t\nor^*_{p/(p-1)}\lesssim_C t^{-2H/p}.
$$
Since $\nabla K_2$ is a Calder\'on-Zygmund kernel, for any $p\in(1,\infty)$, by the $L^p$-boundedness of singular integral operators, we have {(see Remark \ref{rmk:12})}
\begin{align}\label{SU1.0}
\|\nabla u(t)\|_p\leq\int_{\mR^d}\|\nabla (K_2*\rho^y_t)\|_p|\nu_0|(\dif y)\lesssim \sup_y\|\rho^y_t\|_p\|\nu_0\|_{\rm var}\lesssim t^{-2H(p-1)/p}\|\nu_0\|_{\rm var}.
\end{align}
Thus we get \eqref{SU5} for $j=1$.
%In particular, $X^x_t$ uniquely solves the following SDE
%$$X_t^x=x+\int_0^t u(s,X^x_s)\dif s+W^H_t.$$

For higher order derivative estimates, we use the Malliavin calculus. 
We first recall the main ingredients in the Malliavin calculus.
Let $\mu$ be the classical Wiener measure on $\mC_T$ so that the coordinate process $w$ is a $d$-dimensional standard Brownian motion. 
For an absolutely continuous function $h$ with $h(0)=0$ and $\int^T_0|\dot h(s)|^2\dif s<\infty$, the Malliavin derivative of a functional $F:\mC_T\to\mR$ is defined by
$$
D_h F(\omega):=\lim_{\eps\to 0}\frac{F(\omega+\eps h)-F(\omega)}{\eps}\ \ \mbox{in $L^2(\mC_T;\mu)$.}
$$
The following integration by parts formula holds:
$$
\mE^\mu(D_h F)=\mE^\mu\left(F\int^T_0\dot h(s)\dif w_s\right).
$$
Recall that on the classical Wiener space $(\mC_T,\mu)$, the {\rm f}Bm $W^H_t=\int^t_0K_H(t,s)\dif w_s$ can be considered as a Wiener functional.
Below we fix $t\in(0,T]$.
By  Girsanov's construction of weak solutions in Theorem \ref{Th25} we have
\begin{align}\label{SU3}
\mE^{\bP_x} \nabla^j f(w_t)=\mE^\mu (Z^x_t \nabla^j f(W^H_t+x)),
\end{align}
where
$$
Z^x_t:=\exp\left(-\int_0^t(\widetilde{\bf K}_H \sI^x_u)(s)\dif w_s-\frac12 \|\widetilde{\bf K}_H\sI^x_u\|^2_{\mL^2_t}\right),
$$
and
$$
\sI^x_u(t):=\int^t_0 u(s, W^H_s+x)\dif s.
$$
Since the initial point $x$ does not play any role in the following calculations, without loss of generality, we may assume $x=0$ and drop the superscript.
Note that for any $\gamma\geq 1$,
\begin{align}\label{SU1}
\sup_{t\in[0,T]}\mE^\mu|Z_t|^\gamma\leq C(\gamma,T,d,H,p,\|u\|_{\mL^\infty_T\mL^p}).
\end{align}
Fix $v\in\mR^2$. Let $h(s):=vs$ and $\nabla_v f:=\<\nabla f,v\>_{\mR^2}$. By simple calculations, we have for some constant $C_H>0$,
\begin{align}\label{SU8}
D_{h} W^H_t=\int^t_0K_H(t,s)v\dif s=C_Ht^{H+\frac12} v\Rightarrow C_H t^{H+\frac12}\nabla_v f(W^H_t)=D_{h} f(W^H_t),
\end{align}
and thus, by the integration by parts,
\begin{align}\label{SU4}
\begin{split}
&C_Ht^{H+\frac12}\mE^\mu (Z_t \nabla_v f(W^H_t))=\mE^\mu (Z_t D_{h} f(W^H_t))\\
&\quad=\mE^\mu \left(Z_t f(W^H_t)\<w_t,v\>_{\mR^2}\right)-\mE^\mu (D_hZ_t f(W^H_t)).
\end{split}
\end{align}
({\it Claim}:) For any $\gamma>1$, there is a constant $C=C(T,d,H,\gamma)>0$ such that for all $t\in(0,T]$,
\begin{align}\label{SU2}
\mE^\mu|D_hZ_t|^\gamma\leq C|v|^\gamma\|\nu_0\|_{\rm var}^\gamma t^{\gamma}.
\end{align}
By the chain rule, we have
$$
D_hZ_t=-Z_t\left(\int_0^t\<\widetilde{\bf K}_H \sI_u(s) ,v\>_{\mR^2}\dif s+\int_0^t(\widetilde{\bf K}_H D_h\sI_u)(s)\dif w_s+\<\widetilde{\bf K}_HD_h\sI_b,\widetilde{\bf K}_H\sI_u\>_{\mL^2_t}\right)
$$
and
$$
D_h\sI_u(s)=\int^s_0\<\nabla u(r, W^H_r), D_hW^H_r\>_{\mR^2}\dif r=C_H\int^s_0r^{\frac12+H}\nabla_v u(r, W^H_r)\dif r.
$$
By BDG's inequality and Minkowskii's inequality, for any $\gamma\geq 2$, we have
\begin{align*}
\mE^\mu\left|\int_0^t(\widetilde{\bf K}_H D_h\sI_u)(s)\dif w_s\right|^\gamma
&\lesssim
\mE^\mu\left(\|\widetilde{\bf K}_H D_h\sI_u\|_{L^2(0,t)}^{\gamma}\right)\stackrel{\eqref{AS2}}{\lesssim} \mE^\mu\left(\|D_h\sI_u\|_{\mH^{q_\H}_t}^{\gamma}\right)\\
&\lesssim\left(\int_0^t r^{(\frac12+H)q_{\H}}\|\nabla_v u(r, W^H_r)\|^{q_{\H}}_{L^\gamma(\mu)}\dif r\right)^{\gamma/q_{\H}}\\
&\lesssim |v|^\gamma\left(\int_0^t r^{(\frac12+H)q_{\H}{%
-\frac{2Hq_{\H}}{\gamma}}}\|\nabla u(r)\|^{q_{\H}}_{\gamma}\dif r\right)^{\gamma/q_{\H}}\\
&\overset{\eqref{SU1.0}}{\lesssim} |v|^\gamma\|\nu_0\|_{\rm var}^\gamma\left(\int_0^t r^{(\frac12+H)q_{\H}-2H q_{\H}/\gamma'{%\red
-\frac{2Hq_{\H}}{\gamma}}}\dif r\right)^{\gamma/q_{\H}}\\
&\lesssim |v|^\gamma\|\nu_0\|_{\rm var}^\gamma t^{\gamma(q_{\H}+\frac12-H)}%^{\frac\gamma{q_{\H}}+\gamma(\frac12+H)-2H(\gamma-1) }
\lesssim |v|^\gamma\|\nu_0\|_{\rm var}^\gamma t^{\gamma},
\end{align*}
{%\red 
where we used the following observation in the fourth inequality:
\begin{align*}
\|\nabla_v u(r, W^H_r)\|^{\gamma}_{L^\gamma(\mu)}=\int_{\mR^2}|\nabla_v u(r,x)|g^H_r(x)\dif x\le \|\nabla_v u(r)\|_\gamma^\gamma\|g^h_r\|_\infty\lesssim \|\nabla_v u(r)\|_\gamma^\gamma r^{-2H}, 
\end{align*}
with $g^h_r$ is the distributional density of the {\rm f}Bm $W^H_r$.
}
Similarly, we can show 
$$
\mE^\mu\left(\|\<\widetilde{\bf K}_H \sI_b ,v\>_{\mR^2}\|_{\mL^1_t}+|\<\widetilde{\bf K}_HD_h\sI_b,\widetilde{\bf K}_H\sI_b\>_{\mL^2_t}|\right)^\gamma\lesssim |v|^\gamma\|\nu_0\|_{\rm var}^\gamma t^{\gamma}.
$$
Combining the above calculations, by H\"older's inequality and \eqref{SU1}, we obtain \eqref{SU2}. Now by \eqref{SU3}, \eqref{SU4}, \eqref{SU2} and Lemma \ref{LeA1}, there is a constant $C>0$
such that for all $x\in\mR^d$ and $t\in(0,T]$,
$$
|\mE^{\bP_x} \nabla f(w_t)|\lesssim  \nor f\nor_p t^{-2H/p-H}\|\nu_0\|_{\rm var},
$$
which in turn implies that
$$
\sup_x\nor\nabla \rho^x_t\nor^*_{p/(p-1)}\lesssim_C t^{-2H/p-H}\|\nu_0\|_{\rm var}.
$$
By the same argument as in \eqref{SU1.0}, we obtain
$$
\|\nabla^2 u(t)\|_p\leq\int_{\mR^d}\|\nabla K_2*\nabla\rho^y_t\|_p|\nu_0|(\dif y)\lesssim \sup_y\|\nabla\rho^y_t\|_p\|\nu_0\|_{\rm var}\lesssim t^{-2Hp'-H}\|\nu_0\|^2_{\rm var}.
$$
This gives \eqref{SU5} for $j=2$. By induction, one can show \eqref{SU5} for $j=3,4,\cdots$.

Next we show the time regularity of $u$. Let $\chi\in\bC^\infty_0(\mR^d)$ with $\chi(0)=1$ and $p\in(1,2)$. 
By \eqref{0318:00}, one sees that for any $p_1\ge 1/(1-H)$,
\begin{align}
\nor u(t)-u(0)\nor_{p}&\le\left\|\int_{\mR^d}(K_2\chi)*(\mu^y_t-\delta_y)(\cdot)\nu_0(\dif y) \right\|_{p}+\left\|\int_{\mR^d}(K_2(1-\chi))*(\mu^y_t-\delta_y)(\cdot)\nu_0(\dif y) \right\|_{\infty}\no\\
&\le \sup_y\|(K_2\chi)*(\mu^y_t-\delta_y)\|_{p}+\sup_y\|(K_2(1-\chi))*(\mu^y_t-\delta_y)\|_{\infty}\no\\
&\le t^{\gamma H}\left(\|K_2\chi\|_{\bB^{\gamma}_p}+\|K_2(1-\chi)\|_{\bB^{\gamma}_{\infty}}\right)+t^{\beta_H}\|u\|_{\mL^\infty_T\widetilde{\mL}^{p_1}},\label{0318:02}
\end{align}
where $\beta_H:=1-2H/(p_1(1-H)).$ For any $p_1\in(1,2)$, 
by \eqref{Con2}, it is easy to see that 
$$
\sup_{t\in[0,T]}\nor u(t)\nor_{p_1}<\infty.
$$
and by \cite[Lemma A.3-(iv) and Proposition 2.5]{HRZ23},   for any $\gamma\in(0,2/p-1)$, 
$$
K_2\chi\in \bB^{\gamma}_p,\ \ K_2(1-\chi)\in \bC^\infty_b.
$$
Thus, for any $\eps>0$, one can choose $p_1$ close to $2$ so that
\begin{align*}
\nor u(t)-u(0)\nor_{p}\lesssim t^{[H(\frac{2}{p}-1)]\wedge[\frac{1-2H}{1-H}]-\eps}.
\end{align*}
On the other hand, for any $0<s<t\le T$,
\begin{align*}
|u(t,x)-u(s,x)|&\le \left|\int_{\mR^2} \mE (K_2(1-\chi))(x-X_{t}^y)-\mE (K_2(1-\chi))(x-X_{s}^y)\nu_0(\dif y)\right|\\
&+\left|\int_{\mR^2} \mE (K_2\chi)(x-X_{t}^y)-\mE (K_2\chi)(x-X_{s}^y)\nu_0(\dif y)\right|=:\sJ_1+\sJ_2.
\end{align*}
Since $K_2(1-\chi)\in \bC^\infty_b$, we have
\begin{align*}
\sJ_1&\le \|\nabla(K_2(1-\chi))\|_\infty \sup_{y} \mE|X_{t}^y-X_{s}^y|.
\end{align*}
In view of \eqref{0123:41}, we have for any $p\in(1,2)$
\begin{align}
\mE|X_{t}^y-X_{s}^y|&\le\mE\left|\int_s^t u(r,X_r^y)\dif r\right|+\mE|W_t^H-W_s^H|\le (t-s)^{1-\frac{2H}{p}}\nor u\nor_{\mL^\infty_T\widetilde{\mL}^p}+(t-s)^{H},\label{0318:06}
\end{align}
which by taking $p$  close  to $2$ implies that 
\begin{align*}
\sJ_1\lesssim (t-s)^{H},\quad \text{since $H<1/2$.}
\end{align*}
For $\sJ_2$,  let $\rho_\eps$ be the usual mollifiers. Note that
\begin{align*}
\sJ_2\le& \left|\int_{\mR^2} \mE ((K_2\chi)*\rho_\eps)(x-X_{t}^y)-\mE ((K_2\chi)*\rho_\eps)(x-X_{s}^y)\nu_0(\dif y)\right|\\
&+\left|\int_{\mR^2} \mE ((K_2\chi)*\rho_\eps-K_2\chi)(x-X_{t}^y)\nu_0(\dif y)\right|\\
&+\left|\int_{\mR^2} \mE ((K_2\chi)*\rho_\eps-K_2\chi)(x-X_{s}^y)\nu_0(\dif y)\right|.
\end{align*}
By \eqref{0318:06} and \eqref{0123:01}, we have for any $p_1,p_2\in(1,2)$ and $\gamma_2<2/p_2-1$,
\begin{align*}
\sJ_1&\le \|\nabla (K_2\chi)*\rho_\eps\|_\infty\sup_y\mE|X_t^y-X_s^y|+\sum_{r=s,t}\sup_y\mE |(K_2\chi)*\rho_\eps-K_2\chi|(x-X_{r}^y)\\
&\lesssim \eps^{-(1+\frac{2}{p_1})}\|K_2\chi\|_{p_1}|t-s|^{H}+s^{-\frac{2H}{p_2}}\|(K_2\chi)*\rho_\eps-K_2\chi\|_{p_2}\\
&\lesssim \eps^{-(1+\frac{2}{p_1})}|t-s|^{H}+s^{-\frac{2H}{p_2}}\eps^{\gamma_2}\|K_2\chi\|_{\bB^{\gamma_2}_{p_2}}.
\end{align*}
Now for any $\delta>0$, one chooses $p_1$ close to $2$ and $p_2$ close to $1$, and  $\eps=|t-s|^{\frac{H}{3}}$,  
\begin{align*}
\sJ_1\lesssim s^{-2H}|t-s|^{\frac{H}{3}-\delta}.
\end{align*}
This completes the proof.
\end{proof}
\br
We would like to mention the following open questions: 
\begin{enumerate}[$\bullet$]
\item 
 Can we show the limit of $H\to1/2$ and the regularity of $u$ in $t$?
\item  When $H=1/2$, it is well known that $\lim_{t\to\infty}\|u(t)\|_\infty=0$ when $\nu_0$ is a finite measure (see \cite{GG05}). Can we show the same assertion for $H<1/2$?
\item 
In \cite{Ha05}, the ergodicity was obtained for the solution to SDE driven by {\rm f}Bm. Is it possible to  estabilish the ergodicity of \eqref{NSSDE2}?
\end{enumerate}
\er
The above result does not work for $H=\frac12$ since \eqref{Con2} is no longer true for $p>2$. In what follows,
we consider the backward version of {DFSDE related to} Navier-Stokes equation:
\begin{align}\label{NSSDE3}
X_{s,t}^x=x+\int_s^t B_g(X^x_{s,r},\mu^\centerdot_{r,T})\dif r+\sqrt2(W_t-W_s),
\end{align}
where
\begin{align*}
B_g(x,\mu^\centerdot)=(K_2*\mu^\centerdot(g))(x).
\end{align*}
{The statement of the following theorem is already presented as Theorem \ref{thm2} in the introduction. Here, we provide its proof.}
\bt\label{thm:54}
Let $g\in \mL^{1+}=\cup_{p>1}\mL^p$.
For each $s\in[0,T]$ and $x\in\mR^2$, there is a unique strong solution $X^x_{s,t}$ to DFSDE \eqref{NSSDE3}.% {\blue such that $\sup_{0\le s<r\le T}\|B_g(X^x_{s,r},\mu^\centerdot_{r,T})\|_{p}<\infty$ with some $p>2$.}
Moreover,  $u(s,x):=B_g(x,\mu^\centerdot_{s,T})\in C([0,T); C^\infty_b(\mR^2))$ solves the following backward Navier-Stokes equation:
\begin{align}\label{NSE9}
\p_su+\Delta u+u\cdot \nabla u+\nabla p=0,\ \ u(T)=K_2*g .
\end{align}
\et
\begin{proof}
Let $p_0\in(1,2)$ and $g\in\mL^{p_0}$. 
We divide the proof into three steps.
In step 1, we  check the assumption {\bf (H$^s_2$)} with $\beta=0$ for the norm $\|\cdot\|_p$ and show the well-posedness of DFSDE \eqref{NSSDE3}. 
In step 2, we show the stability of the solution with respect to the initial value.
%use approximation to show that for any $k\in\mN_0$, there is a constant $C>0$ such that for all $t\in[0,T)$,
%\begin{align}\label{0301:07}
%\|\nabla^{k+1}u(t)\|_{p_0}\lesssim_C (T-t)^{-\frac{k}{2}}(1+\|g\|_{p_0}^2).
%\end{align}
In step 3, we show that $u$ is smooth and solves the 2D Navier-Stokes equation.

 {\bf (Step 1).} 
Let $p_1\in(2,\infty)$ satisfy $\frac1{p_1}+\frac12=\frac1{p_0}$.
For any $\mu^\centerdot, \nu^\centerdot\in\cL^{p_0}\cP_s$,
by Hard-Littlewood's inequality (see \cite[Theorem 1.7]{BCD}), we have
$$
\|{B}_g(\cdot,\mu^\centerdot)\|_{p_1}\lesssim 
\|\mu^\centerdot(g)\|_{p_0}\leq  \|\mu^\centerdot \|_{p_0}\|g\|_{p_0}
$$
and
$$
\|{B}_g(\cdot,\mu^\centerdot)-{B}_g(\cdot,\nu^\centerdot)\|_{p_1}\lesssim 
\|\mu^\centerdot(g)-\nu^\centerdot(g)\|_{p_0}\leq  \|\mu^\centerdot -\nu^\centerdot \|_{p_0}\|g\|_{p_0}.
$$
Since $\|{B}_g(\cdot,\mu^\centerdot)\|_{p_1}$ is not bounded in $\|\mu^\centerdot \|_{p_0}$, we need to truncate it. Define
\begin{align*}
\tilde{B}_g(x,\mu^{\centerdot}):=B_g(x,\mu^{\centerdot})\1_{\|\mu^{\centerdot}\|_{p_0}\leq1}+\frac{1}{\|\mu^{\centerdot}\|_{p_0}}B_g(x,\mu^{\centerdot})\1_{\|\mu^{\centerdot}\|_{p_0}>1}.
\end{align*}
Then it is easy to see that
\begin{align}\label{0229:01}
\|\tilde{B}_g(\cdot,\mu^\centerdot)\|_{p_1}\lesssim  \|g\|_{p_0},
\end{align}
and by the fact $\cL^{p_0}\cP_s\subset \cL(\mL^{p_0},\mL^{p_0})$ and Lemma \ref{lem:0301},
$$
\|\tilde{B}_g(\cdot,\mu^\centerdot)-\tilde{B}_g(\cdot,\nu^\centerdot)\|_{p_1}\lesssim 
\|\mu^\centerdot(g)-\nu^\centerdot(g)\|_{p_0}\leq  \|\mu^\centerdot -\nu^\centerdot \|_{p_0}\|g\|_{p_0}.
$$
Thus, by Theorem \ref{Th410},  for any $T>0$ there is a unique strong solution to
\begin{align*}
X_{s,t}^x=x+\int_s^t \tilde{B}_g(X^x_{s,r},\mu^\centerdot_{r,T})\dif r+\sqrt2(W_t-W_s).
\end{align*}
Noting that $\div \tilde{B}_g(\cdot,\mu^\centerdot)=0$,  by \eqref{0229:01} and \eqref{AZ18}, we in fact have
\begin{align*}
\tilde{B}_g(x,\mu^\centerdot_{r,T})=B_g(x,\mu^\centerdot_{r,T}).
\end{align*}
Then the strong well-posedness holds for DFSDE \eqref{NSSDE3}.

{\bf (Step 2).} 
Let $g_n\in C^\infty_c(\mR^2)$ be the smooth approximation of $g$ with $\|g_n-g\|_{\mL^{p_0}}\to 0$ as $n\to\infty$. 
Let $X^{x,n}_{s,t}$ be the unique solution to the following DFSDE
\begin{align*}
X_{s,t}^{x,n}=x+\int_s^t B_{g_n}(X^{x,n}_{s,r},\mu^{\centerdot,n}_{r,T})\dif r+\sqrt2(W_t-W_s),
\end{align*}
where $\mu^{x,n}_{s,t}$ is the law of $X^{x,n}_{s,t}$.
Let
$$
u_n(s,x):=B_{g_n}(x,\mu^{\centerdot,n}_{s,T})=\int_{\mR^2}K_2(x-y)\mu^{y,n}_{s,T}(g)\dif y.
$$
By \eqref{0228:02} and Remark \ref{rmk:0228}, one sees that
\begin{align*}
\|\mu^{\centerdot}_{s,t}-\mu^{\centerdot,n}_{s,t}\|_{p_0}&\lesssim \int_s^t (t-r)^{-\frac{1+d/p_1}{2}}\|B_g(\cdot,\mu^\centerdot_{r,T})-B_{g_n}(\cdot,\mu^{\centerdot,n}_{r,T})\|_{p_1}\dif r\\
&\lesssim \int_s^t (t-r)^{-\frac{1+d/p_1}{2}}\left(\|\mu^\centerdot_{r,T} -\mu^{\centerdot,n}_{r,T}\|_{p_0}\|g_n\|_{p_0}+ \|g_n-g\|_{p_0}\right)\dif r\\
&\lesssim \int_s^t (t-r)^{-\frac{1+d/p_1}{2}}\left(\|\mu^\centerdot_{r,T} -\mu^{\centerdot,n}_{r,T}\|_{p_0}+\|g_n-g\|_{p_0}\right)\dif r.
\end{align*}
Taking $t=T$, by Gronwall's inequality, we get
\begin{align*}
\sup_{s\in[0,T]}\|\mu^{\centerdot}_{s,T}-\mu^{\centerdot,n}_{s,T}\|_{p_0}\lesssim  \|g_n-g\|_{p_0}.
\end{align*}
Therefore, 
\begin{align}
\|u_n-u\|_{\mL^\infty_T\mL^{p_1}}&=\sup_{s\in[0,T]}\|B_{g_n}(\cdot,\mu^{\centerdot,n}_{s,T})-B_{g}(\cdot,\mu^{\centerdot}_{s,T})\|_{p_1}\no\\
&\lesssim \sup_{s\in[0,T]}\|\mu^\centerdot_{s,T} -\mu^{\centerdot,n}_{s,T}\|_{p_0}\|g_n\|_{p_0}+ \|g_n-g\|_{p_0}\to 0\quad \text{as $n\to\infty$}.\label{0229:03}
\end{align}

{\bf (Step 3).}  Since $g_n\in C^\infty_c(\mR^2)$ is smooth, it is well-known that $u_n\in C([0,T); C^\infty_b(\mR^2))$ solves
$$
\p_su_n+\Delta u_n+u_n\cdot \nabla u_n+\nabla p_n=0,\ \ u_n(T)=K_2*g_n,
$$
and for any $T_0<T$,
$$
\sup_{s\in[0,T_0]}\sup_n\|\nabla^k u_n(s)\|_\infty<\infty,
$$
which together with \eqref{0229:03} implis that $u\in C([0,T); C^\infty_b(\mR^2))$ and solves \eqref{NSE9}.
\end{proof}

\begin{appendix}
\renewcommand{\thetable}{A\arabic{table}}
\numberwithin{equation}{section}

\section{Proofs of Propositions \ref{Le20} and \ref{Le21}}

\begin{proof}[Proof of Proposition \ref{Le20}]

Equivalences \eqref{AW1} are proven in \cite{RZ21}. Let us prove \eqref{AW2}. For $r=1$, it follows by \eqref{AW1} and H\"older's inequality.
Next we assume $r\in(1,\infty]$. Let $\frac1r+\frac1{r'}=1$. By \eqref{AW1}, it suffices to prove that for any $h\in\bar\mL^{r'}$,
$$
\sI:=\int_{\mR^d}\int_{\mR^d} h(x) f(x-y) g(y)\dif x\dif y\leq \nor h\nor^*_{r'}\nor f\nor_{p}\nor g\nor^*_q.
$$
Noting that $\frac1{p'}+\frac1r=\frac1q$, by H\"older's inequality we have
\begin{align*}
\sI&=\sum_{i,j}\int_{D_i}\int_{D_j} (h(x)^{r'} f(x-y)^p)^{\frac{1}{q'}} (f(x-y)^p g(y)^q)^{\frac1{r}}(h(x)^{r'} g(y)^q)^{\frac1{p'}}\dif x\dif y\\
&\leq\sum_{i,j}\left(\int_{D_i}\int_{D_j}h(x)^{r'} f(x-y)^p\dif x\dif y\right)^{\frac1{q'}}\left(\int_{D_i}\int_{D_j} f(x-y)^p g(y)^q\dif x\dif y\right)^{\frac1{r}}\\
&\qquad\qquad\times\left(\int_{D_i}\int_{D_j} h(x)^{r'} g(y)^q\dif x\dif y\right)^{\frac1{p'}}\\
&\leq\sum_{i,j}\left(\|\1_{D_i}h\|^{\frac{r'}{q'}}_{r'}\nor f\nor_p^{\frac p{q'}}\right)\left(\nor f\nor_p^{\frac p{r}}\|\1_{D_j}g\|^{\frac q{r}}_q\right)
\left(\|\1_{D_i}h\|^{\frac{r'}{p'}}_{r'}\|\1_{D_j}g\|^{\frac q{p'}}\right)\\
&=\nor f\nor_p\sum_{i,j} \|\1_{D_i}h\| _{r'} \|\1_{D_j}g\| _q=\nor f\nor_{p}\nor h\nor^*_{r'}\nor g\nor^*_q.
\end{align*}
The proof is finished.
\end{proof}

\begin{proof}[Proof of Proposition \ref{Le21}]
(i) We only show \eqref{AW02} for $p\in[1,\infty)$ since \eqref{AW202}  is similar. Suppose that for some $C_0>0$ and any $\phi\in C_c(\mR^d)$,
\begin{align}\label{DA1}
\nor \mu^\centerdot(\phi)\nor_p\leq C_0\nor\phi\nor_p.
\end{align}
To show  it  for all $\phi\in \tilde{\mL}^p$, we divide the proof into four steps. Note that Fatou's lemma can not be used directly.
\begin{enumerate}[$\bullet$]
\item First we show \eqref{DA1} holds for any $\phi=\1_O$ with $O$ being a bounded open set. For $n\in\mN$, define
$$
\phi_n(x):=1-1/(1+\dis(x, O^c))^n.
$$
Clearly, $\phi_n\in C_c(\mR^d)$ and $\phi_n\uparrow \1_O$. Now by the monotone convergence theorem and Fatou's lemma, we have
$$
\nor \mu^\centerdot(\1_O)\nor_p=\nor \lim_{n\to\infty}\mu^\centerdot(\phi_n)\nor_p\leq \varliminf_{n\to\infty}\nor \mu^\centerdot(\phi_n)\nor_p
\leq C_0\varliminf_{n\to\infty}\nor\phi_n\nor_p\leq C_0\nor\1_O\nor_p.
$$
\item Next we show \eqref{DA1} holds for any $\phi=\1_A$ with $A$ being any Borel subset of $O=(-m,m]^d$. Define
$$
\cE:=\Big\{A\in O\cap \sB(\mR^d): \nor \mu^\centerdot(\1_A)\nor_p\leq C_0\nor\1_A\nor_p\Big\}.
$$
Let $(A_n)_{n\in\mN}\subset\cE$ and $A_n\downarrow A$. By Fatou's lemma, we have
$$
\nor \mu^\centerdot(\1_A)\nor_p=\nor \lim_{n\to\infty}\mu^\centerdot(\1_{A_n})\nor_p\leq \varliminf_{n\to\infty}\nor \mu^\centerdot(\1_{A_n})\nor_p
\leq C_0\varliminf_{n\to\infty}\nor\1_{A_n}\nor_p=C_0\nor\1_A\nor_p,
$$
where the last equality is due to $\lim_{n\to\infty}\nor\1_{A_n-A}\nor_p\lesssim \lim_{n\to\infty}\|\1_{\{A_n-A\}\cap O}\|_p=0.$ So, $A\in\cE$.
Similarly, if $A_n\uparrow A$, then $A\in\cE$. Thus $\cE$ is a monotone class. Let
$$
\cA:=\Big\{\Pi_{i=1}^d(a_i,b_i]\cap(-m,m]^d, a_i<b_i\Big\}
$$
and $\cA_\Sigma$ be the algebra generated by $\cA$ through finite disjoint unions. For given $A\in\cA_\Sigma$, there is a family of bounded open sets $A_n$ so that
$A_n\downarrow A$. By Fatou's lemma again, we have
$$
\nor \mu^\centerdot(\1_A)\nor_p=\nor \lim_{n\to\infty}\mu^\centerdot(\1_{A_n})\nor_p\leq \varliminf_{n\to\infty}\nor \mu^\centerdot(\1_{A_n})\nor_p
\leq C_0\varliminf_{n\to\infty}\nor\1_{A_n}\nor_p\leq C_0\nor\1_A\nor_p.
$$
Hence, $\cA_\Sigma\subset\cE$. By the monotone class theorem, we have
$$
\sB(O)\subset\sigma(\cA_\Sigma)\subset \cE\subset \sB(O).
$$

\item Now we show \eqref{DA1} holds for any nonnegative bounded measurable function $\phi$ with support in $O=(-m,m]^d$. By Lusin's theorem, for any $\eps>0$, there is a 
continuous function $\phi_\eps$ with support in $O$ so that
$$
\|\varphi_\eps\|_\infty\leq\|\varphi\|_\infty,\ \ \lim_{\eps\to 0}|\{x: \phi(x)\not=\phi_\eps(x)\}|=0.
$$
Let $A_\eps:=\{x: \phi(x)\not=\phi_\eps(x)\}$. By what we have proved, as $\eps\to 0$, we have 
$$
\nor \mu^\centerdot(\phi-\phi_\eps)\nor_p\leq 2\|\phi\|_\infty\nor \mu^\centerdot(\1_{A_\eps})\nor_p\leq 2\|\phi\|_\infty C_0\nor\1_{A_\eps}\nor_p\leq C\|\1_{A_\eps}\|_p\to 0.
$$
Therefore, by the dominated convergence theorem,
$$
\nor \mu^\centerdot(\phi)\nor_p= \lim_{\eps\to 0}\nor\mu^\centerdot(\phi_\eps)\nor_p 
\leq C_0\varliminf_{n\to\infty}\nor\phi_\eps\nor_p=C_0\nor\phi\nor_p.
$$
\item Finally, for general nonnegative $\phi\in\tilde\mL^p$, let $\phi_n(x):=(\phi(x)\wedge n)\1_{\{|x|<n\}}$. By the monotone convergence theorem and Fatou's lemma, we have
$$
\nor \mu^\centerdot(\phi)\nor_p=\nor \lim_{n\to\infty}\mu^\centerdot(\phi_n)\nor_p\leq \varliminf_{n\to\infty}\nor \mu^\centerdot(\phi_n)\nor_p
\leq C_0\varliminf_{n\to\infty}\nor\phi_n\nor_p\leq C_0\nor\phi\nor_p.
$$
\end{enumerate}

(ii) Let $\mX$ denote $\tilde{\mL}^p$ or $\mL^p$ and let $\mX\cP_s$ denote $\tilde{\cL}^p\cP_s$ or ${\cL}^p\cP_s$.
Suppose that $(\mu^\centerdot_n)_{n\in\mN}$ is a Cauchy sequence in $\mX\cP_s$. Since the space $\cL(\mX,\mX)$
of all bounded linear operators from $\mX$ to $\mX$ is complete  with respect to the operator norm, 
and $(\mu^\centerdot_n)_{n\in\mN}$ can be regarded as a Cauchy sequence in $\cL(\mX,\mX)$ in a natural way, there is an operator  $T\in \cL(\mX,\mX)$ such that
\begin{align}\label{ZW3}
\lim_{n\to\infty}\|\mu^\centerdot_n-{T}\|_{\cL(\mX,\mX)}=\lim_{n\to\infty}\sup_{\|\phi\|_\mX\leq 1}\|\mu^\centerdot_n(\phi)-{T}(\phi)\|_\mX=0.
\end{align}
By (i), it suffices to show that there is a sub-probability kernel $\mu^\centerdot\in\mX\cP_s$ so that for each $\phi\in C_c(\mR^d)$,
\begin{align}\label{ZR1}
{T}(\phi)(x)=\mu^x(\phi) \ \text{for Lebesgue almost all $x\in\mR^d$}.
\end{align}
Note that for each $\phi\in\mX$, there is a null set $A_\phi$ and a subsequence $n_k$ so that for each $x\notin A_\phi$,
$$
\lim_{k\to\infty}|\mu^x_{n_k}(\phi)-{T}(\phi)(x)|=0.
$$
Let $\{\phi_m\}_{m\in\mN}$ be a dense subset of $C_c(\mR^d)\subset \tilde\mL^p\cap \mL^p$. By a standard diagonalization method, one can 
find a common null set $A\subset\mR^d$  and a subsequence $n'_k$ so that for each $x\notin A$ and $m\in\mN$,
\begin{align}\label{ZR2}
\lim_{k\to\infty}|\mu^x_{n'_k}(\phi_m)-{T}(\phi_m)(x)|=0,\ \
\end{align}
Moreover, since the dual space of $C_c(\mR^d)$ is the space of all finite Borel measures,
by the Banach--Alaoglu theorem, for any $x\in\mR^d$, there is a sub-probability measures $\mu^x$ and a subsequence $n''_k(x)$ of $n'_k$ such that
\begin{align*}
\lim_{k\to\infty}\mu^x_{n''_k(x)}(\phi)=\mu^x(\phi),
\end{align*}
 which together with \eqref{ZR2} implies that for any $x\notin A$ and $m\in\mN$,
$\mu^x(\phi_m)={T}(\phi_m)(x).$
From this and by the density of $\{\phi_m,m\in\mN\}$ in $C_c(\mR^d)$, we derive \eqref{ZR1}. The completeness of $\mX\cP_s$ is obtained.

Next we show the completeness of $\tilde{\mL}^p\cP$. Let  $(\mu^\centerdot_n)_{n\in\mN}$ be a family of probability kernels. We need to  show that $\mu^x$ is a probability measure. 
For any $m\in\mN$, we define $\psi_m\in C_c(\mR^d)$ by $|\psi_m|\le 1$,
\begin{align*}
\psi_m(y)=1\ \text{for $|y|\le m$}\quad \text{and}\quad \psi_m(y)=0 \ \text{for $|y|>2m$}.
\end{align*}
It follows from \eqref{ZR1} that there is a common Lebesgue null set $A'$ such that for all $x\notin A'$,
\begin{align}\label{0227:00}
\mu^x(\psi_m)={T}(\psi_m)(x),\quad \text{for all $m\in\mN$.}
\end{align}
We note that by the fact $\sup_m\nor\psi_m\nor_{p}\lesssim\sup_m\|\psi_m\|_{\infty}=1$ and \eqref{ZW3}, for each bounded domain $D$, we also have
\begin{align}\label{ZW4}
\lim_{n\to\infty}\sup_{m}\int_D|\mu^x_n(\psi_m)-{T}(\psi_m)(x)|\dif x=0,
\end{align}
which implies that
\begin{align*}
\int_D|{T}(\psi_m)(x)|\dif x\le \lim_{n\to\infty}\sup_{m}\int_D|\mu^x_n(\psi_m)|\dif x\le |D|.
\end{align*}
Moreover, since for each $n$, 
$\lim_{m\to\infty}\int_D\mu^x_n(\psi_m)\dif x=|D|,$
and by \eqref{ZW4}, we have 
$$\lim_{m\to\infty}\int_D {T}(\psi_m)(x)\dif x=|D|.$$
This in turn implies that there is a null set $A''$ and subsequence $m_k$ so that for each $x\notin A''$,
$$
\lim_{k\to\infty}{T}(\psi_{m_k})(x)=1.
$$
This together with \eqref{0227:00} implies that for each $x\notin A'\cup A''$, $\mu^x(\mR^d)=1$.

Finally, we show the uncompleteness of ${\mL}^p\cP$ through a counterexample. Consider $d=1$ and for $n\in\mN$,
$$
\mu^x_n(\dif y)=\1_{[0,1]}(x)\1_{[0,n]}(y)\dif y/n.
$$
It is easy to see that for any $p\in[1,\infty)$,
$$
\|\mu^\centerdot_n\|_p=\|n^{-1}\1_{[0,n]}\|_{p/(p-1)}=n^{-1/p}\to 0,\ \ n\to\infty.
$$
The proof is complete.
\end{proof}

\section{Technical Lemmas}

\bl\label{LeA1}
Let $\xi\sim N(0,\sigma^2)$ be a $d$-dimensional normal random variable with mean zero and variance $\sigma^2$. 
For any $1\leq p\leq q\leq\infty$ and $j\in\mN_0$, there is a constant $C=C(j,q,p,d)>0$ such that for all $f\in\tilde\mL^p$,
\begin{align}\label{AZ09}
\nor\mE \nabla^jf(\xi+\cdot)\nor_q\lesssim_C(\sigma^{-j}+\sigma^{d/q-d/p-j})\nor f\nor_p.
\end{align}
\el
\begin{proof}
Note that by the integration by parts,
\begin{align*}
|\mE \nabla^j f(\xi+x)|&=(2\pi\sigma^2)^{-d/2}\left|\int_{\mR^d}f(y+x)\nabla^j\e^{-|y|^2/2\sigma^2}\dif y\right|\\
&\leq(2\pi\sigma^2)^{-d/2}\int_{\mR^d}|f(y+x)|\,|\nabla^j\e^{-|y|^2/2\sigma^2}|\dif y\\
&\lesssim\sigma^{-d-j}\int_{\mR^d}|f(y+x)|\e^{-c|y|^2/\sigma^2}\dif y=\sigma^{-d-j}|f|*\phi_\sigma(x),
\end{align*}
where we have used that for some $c>0$,
$$
|\nabla^j\e^{-|y|^2/2\sigma^2}|\lesssim\sigma^{-j}\e^{-c|y|^2/\sigma^2}=:\phi_\sigma(y).
$$ 
Let $1+\frac1q=\frac1p+\frac1r$. By Young's convolution inequality \eqref{AW2}, we get
$$
\nor\mE \nabla^jf(\xi+\cdot)\nor_q\lesssim \sigma^{-d-j}\nor f\nor_p \nor\phi_\sigma\nor^*_r.
$$
By the definition of $\nor\cdot\nor^*_r$, we have
\begin{align*}
\nor\phi_\sigma\nor^*_r=\sum_i\left(\int_{D_i}\e^{-cr|x|^2/\sigma^2}\dif x\right)^{1/r}\lesssim\int_{\mR^d}\left(\int_{D_z}\e^{-cr|x|^2/\sigma^2}\dif x\right)^{1/r}\dif z,
\end{align*}
where $D_z$ is the unit cube with center $z\in\mR^d$.
Noting that for $|z|\geq \sqrt d$ and $x\in D_z$,
$$
|x|\geq|z|-|x-z|\geq|z|-\sqrt d/2\geq |z|/2,
$$
we have
$$
\int_{|z|\geq \sqrt d}\left(\int_{D_z}\e^{-cr|x|^2/\sigma^2}\dif x\right)^{1/r}\dif z
\leq \int_{\mR^d}\e^{-c|z|^2/4\sigma^2}\dif z\lesssim \sigma^d.
$$
On the other hand, we clearly have
$$
\int_{|z|\leq \sqrt d}\left(\int_{D_z}\e^{-cr|x|^2/\sigma^2}\dif x\right)^{1/r}\dif z
\lesssim\left(\int_{\mR^d}\e^{-cr|x|^2/\sigma^2}\dif x\right)^{1/r}\lesssim \sigma^{d/r}.
$$
Hence,
$$
\nor\mE \nabla^jf(\xi+\cdot)\nor_q\lesssim \sigma^{-d-j}(\sigma^d+\sigma^{d/r})\nor f\nor_p,
$$
which in turn gives the desired estimate.
\end{proof}

 Recall
\begin{align*}
\cB(\alpha,\beta):=\int_0^1 r^{\alpha-1}(1-r)^{\beta-1}\dif r,\quad  \text{for $\alpha,\beta\ge0$}.
\end{align*}
\bl[Estimate for Beta functions]\label{Beta}
For any $\alpha\in(0,1]$ and $\beta>0$, we have for any $k\in\mN$,
\begin{align*}
\cB(\alpha, k\beta+1)\le \left(\frac1\alpha+\frac1\beta\right)k^{-\alpha}.
\end{align*}

\el
\begin{proof}
For any $h\in(0,1)$, one sees that
\begin{align*}
\cB(\alpha,k\beta+1)&\le \int_0^h r^{\alpha-1}\dif r+h^{\alpha-1}\int_h^1(1-r)^{k\beta}\dif r
\le \frac1\alpha h^\alpha+\frac1{k\beta+1}h^{\alpha-1}\le \left(\frac1\alpha+\frac1{k\beta}h^{-1}\right)h^\alpha.
\end{align*}
Taking $h=k^{-1}$, we complete the proof.
\end{proof}

\bl\label{Seri}
For any $\alpha>0$, there is a constant $C=C(\alpha)>0$ such that for all $\lambda\geq 1$,
\begin{align*}
\sum_{m=0}^\infty\frac{\lambda^m}{(m!)^\alpha}\leq \e^{C\lambda^{1/\alpha}\ln\lambda}.
\end{align*}
\el
\begin{proof}
%Without loss of generality, we assume $\lambda\geq 1$.
By Stirling's formula, we have
\begin{align*}
\sum_{m=0}^\infty\frac{\lambda^m}{(m!)^\alpha}
\leq 1+C\sum_{m=1}^\infty\frac{\lambda^m}{m^{m\alpha}}
&\leq 1+C\int^\infty_1\frac{\lambda^x}{x^{x\alpha}}\dif x
=1+C\int^\infty_1\e^{\alpha x\ln(\lambda^{1/\alpha}/x)}\dif x\\
&\leq 1+C\int^{2\lambda^{1/\alpha}}_1\e^{x\ln\lambda}\dif x+C\int^\infty_{2\lambda^{1/\alpha}}\e^{-\alpha x\ln 2}\dif x.
\end{align*}
From this we derive the desired estimate.
\end{proof}

\bl[Gronwall's inequality]\label{lemA03}
Let $f(s,t), g(s,t):\mD_T\to[0,\infty)$ and $h:[0,T]\to[0,\infty)$ be measurable functions. Assume that for all $(s,t)\in\mD_T$,
$$
f(s,t)\leq g(s,t)+\int^t_sh(r)(f(s,r)+f(r,T))\dif r.
$$
Then we have
$$
f(s,t)\leq G(s,t)+\int^t_sH(s',t)\left(G(s',T)+\int^T_{s'} G(r,T)H(r)\e^{\int^r_{s'} H(r')\dif r'}\dif r\right)\dif s',
$$
where 
$$
G(s,t):=g(s,t)+\int^t_s g(s,r)h(r)\e^{\int^t_r h(r')\dif r'}\dif r
$$
and
$$
H(r,t):=h(r)\left(1+\int^t_r h(r')\e^{\int^t_{r'} h(r'')\dif r''}\dif r'\right).
$$
\el
\begin{proof}
For fixed $s\in[0,T]$, by the assumption we have
$$
f(s,t)\leq F(s,t)+\int^t_sh(r)f(s,r)\dif r,
$$
where
$$
F(s,t):=g(s,t)+\int^t_sh(r)f(r,T)\dif r.
$$
By the usual Gronwall's inequality we get
\begin{align*}
f(s,t)&\leq F(s,t)+\int^t_s F(s,r)h(r)\e^{\int^t_rh(r')\dif r'}\dif r\\
&=g(s,t)+\int^t_sh(r)f(r,T)\dif r
+\int^t_s g(s,r)h(r)\e^{\int^t_rh(r')\dif r'}\dif r\\
&\quad+\int^t_s\left( \int^{r}_sh(r')f(r',T)\dif r'\right)h(r)\e^{\int^t_rh(r')\dif r'}\dif r\\
&=G(s,t)+\int^t_sH(r,t)f(r,T)\dif r,
\end{align*}
where $G$ and $H$ are defined in the lemma.
In particular,
$$
f(s,T)\leq G(s,T)+\int^T_s H(r,T)f(r,T)\dif r,
$$
By Gronwall's inequality again, we have
$$
f(s,T)\leq G(s,T)+\int^T_s G(r,T)H(r)\e^{\int^r_s H(r')\dif r'}\dif r,
$$
Combining the above calculations, we obtain the desired estimate.
\end{proof}

\bl\label{lem:0301}
Let $E_i,$ be two Banach spaces with norms $\|\cdot\|_i$, $i=1,2$. Let $G: E_1\to E_2$ be a Lipschitz mapping with $G(0)=0$ and define
\begin{align*}
F(x):=G(x)\1_{\|x\|_1\le1}+G(x)\1_{\|x\|_1>1}/|x|_1.
\end{align*}
Then for any $x,y\in E$
\begin{align*}
\|F(x)-F(y)\|_2\le 2\|G\|_{\rm Lip}\|x-y\|_1.
\end{align*}
\el
\begin{proof}
We consider three cases: (i) $\|x\|_1\wedge \|y\|_1\le 1$; (ii) $\|x\|_1\leq1<\|y\|_1$; (iii) $\|x\|_1\wedge \|y\|_1>1$.
 In case (i), $F(x)=G(x)$, it is trivial.
 In case (ii), one sees that
 \begin{align*}
\|F(y)-F(x)\|&=\left\|G(y)/\|y\|_1-G(x)\right\|_2\le\left\|(G(y)-G(x))/\|y\|_1\right\|_2+\|G(x)/\|y\|_1-G(x)\|_2\\
&\le\left(\|G\|_{\rm Lip}\|x-y\|_1+\|G(x)\|_2(\|y\|_1-1)\right)\\
&\le \|G\|_{\rm Lip}\left(\|x-y\|_1+\|x\|_1(\|y\|_1-\|x\|_1)\right)\le 2\|G\|_{\rm Lip}\|x-y\|_1.
\end{align*}
 In case (iii), we have
 \begin{align*}
\|F(y)-F(x)\|&=\left\|\frac{G(x)\|y\|_1-G(y)\|x\|_1}{\|x\|_1\|y\|_1}\right\|_2\le \frac{\|(G(x)-G(y))\|y\|_1-G(y)(\|x\|_1-\|y\|_1)\|_2}{\|x\|_1\|y\|_1}\\
&\le\|G\|_{\rm Lip}\|x-y\|_1+\|G(y)\|_2\|x-y\|_1/\|y\|_1\leq 2\|G\|_{\rm Lip}\|x-y\|_1.
\end{align*}
The proof is complete.
\end{proof}

Consider the following PDE:
\begin{align}\label{PDE0}
\p_t u=\Delta u+b\cdot \nabla u, \quad u_0=\phi.
\end{align}

\bt\label{ThA}
Let $q_1,p_1,p_0\in(1,\infty]$ satisfy $\tfrac2{q_1}+\tfrac d{p_1}<1.$
Assume $b\in\mL^{q_1}_T\tilde\mL^{p_1}$ and $\phi\in C^\infty_b(\mR^d)$. 
For any $p\in[p_0\vee \frac{p_1}{p_1-1},\infty]$ with $\frac2{q_1}+\frac d{p_0}<1+\frac dp$, 
there is a unique solution $u$ to PDE \eqref{PDE0} with 
$$
\nor\nabla^j u(t)\nor_p\lesssim t^{-\frac{j+d/p_0-d/p}{2}}\nor\phi\nor_{p_0}.
$$
\et
\begin{proof}
Let $(P_t)_{t\geq 0}$ be the Gaussian heat semigroup. By Duhamel's formula, we have
\begin{align*}
u(t)=P_t\phi+\int_0^t P_{t-s}(b\cdot \nabla u)(s)\dif s.
\end{align*}
Let $p\in[p_0\vee p_3,\infty]$ satisfy $\tfrac1{p_3}=\tfrac1{p_1}+\tfrac1{p}\leq 1.$ 
For $j=0,1$, by Lemma \ref{LeA1}, we have
\begin{align}
\nor\nabla^j u(t)\nor_p&\lesssim t^{-\frac{j+d/p_0-d/p}{2}}\nor\phi\nor_{p_0}
+\int_0^t (t-s)^{-\frac{j+d/p_3-d/p}{2}}\nor b(s)\cdot\nabla u(s)\nor_{p_3}\dif s\no\\
&\lesssim t^{-\frac{j+d/p_0-d/p}{2}}\nor\phi\nor_{p_0}
+\int_0^t (t-s)^{-\frac{j+d/p_1}{2}}\nor b(s)\nor_{p_1}\nor\nabla u(s)\nor_{p}\dif s.\label{AZ3}
\end{align}
Suppose $\tfrac 2{q_1}+\tfrac d{p_0}<1+\tfrac dp.$
By H\"older's inequality, we have
\begin{align*}
\nor\nabla u(t)\nor^{q'_1}_p
&\lesssim t^{-q'_1\frac{1+d/p_0-d/p}{2}}\nor\phi\nor_{p_0}
+\nor b\nor^{q'_1}_{\mL^{q_1}_T\tilde\mL^{p_1}}\int_0^t (t-s)^{-q_1'\frac{1+d/p_1}{2}}\nor\nabla u(s)\nor^{q_1'}_{p}\dif s,
\end{align*}
which implies that by Gronwall's inequality of Volterra's type,
$$
\nor\nabla u(t)\nor_p\lesssim t^{-\frac{1+d/p_0-d/p}{2}}\nor\phi\nor_{p_0}.
$$
The proof is complete.
\end{proof}

\end{appendix}

\end{document}